\documentclass[10pt]{amsart}

\usepackage{color,graphicx,amssymb}
\usepackage{amsthm}
\usepackage{dsfont}

\theoremstyle{definition}

\newtheorem{theorem}{Theorem}[section]
\newtheorem{definition}[theorem]{Definition}
\newtheorem{lemma}[theorem]{Lemma}
\newtheorem{corollary}[theorem]{Corollary}
\newtheorem{proposition}[theorem]{Proposition}
\newtheorem{claim}[theorem]{Claim}
\newtheorem{example}[theorem]{Example}

\newtheorem{remark}[theorem]{Remark}

\newtheorem{acknowledgments}[theorem]{Acknowledgments}

\DeclareMathOperator{\Con}{Con}

\DeclareMathOperator{\Hom}{Hom}
\DeclareMathOperator{\Cg}{Cg}

\DeclareMathOperator{\id}{id}

\DeclareMathOperator{\im}{im}
\DeclareMathOperator{\ar}{ar}

\begin{document}

\title{Extensions realizing affine datum : central extensions}
\author{Alexander Wires}
\address{School of Mathematical Sciences, University of Electronic Science and Technology\\Chengdu 611731, Sichuan, China}
\email{awires81@uestc.edu.cn}
\date{October 1, 2023}


\begin{abstract}
The study of extensions realizing affine datum is specialized to central extensions in varieties with a difference term which leads to generalizations of several classical theorems on central extensions from group theory. We establish a 1-dimensional Hochschild-Serre sequence for a central extension equipped with affine datum. This is used to develop a Schur-Hopf formula which characterizes the $2^{\mathrm{nd}}$-cohomology group of regular datum in terms of the transgression map and commutators in free presentations. We prove, assuming the existence of an idempotent, the existence of covers and provide a cohomological characterization of perfect algebras. The class of varieties with a difference term contain all varieties of algebras with modular congruence lattices; for example, any variety of groups with multiple operators in the parlance of P.J. Higgins or algebras of Loday-type - analogous results recently established for these algebras can be recovered by specialization.
\end{abstract}

\maketitle


\section{Introduction}\label{sec:1}
\vspace{0.3cm}

We continue the study of extensions realizing affine datum in varieties of universal algebras which was inaugurated in Wires \cite{wiresI}. In the present manuscript, we focus on central extensions in varieties with a difference term and broaden to this more general class of algebras several classical results: a 1-dimensional Hochschild-Serre (inflation/restriction) exact sequence for central extensions, characterizations of injectivity and surjectivity of the transgression map, a Schur-formula relating commutators in free presentations to $2^{\mathrm{nd}}$-cohomology of regular datum, the existence of covers and a cohomological characterization of perfect algebras.

A variety $\mathcal V$ has a difference term if there is a ternary term $m$ in the signature such that for all algebras $A \in \mathcal V$,
\begin{align}\label{eqn:diffterm}
m(x,x,y) &= y  \quad \quad \text{and} \quad \quad x \mathrel{[\alpha,\alpha]}  m(x,y,y)    &(x,y \in \alpha \in \Con A).
\end{align}
The commutator referenced here is the term-condition commutator (or TC-commutator). The class of varieties with a difference term is a Mal'cev condition (see Kearnes \cite[Thm 3.3]{vardiff}) which properly contains varieties of algebras which have modular congruence lattices and so provides a common and vast generalization for many classes of algebras like groups, quasigroups, braces and conformal algebras (see Smith \cite{smith}), but also modules expanded by multilinear operations (or $\Omega$-algebras in the language of Kurosh \cite{kurosh}) such as rings, Leibniz algebras, diassociative algebras and Rota-Baxter algebras to give just a few examples. It is also a class for which the TC-commutator has good formal properties such as commutativity, join-additivity and the homomorphism property \cite[Thm 2.10]{vardiff}. From the defining property Eq~\eqref{eqn:diffterm}, it follows that the difference term interprets as a Mal'cev operation in any abelian algebra; therefore, the abelian algebras form a coslice variety term-equivalent to the variety of unital modules over a fixed ring with unity (see Herrmann \cite{herrmann}, Freese and McKenzie \cite[Ch. 9]{commod} and Bergman \cite[Thm 7.38]{bergman}).

Recent work has examined the notion of multipliers, liftings of central extensions, and covers of perfect algebras and established analogues of the 1-dimensional Hochschild-Serre sequences and Schur-Hopf formulas for central extensions in some classes of algebras of Loday-type such as Lie algebras in Batten \cite{batten1}, Leibniz algebras in Elyse \cite{elyse1} and Mainellis \cite{mainellis1}, and diassociative algebras in Mainellis \cite{mainellis2,mainellis3} to give just a very few examples. We confirm these as special cases of the more general study of central extensions in varieties with a difference term; in each particular class of algebras, the previous result can be recovered by specialization of the terms involved.

In varieties with a difference term, the general representation of central extensions realizing affine datum \cite[Def 3.12]{wiresI} takes an alternative form \cite[Cor 2.7]{wiresI} which is more closely analogous to the representation of central extensions over a product universe used in group theory. Instead of redeveloping the apparatus of second-cohomology, we use Section~\ref{sec:3} to give a self-contained review of the machinery for extensions realizing affine datum in arbitrary varieties and then derive the special form terms like 2-cocycles, 2-coboundaries and stabilizing isomorphisms take in light of the representation \cite[Cor 2.7]{wiresI} in varieties with a difference term. In Section~\ref{sec:4}, we prove a 1-dimensional Hochschild-Serre sequence for central extensions in varieties with a difference term. We then examine more closely the action of the inflation map on second-cohomology in Lemma~\ref{lem:30} and Lemma~\ref{lem:32}. In Section~\ref{sec:5}, we define the Schur multiplier (Definition~\ref{def:schurmult}) of an algebra in any variety with a difference term as a a certain subalgebra of a kernel algebra derived from commutators in a free presentation. We then establish Theorem~\ref{thm:schurmult} which is a generalization of the Schur-Hopf formula which relates the second-cohomology of regular datum to the homeset of the Schur multiplier and the regular algebra. We remark how to recover the classic theorem of Schur \cite{schur2} in the case of groups from the general theorem, but also consider examples of varieties which are modules expanded by multilinear operations. We continue by characterizing injectivity and surjectivity of the transgression map related to second-cohomology of regular datum. The existence of covers is established in Theorem~\ref{thm:coversexist}. In Section~\ref{sec:6}, we consider universal central extensions and perfect algebras. Following Milnor \cite{milnor} in the case of groups, we give in Theorem~\ref{thm:perfectchar} a characterization of perfect algebras and universal central extensions by second-cohomology; consequently, this yields a universal central extension from the cover of a perfect algebra constructed from any free presentation. We reserve to Section~\ref{sec:7} for the proof that the Schur multiplier of an algebra is invariant under free-presentations.

Anyone familiar with central extensions and Schur multipliers of groups will observe how the outline of the present manuscript is indebted to the monograph of Karpilovski \cite{karpil} for the presentation of those topics.

\vspace{0.3cm}

\section{Preliminaries}\label{sec:2}
\vspace{0.3cm}

Consult Bergman \cite{bergman} for fundamentals of universal algebra; however, for this manuscript we deviate with the use of boldface type to denote the algebra structure as opposed to the universe of the algebra which is in plain type. We use plain type for both the universe of the algebra and algebra structure itself in the hope that context will clarify any potential confusion. We make constant use of basic properties of the term-condition commutator in varieties of universal algebras. For the basic definition and general properties of the term-condition commutator, Kearnes and Kiss \cite[Cha 2.5]{shape} is useful and for a thorough development of the theory of the commutator in congruence modular varieties, consult Freese and McKenzie \cite{commod} and McKenzie and Snow \cite{mckenziesnow}. The properties of the term-condition commutator in varieties with a difference term is developed in Kearnes \cite{vardiff}; in particular, familiarity with \cite[Thm 2.10]{vardiff} which compiles the basic properties developed in that manuscript is recommended. For the most part, any special definitions are introduced in the text in the place they are needed.

It will be beneficial to be explicit about the connection between abelian Mal'cev varieties and varieties of modules. The constructions involved will be referenced by Proposition~\ref{prop:regexists}. Let $\mathcal A$ be an abelian Mal'cev variety in the signature $\tau$ with Mal'cev term $m(x,y,z)$. Consider the set of terms $R = \left\{ r(x,y) \in F_{\mathcal A}(x,y) : r(y,y) = y \right\}$ in the 2-generated free algebra in $\mathcal A$. Then $R=\left\langle R, +, - , \cdot, y,x \right\rangle$ is a unital ring under the definitions
\begin{align*}
r_{1} + r_{2} &:= m(r_{1},y,r_{2}) , &-r := m(y,r,y) &\quad \quad \text{and} &r_{1}(x,y) \cdot r_{2}(x,y) := r_{1}(r_{2}(x,y),y)
\end{align*}
where $y$ is the zero element for addition and $x$ is the identity for multiplication. Every abelian algebra in $A \in \mathcal A$ is polynomially equivalent to an $R$-module effected in the following manner: for any choice $0 \in A$, then the abelian group operations are given by
\begin{align*}
a + b &:= m^{A}(a,0,b) &-a := m(0,a,0)
\end{align*} 
with zero element $0$ and the action of the ring $R$ is given by 
\begin{align*}
r(x,y) \cdot a := r^{A}(a,0) 
\end{align*}
for $r(x,y) \in R$, $a \in A$. If $J=F_{\mathcal A}(x)$ is the 1-generated free algebra in $\mathcal A$, then we construct the $R$-module $\tilde{J}$ were we make the choice of the free generator $x$ to be the zero of the module.

The abelian variety $\mathcal A$ and the variety $_R \mathcal Mod$ of $R$-modules are very closely related if the choice of the zero element in the construction of the associated modules is folded into the correspondence - this is accomplished by passing to the coslice varieties $(J \downarrow \mathcal A)$ and $(\tilde{J} \downarrow \,  _R \mathcal Mod)$. Given an abelian algebra $A \in \mathcal A$ and homomorphism $\mu: J \rightarrow A$, $M(A,\mu)$ is defined as the $R$-module over the universe of $A$ where the abelian group operation is given by
\begin{align*}
a + b &:= m^{A}(a,\mu(x),b) &-a := m(\mu(x),a,\mu(x))
\end{align*} 
with zero element $\mu(x)$ and the action of the ring $R$ is given by 
\[
r(x,y) \cdot a := r^{A}(a,\mu(x)) 
\]
for $r(x,y) \in R$, $a \in M$. Given $f \in \tau$, define the terms
\begin{align*}
r_{f,i}(x,y) &:= m(f(y,\ldots,y,x,y,\ldots,y),f(y,\ldots,y),y) \\
d_{f} &:= f^{J}(x,\ldots,x)
\end{align*}
for $1 \leq i \leq \ar f$. We see that $r_{f,i} \in R$. Given an $R$-module $M$ and homomorphism $\mu :\tilde{J} \rightarrow M$, then $A(M,\mu)$ is an algebra with fundamental operations defined by
\begin{align}\label{eqn:idemmodule}
f^{A(M,\mu)}(x_{1},\ldots,x_{n}) := \sum_{i=1}^{n} r_{f,i} \cdot x_{i} + \mu(d_{f})
\end{align}
for $f \in \tau$ with $n = \ar f$. With these interpretations, we see that $\tilde{J} = M(J,\id)$ and the correspondence 
\begin{align*}
(A,\mu) &\longmapsto (M(A,\mu),\mu) \\
(M,\mu) &\longmapsto (A(M,\mu),\mu)
\end{align*}
witnesses the term-equivalence of the coslice varieties $(J \downarrow \mathcal A)$ and $(\tilde{J} \downarrow \,  _R \mathcal Mod)$.

\begin{theorem}\cite[Thm 7.38]{bergman}
Let $\mathcal A$ be a nontrivial abelian Mal'cev variety. There is ring $R$ and an algebra $J \in \mathcal A$ and an $R$-module $\tilde{J}$ such that the coslice varieties $(J \downarrow \mathcal A)$ and $(\tilde{J} \downarrow \,  _R \mathcal Mod)$ are term-equivalent.
\end{theorem}

If we let $\hat{0}: \tilde{J} \rightarrow M$ be the homomorphism which maps every element to the zero of the module, then it follows from Eq~\eqref{eqn:idemmodule} that $\left\{ A(M,\hat{0}) : (M, \hat{0}) \in (\tilde{J} \downarrow \,  _R \mathcal Mod) \right\}$ denotes the class of abelian algebras which have an idempotent element.

\vspace{0.3cm}

\section{Central extensions in varieties with a difference term}\label{sec:3}
\vspace{0.3cm}

We refer to \cite{wiresI} for the deconstruction/reconstruction (or cohomology) of extensions realizing affine datum. When restricted to central extensions in varieties with a difference term, the machinery can take an equivalent modified form which is more useful for the ensuing development. Let us first review the definitions for extensions realizing affine datum and then derive the modified notions for central extensions in varieties with a difference term.

Some notation and constructions. As was seen in the previous section, if we fix a ternary operation $m: A^{3} \rightarrow A$ on the set $A$, then any choice $u \in A$ determines a binary operation by $x +_{u} y : = m(x,u,y)$. For $0 < n \in \mathds{N}$, write $[n] = \{1,\ldots,n\}$ for the initial segment of positive integers and $[n]^{\ast} = \{ s \subseteq \{ 1,\ldots,n \}: 0 < |s| < n \}$ for the non-empty proper subsets of $[n]$. Given sets $Q$ and $A$ and $s \in [n]^{\ast}$, we define the set of coordinates $[Q,A]^{s} = \{ \vec{x} \in Q \cup A : \vec{x}(i) \in A \text{ for } i \in s, \vec{x}(i) \in Q \text{ for } i \not\in s \}$. Let $A$ be an algebra with $\alpha,\beta \in \Con A$.  
\begin{itemize}

	\item $A(\alpha) = \left\{ (x,y) \in A \times A : (x,y) \in \alpha \right\}$ is the congruence $\alpha$ as a subalgebra of $A \times A$.

	\item $\Delta_{\alpha \beta} = \Cg^{A(\alpha)}\left( \left\{ \left\langle \begin{bmatrix} u \\ u \end{bmatrix} , \begin{bmatrix} v \\ v \end{bmatrix} \right\rangle : (u,v) \in \beta \right\} \right) \in \Con A(\alpha)$.
	
	\item The diagonal homomorphism is the map $\delta: A \rightarrow A(\alpha)/\Delta_{\alpha \alpha}$ given by $\delta(x) = \begin{bmatrix} x \\ x \end{bmatrix}/\Delta_{\alpha \alpha}$.
	
	\item Given an equivalence relation $\alpha$ on a set $A$, there is an equivalence relation $\hat{\alpha}$ on the set $A(\alpha)$ defined by $\begin{bmatrix} x \\ y \end{bmatrix} \mathrel{\hat{\alpha}} \begin{bmatrix} u \\ v \end{bmatrix} \Leftrightarrow x \mathrel{\alpha} y \mathrel{\alpha} u \mathrel{\alpha} v$.
	
	\item If $\pi : A \rightarrow A/\alpha$ is the canonical epimorphism, then there is an induced extension $\rho: A(\alpha)/\Delta_{\alpha \alpha} \rightarrow A/\alpha$ given by $\rho \left( \begin{bmatrix} x \\ y \end{bmatrix}/\Delta_{\alpha \alpha} \right) = \pi(x)$. We note that $\hat{\alpha}/\Delta_{\alpha \alpha} = \ker \rho$.

\end{itemize}

For a surjective map $\pi : A \rightarrow Q$ with $\alpha = \ker \pi$, a \emph{lifting} of $\pi$ is a map $l: Q \rightarrow A$ such that $\pi \circ l = \id_{Q}$. Note that a lifting $l: Q \rightarrow A$ for $\pi$ has the property that $(x,l(q)) \in \alpha \Leftrightarrow \rho\left(\begin{bmatrix} l(q) \\ x \end{bmatrix} \right)=q$. We see that $l$ is no longer a right-inverse for $\rho$ as a set map, but that $\rho \circ (\delta \circ l) = \id$; therefore, while not technically correct, we still refer to $l$ as a lifting associated to $\rho$. The $\alpha$-trace $r: A \rightarrow A$ associated to the lifting is defined by $r = l \circ \rho$.

We will be stating results for algebras whose universe are given by sets of the form $A(\alpha)/\Delta_{\alpha \alpha}$. Given a ternary operation $m: A^{3} \rightarrow A$ and $u \in A$, we use the same convenient notation $x +_{u} y : = m(x,\delta(u),y)$ for the derived binary operation on $A(\alpha)/\Delta_{\alpha \alpha}$ determined by the $\Delta_{\alpha \alpha}$-class of the diagonal $\delta(u) = \begin{bmatrix} u \\ u \end{bmatrix}/\Delta_{\alpha \alpha}$. This is the context in which this notation will most often appear. Let us define the structure which encodes the information provided by the kernel of a possible extension.

\begin{definition}
Fix signature $\tau$ and ternary operation $m$. Define $A^{\alpha,\tau} = \left\langle  A,\alpha, \{f^{\Delta} : f \in \tau\} \right\rangle$ where
\begin{itemize}

	\item $A = \left\langle A, m \right\rangle$ is an algebra in the single operation symbol $m$;
	
	\item $\alpha \in \Con A$;
	
	\item $\{f^{\Delta} : f \in \tau\}$ is a sequence of operations $f^{\Delta} : A(\alpha)/\Delta_{\alpha \alpha} \times \delta(A)^{\ar f - 1} \rightarrow A(\alpha)/\Delta_{\alpha \alpha}$ where $\delta: A \rightarrow  A(\alpha)/\Delta_{\alpha \alpha}$ is the diagonal map.

\end{itemize}
\end{definition}

Note $A(\alpha)/\Delta_{\alpha \alpha}$ is constructed from the algebra $\left\langle A, m\right\rangle$ in the single ternary operation $m$. The partial structure $A^{\alpha,\tau}$ is \emph{homomorphic} if for each $f \in \tau$ and $a \mathrel{\hat{\alpha}/\Delta_{\alpha \alpha}} \delta(u) \mathrel{\hat{\alpha}/\Delta_{\alpha \alpha}} b$, $\vec{x} \in A^{\ar f - 1}$ we have 
\[
f^{\Delta} \left( a +_{u} b, \delta(\vec{x}) \right) = f^{\Delta} \left( a , \delta(\vec{x}) \right) +_{v} f^{\Delta} \left( b, \delta(\vec{x}) \right)
\]
where $v=l \circ \rho \left( f^{\Delta}(\delta(u), \delta(\vec{x}) \right)$ for any lifting $l: A/\alpha \rightarrow A$ associated to $\rho : A(\alpha) \rightarrow Q$.

\begin{definition}\label{def:pairing}
Let $Q,A$ be sets and $f$ an n-ary operation on $Q$. A \emph{pairing} of $Q$ on $A$ with respect to $f$ is a choice of subsets $\sigma(f) \subseteq [n]^{\ast}$ and a sequence of functions $a(f,s) : [Q,A]^{s} \rightarrow A$ with $s \in \sigma(f)$.

Let $Q$ be an algebra in the signature $\tau$. For a set $A$, an \emph{action} $Q \ast A$ is a sequence of pairings $\{a(f,s) : f \in \tau, s \in \sigma(f) \subseteq [\ar f]^{\ast} \}$.
\end{definition}

For affine extensions, all the actions considered are said to be \emph{unary}; that is, $\sigma(f) = \left\{ \{i\}: i \in [\ar f] \right\}$ is the set of all singleton subsets for the coordinates of $f$. The pairings of the action will then be functions of the form $a(f,i): Q^{i - 1} \times A \times Q^{\ar f - i} \rightarrow A$. In the following, the actions will all be unary.

An action $Q \ast A(\alpha)/\Delta_{\alpha \alpha}$ is \emph{homomorphic} if for all $f \in \tau$ with $\ar f \geq 2$ and $a \mathrel{\hat{\alpha}/\Delta_{\alpha \alpha}} \delta(u) \mathrel{\hat{\alpha}/\Delta_{\alpha \alpha}} b$, $\vec{q} \in Q^{\ar f}$ we have 
\[
a(f,i)(q_{1},\ldots,a +_{u} b,\ldots,q_{n}) = a(f,i)(q_{1},\ldots,a,\ldots,q_{n}) +_{v} a(f,i)(q_{1},\ldots,b,\ldots,q_{n})
\]
where $v=l(a(f,i)(q_{1},\ldots,\delta(u),\ldots,q_{n}))$ for any lifting $l:Q \rightarrow A$ associated to $\rho$.

\begin{definition}\label{def:datum}
Fix a signature $\tau$. A triple $\left( Q, A^{\alpha,\tau},\ast \right)$ is \emph{datum} if the following holds:
\begin{enumerate}
		
		\item[(D1)] $Q$ is an algebra in the signature $\tau$;
		
		\item[(D2)] The ternary operation symbol $m$ referenced in $A^{\alpha,\tau}$ also has an interpretation in $Q$ which is idempotent and there is a surjective homomorphism $\rho: A(\alpha) \rightarrow \left\langle Q, m \right\rangle$ such that $\hat{\alpha} = \ker \rho$;
		
		\item[(D3)] The action $Q \ast A(\alpha)/\Delta_{\alpha \alpha}$ and $A^{\alpha,\tau}$ have the following property: for any $f \in \tau$ with $n=\ar f \geq 2$ and $\rho \left(\begin{bmatrix} y_i \\ x_i \end{bmatrix} \right) =q_i$ for $i=1,\ldots,n$, we have
		\[
		f^{\Delta}\left(\begin{bmatrix} y_1 \\ x_1 \end{bmatrix}/\Delta_{\alpha \alpha},\delta(x_2),\ldots,\delta(x_n) \right) \mathrel{\hat{\alpha}/\Delta_{\alpha \alpha}} a(f,s)\left( \vec{z} \right)
		\]
for all $s \in \sigma(f)$ and $\vec{z} \in [Q, A(\alpha)/\Delta_{\alpha \alpha}]^{s}$ such that $z_{i} = q_{i}$ for $i \not\in s$ and $z_{i} = \begin{bmatrix} y_i \\ x_i \end{bmatrix}/\Delta_{\alpha \alpha}$ for $i \in s$. If $1 \in \sigma(f)$ for $f \in \tau$ with $\ar f \geq 2$, then 
		\[
		f^{\Delta}\left(\begin{bmatrix} y_1 \\ x_1 \end{bmatrix}/\Delta_{\alpha \alpha},\delta(x_2),\ldots,\delta(x_n) \right) = a(f,1) \left( \begin{bmatrix} y_1 \\ x_1 \end{bmatrix}/\Delta_{\alpha \alpha}, q_{2},\ldots,q_{n} \right) .
		\]

\end{enumerate}
We say $(Q,A^{\alpha,\tau},\ast)$ is \emph{affine datum} if in addition
\begin{enumerate}

	\item[(AD1)] $\alpha \in \Con \left\langle A, m \right\rangle$ is an abelian congruence and $m$ is a ternary abelian group operation when restricted on each $\alpha$-class;
	
	\item[(AD2)] the action $Q \ast A(\alpha)/\Delta_{\alpha \alpha}$ is unary and both the action and $A^{\alpha,\tau}$ are homomorphic. 

\end{enumerate}
\end{definition}

A \emph{2-cocycle} appropriate for affine datum $(Q,A^{\alpha,\tau},\ast)$ is a sequence of functions $T= \{ T_{f}: f \in \tau \}$ indexed by the signature such that each $T_{f}: Q^{\ar f} \rightarrow A(\alpha)/\Delta_{\alpha \alpha}$. If we fix a lifting $l:Q \rightarrow A$ associated to the datum, then we can define an algebra $A_{T}(Q,A^{\alpha,\tau},\ast)$ on the universe $A(\alpha)/\Delta_{\alpha \alpha}$ of the datum with new operations
\begin{align}\label{eqn:algebrarep}
\begin{split}
F_{f} \left(\begin{bmatrix} a_1 \\ b_1 \end{bmatrix}/\Delta_{\alpha \alpha}, \ldots, \begin{bmatrix} a_n \\ b_n \end{bmatrix}/\Delta_{\alpha \alpha} \right) &:= 
 \sum_{i=1}^{n} a(f,i)\left(q_1,\ldots,q_{i-1}, \begin{bmatrix} a_{i} \\ b_{i} \end{bmatrix}/\Delta_{\alpha \alpha}, q_{i+1},\ldots,q_n \right) \\
&\quad +_{u} \ T_{f}(q_1,\ldots,q_n) \quad \quad \quad \quad \quad \quad \quad \quad \quad \quad \ (f \in \tau).
\end{split}
\end{align}
a left-associated composition where $u= l(f^{Q}(q_1,\ldots,q_n))$ and $\rho \left(\begin{bmatrix} a_i \\ b_i \end{bmatrix}/\Delta_{\alpha \alpha} \right)=q_i$. Note we are using the identification $f^{\Delta} = a(f,1)$ from (D3).

A plentiful source of extensions which determine affine datum are given by abelian congruences in varieties with a weak difference term \cite[Thm 3.19]{wiresI}. Let $\alpha \in \Con A$ be an abelian congruence and $A \in \mathcal V$ a variety with a weak-difference term $m$. Set $Q = A/\alpha$ and choose a lifting $l: Q \rightarrow A$ for the canonical epimorphism $\pi: A \rightarrow Q$. The following definitions 
\begin{align}
\begin{split}
a(f,i) \left(x_{1}, \ldots, a_{i}, \ldots,x_{n} \right) &:= f^{A(\alpha)/\Delta_{\alpha \alpha}}\left( \delta \circ l(x_{1}),\ldots, a_{i}, \ldots, \delta \circ l(x_{n}) \right) \\
f^{\Delta} \left(a_{1}, \delta(b_{2}), \ldots,\delta(b_{n}) \right) &:= f^{A(\alpha)/\Delta_{\alpha \alpha}}\left( a_{1},\delta(b_{2}),\ldots, \delta(b_{n}) \right) \\
T_{f}(x_1,\ldots,x_{n}) &:= \begin{bmatrix} l(f^{Q}(x_{1},\ldots,x_{n})) \\ f^{A}(l(q_{1}),\ldots,l(q_{n})) \end{bmatrix}/\Delta_{\alpha \alpha}
\end{split}
\end{align}
determine affine datum $(Q, A^{\alpha,\tau},\ast)$ where 
\begin{align}
\phi(x) := \begin{bmatrix} l \circ \pi(x) \\ x \end{bmatrix}/\Delta_{\alpha \alpha} 
\end{align}
yields an isomorphism $\phi: A \rightarrow A_{T}(Q, A^{\alpha,\tau},\ast)$.

This motivates the definition of realization. An algebra $B$ \emph{realizes} the affine datum $(Q,A^{\alpha,\tau},\ast)$ if there is an extension $\pi: B \rightarrow Q$ with $\beta = \ker \pi$, a bijection $i : B(\beta)/\Delta_{\beta \beta} \rightarrow A(\alpha)/\Delta_{\alpha \alpha}$ and a lifting $l: Q \rightarrow B$ such that for all $f \in \tau$: if $n = \ar f$, $q_{1},\ldots,q_n \in Q$ and $x_{1},\ldots,x_{n} \in B(\beta)/\Delta_{\beta \beta}$, then
\begin{align*}
a(f,i)\left(q_1,\ldots, i \circ x_{i},\ldots, q_{n} \right) = i \circ f^{B(\beta )/\Delta_{\beta \beta}} \left( \delta \circ l(q_1) ,\ldots, x_{i},\ldots,\delta \circ l(q_{n}) \right).
\end{align*}

Terms in the algebra $A_{T}(Q,A^{\alpha,\tau},\ast)$ can be thought of as interpretations of terms from the multisorted signature $T \cup \{f^{\Delta}: f \in \tau \} \cup \{ a(f,i): f \in \tau, \ar f \geq 2, i \in [\ar f] \}$. Using the homomorphic property to distribute the operations over each other, we see from \cite[Lem 3.18]{wiresI} that for any term $t(\vec{x})$ in the signature $\tau$ there are terms $t^{\partial, T}$ and $(t^{\sigma})^{\ast}$ in the multisorted signature such that in the algebra $A_{T}(Q,A^{\alpha,\tau},\ast)$ we can represent the term as
\begin{align}\label{eqn:1}
F_{t} \left( \epsilon(\vec{x}) \right) =\sum_{\mu \in L(t,\epsilon(\vec{x}))} (t^{\sigma})^{\ast}(\mu ( \mathrm{var} \, t^{\sigma} ))  \ +_{u} \ t^{\partial, T}\left( \rho \circ \epsilon (\vec{x}) \right)
\end{align}
for any evaluation $\epsilon : \mathrm{var} \, t \rightarrow A(\alpha)/\Delta_{\alpha \alpha}$, $u = l \left( t^{Q}(\rho \circ \epsilon (\vec{x})) \right)$ for any lifting $l$ associated to the datum and $L(t,\epsilon(\vec{x}))$ is a set of evaluations related to $\epsilon$. We should note here that $t^{\partial, T}$ is the part of the interpretation which involves the operations of the 2-cocycle $T$ and $(t^{\sigma})^{\ast}$ only uses the pairings of the action and operations of $A^{\alpha,\tau}$. The representation in Eq~\eqref{eqn:1} allows us to discuss the equational theories of extensions.

The action in affine datum $(Q,A^{\alpha,\tau},\ast)$ is \emph{compatible} with a variety $\mathcal U$ in the signature $\tau$ if for all $f = g \in \mathrm{Id} \, \mathcal U$, 
\begin{align}\label{eqn:weakcomp}
 \sum_{\mu \in L(f,\epsilon(\vec{x}))} (f^{\sigma})^{\ast}(\mu ( \mathrm{var} \, f^{\sigma} )) = \sum_{\mu \in L(g,\epsilon(\vec{x}))} (g^{\sigma})^{\ast}(\mu ( \mathrm{var} \, g^{\sigma} )). 
\end{align} 

The 2-cocycle $T$ associated to affine datum $(Q,A^{\alpha,\tau},\ast)$ is \emph{compatible} with a variety $\mathcal U$ in the signature $\tau$ if  
\begin{itemize}

	\item for all $f \in\tau$, $T_{f}(q_1,\ldots,q_{\ar f}) \mathrel{\hat{\alpha}/\Delta_{\alpha \alpha}} a(f,1) \left( \begin{bmatrix} a \\ b \end{bmatrix}/\Delta_{\alpha \alpha}, q_2,\ldots, q_{\ar f} \right)$ where $\rho \left( \begin{bmatrix} a \\ b \end{bmatrix} \right) = q_{1}$;
	
	\item for all $t(\vec{x})=g(\vec{y}) \in \mathrm{Id} \, \mathcal U$ and evaluations $\epsilon : \mathrm{var} \, t \cup \mathrm{var} \, g \rightarrow Q$ we have 
\begin{align}
	t^{\partial,T}(\epsilon(\vec{x})) = g^{\partial,T}(\epsilon(\vec{y})).
\end{align}

\end{itemize}
It follows from \cite[Thm 3.21]{wiresI} that when the action in the affine datum $(Q,A^{\alpha,\tau},\ast)$ is compatible with $\mathcal U$, then algebras of the form $A_{T}(Q,A^{\alpha,\tau},\ast)$ for $\mathcal U$-compatible 2-cocycles $T$ are exactly the algebras in $\mathcal U$ which realize the given affine datum. In the case $\pi : A \rightarrow Q$ realizes affine datum and $A \approx A_{T}(Q,A^{\alpha,\tau},\ast)$, then $T$ is the 2-cocycle \emph{associated} to the extension $A$. We say the 2-cocycle $T$ is \emph{trivial} and write $T=0$ if each $T_{f}(\vec{x}) = \begin{bmatrix} l(f^{Q}(\vec{x})) \\ l(f^{Q}(\vec{x})) \end{bmatrix}/\Delta_{\alpha \alpha}$ for a lifting $l$ associated to the affine datum. Then for a trivial 2-cocycle, the 2-cocycle terms in the operations of an extensions which is displayed in Eq~\eqref{eqn:algebrarep} do not appear since they are zeros for the induced binary operation; therefore, $A_{0}(Q,A^{\alpha,\tau},\ast)$ becomes the unique semidirect product realizing the datum \cite[Prop 3.22]{wiresI}. This means that by realization the syntactically complicated condition that the action is compatible with $\mathcal U$ displayed in Eq~\eqref{eqn:weakcomp} is equivalent to $A_{0}(Q,A^{\alpha,\tau},\ast) \in \mathcal U$ which is membership of the semidirect product in the variety.

Different choices of lifting for an extension leads to a combinatorial equivalence on isomorphic extensions. A sequence of operations $G=\{G_{f}: f \in \tau\}$  where $G_{f}: Q^{\ar f} \rightarrow A(\alpha)/\Delta_{\alpha \alpha}$ is a \emph{2-coboundary} of datum $(Q,A^{\alpha,\tau},\ast)$ if there is a function $h: Q \rightarrow A(\alpha)/\Delta_{\alpha \alpha}$ such that for any lifting $l:Q \rightarrow A$ associated to the datum
\begin{itemize}

	\item $h(x) \mathrel{\hat{\alpha}/\Delta_{\alpha \alpha}} \delta \circ l(x)$;
	
	\item for each $f \in \tau$,
\begin{align}
G_{f}(x_1,\ldots,x_n) &= \sum_{i=1}^{n} a(f,i) ( x_1,\ldots,x_{i-1},h(x_{i}),x_{i+1},\ldots,x_n ) -_{u} h(f^{Q}(x_1,\ldots,x_n))
\end{align}

\end{itemize}
where $u = l(f^{Q}(x_1,\ldots,x_n))$. The function $h : Q \rightarrow A(\alpha)/\Delta_{\alpha \alpha}$ is said to \emph{witness} the 2-coboundary. For different liftings $l,l'$ of an extension $\pi: A \rightarrow Q$ realizing affine datum, $h(x) = \begin{bmatrix} l(x) \\ l'(x) \end{bmatrix}/\Delta_{\alpha \alpha}$ witnesses a 2-coboundary. Let $Z^{2}_{\mathcal U}(Q,A^{\alpha,\tau},\ast)$ denote the set of $\mathcal U$-compatible 2-cocycles and $B^{2}(Q,A^{\alpha,\tau},\ast)$ the set of 2-coboundaries. The definition 
\begin{align*}
(T_{f} + T'_{f})(\bar{x}) &:= T_{f}(\bar{x}) +_{l(f(\bar{x}))} T'_{f}(\bar{x})   &(f \in \tau, \bar{x} \in Q^{\ar f}) \\
\end{align*}
for any choice of lifting associated to the datum makes $Z^{2}_{\mathcal U}(Q,A^{\alpha,\tau},\ast)$ into an abelian group with subgroup $B^{2}(Q,A^{\alpha,\tau},\ast) \leq Z^{2}_{\mathcal U}(Q,A^{\alpha,\tau},\ast)$ \cite[Lem 3.26]{wiresI}. This determines the abelian $2^{\mathrm{nd}}$-cohomology group 
\begin{align}\label{eqn:2ndcohom}
H^{2}_{\mathcal U}(Q,A^{\alpha,\tau},\ast) := Z^{2}_{\mathcal U}(Q,A^{\alpha,\tau},\ast)/B^{2}(Q,A^{\alpha,\tau},\ast).
\end{align}
Two extensions realizing the datum are said to be \emph{equivalent} if their associated 2-cocycles determine the same cohomology class as defined in Eq~\eqref{eqn:2ndcohom}. The equivalence classes of extensions in $\mathcal U$ realizing the affine datum is in bijective correspondence with the abelian  $2^{\mathrm{nd}}$-cohomology group in which the zero corresponds to the class of the semidirect product.

The equivalence of extensions determined by 2-coboundaries can be characterized by a restricted class of isomorphisms \cite[Th 3.31]{wiresI}. Two extensions $\pi: A \rightarrow Q$ and $\pi' : A' \rightarrow Q$ realizing fixed affine datum are equivalent if and only if there is a \emph{stabilizing isomorphism} $\gamma: A \rightarrow A'$; that is, $\gamma$ is an isomorphism satisfying 
\begin{enumerate}

	\item $\pi' \circ \gamma = \pi$, and
	
	\item $\gamma = m(\gamma \circ r, r , \id)$ for all $\alpha$-traces $r: A \rightarrow A$.

\end{enumerate}
A \emph{derivation} (or 1-cocycle) of affine datum is a witness of the trivial 2-coboundary and the set of derivations forms an abelian group under the sum $(d + d')(x) = d(x) +_{l(x)} d'(x)$ for any lifting $l$ associated to the datum. The group of stabilizing automorphisms of any fixed extension realizing the affine datum is isomorphic to the abelian group of derivations \cite[Thm 3.34]{wiresI}.

In varieties $\mathcal V$ with a difference term, central extensions are characterized by \emph{trivial} actions \cite[Prop 3.37]{wiresI}; thus, for these actions $H^{2}_{\mathcal V}(Q,A^{\alpha,\tau},\ast)$ parametrizes the class of central extensions, but an alternative description is possible since in these varieties central extensions have a representation slightly different from Eq~\eqref{eqn:algebrarep}. We now derive modifications of the preceding machinery in the special case of central extensions in varieties with a difference term. We begin with the \emph{basic construction} for producing central extensions. Suppose we have two algebras $B$ and $Q$ in the same signature $\tau$ and a binary operation on $B$ denoted by $x + y$. Suppose further that for every operation symbol $f \in \tau$ we have an operation $T_{f} : Q^{\ar f} \rightarrow B$. We define a new algebra $B \otimes^{T} Q$ over the universe of the direct product $B \times Q$ where each operation symbol $f \in \tau$ is interpreted by the rule
\[
F_{f} \left( (b_1,q_1),\ldots, (b_n,q_n) \right) := \left\langle  f(b_1,\ldots,b_n) + T_{f}(q_1,\ldots,q_n) , f(q_1,\ldots,q_n) \right\rangle.
\]
We have an extension over $Q$ given by second-projection $p_{2} : B \otimes^{T} Q \rightarrow Q$. If $B$ is an abelian algebra with a Mal'cev term $m$ which is used to define the binary operation by $x + y:=m(x,0,y)$ for some choice $0 \in B$, then according to \cite[Lem 2.9]{wiresI} this is a central extension. In varieties with a difference term, special cases of the basic construction recovers all central extensions.

Let $A \in \mathcal V$ a variety with a difference term $m$ and $\alpha \in \Con A$ a central congruence. Since $\alpha$ is central, the diagonal elements of $A(\alpha)$ form a single $\Delta_{\alpha 1}$-class which allows us to conclude that  $\hat{\delta} = \left\{ \begin{bmatrix} a \\ a \end{bmatrix}/\Delta_{\alpha 1} : a \in A \right\}$ is an idempotent element in $A(\alpha)/\Delta_{\alpha 1}$. According to \cite[Lem 2.2]{wiresI}, $A(\alpha)/\Delta_{\alpha 1}$ is an abelian algebra and so if we make the definitions 
\begin{align*}
x+y := m(x,\hat{\delta},y) \quad \quad \text{ and } \quad \quad -x := m(\hat{\delta},x,\hat{\delta}), 
\end{align*}
then $\left\langle A(\alpha)/\Delta_{\alpha 1}, + , - ,\hat{\delta} \right\rangle$ is an abelian group \cite[Thm 7.34]{bergman} will neutral element $\hat{\delta}$; furthermore, the algebra $A(\alpha)/\Delta_{\alpha 1}$ is polynomially-equivalent to a module over a unital ring \cite[Thm 7.35]{bergman} in which the terms of the signature are interpreted by module terms because $\hat{\delta}$ is an idempotent element. Set $Q=A/\alpha$ and choose a lifting $l$ for the canonical epimorphism $\pi: A \rightarrow Q$. If we make the definition
\begin{align}\label{eqn:2cocycle}
T_{f}(x_{1},\ldots,x_{n}) := \begin{bmatrix} l(f^{Q}(x_{1},\ldots,x_{n})) \\ f^{A}(l(x_{1}),\ldots,l(x_{n})) \end{bmatrix}/\Delta_{\alpha 1} 
\end{align}
for each $f \in \tau$, then according to \cite[Cor 2.7]{wiresI} there is an isomorphism $\psi : A \rightarrow A(\alpha)/\Delta_{\alpha 1} \otimes^{T} Q$ of $A$ with the basic construction given by 
\begin{align}\label{eqn:centriso}
\psi(x) := \left\langle \begin{bmatrix} l \circ \pi (x) \\ x \end{bmatrix}/\Delta_{\alpha 1} , x/\alpha \right\rangle.
\end{align}
This motivates the notion of central datum and the following refinements to the cohomological machinery.

\begin{definition}\label{def:centdatum}
A pair of algebras $(Q,B^{\tau})$ is \emph{central datum} if
\begin{enumerate}

	\item[(CD1)] $B$ is a unital $R$-module for some ring $R$ with identity in which $B$ carries an interpretation of $\tau$ by module terms, denoted by $B^{\tau}$, and there is a ternary operation $m$ such that $m^{B}(x,y,z) = x - y + z$;
	 
	\item[(CD2)] $Q$ is an algebra in the signature $\tau$ in which there is an interpretation $m^{Q}$ which is idempotent.

\end{enumerate}
\end{definition}

By writing $B^{\tau}$, we stress the $\tau$-structure $\left\langle B, \{f^{B}: f \in \tau \} \right\rangle$ with the additional information that the interpretations of the signature can be calculated by module terms $f^{B}(x_{1},\ldots,x_{n}) = r_{1}x_{1} + \cdots + r_{n}x_{n}$ given by a unital module $_{R} B$ for a ring $R$ with identity; consequently, the zero of the module is an idempotent element for the $\tau$-operations. We refer to $B^{\tau}$ as the \emph{kernel algebra} in central datum. We say central datum $(Q,B^{\tau})$ is contained in a variety $\mathcal V$ if both $Q, B^{\tau} \in \mathcal V$.

A sequence of functions $T=\{ T_{f}: f \in \tau \}$ is a \emph{2-cocycle} appropriate to central datum $(Q,B^{\tau})$ if $T_{f}: Q^{\ar f} \rightarrow B$ for each $f \in \tau$. Given central datum $(Q,B^{\tau})$ and appropriate 2-cocycle $T$, we have an extension $p_{2} : B^{\tau} \otimes^{T} Q \rightarrow Q$ in which $\alpha = \ker p_{2}$ is a central congruence \cite[Lem 2.9]{wiresI}. Set $A := B^{\tau} \otimes^{T} Q$ and choose the particular lifting $l : Q \rightarrow A$ of $p_{2}$ given by $l(x) = \left\langle 0, x \right\rangle$. If $A := B^{\tau} \otimes^{T} Q$ is in a variety with a difference term, then $A \approx A(\alpha)/\Delta_{\alpha 1} \otimes^{T'} Q$ for a 2-cocycle $T'$ appropriate to central datum $(Q,A(\alpha)/\Delta_{\alpha 1})$ which is defined by Eq~\eqref{eqn:2cocycle}. By the unique representation of $\Delta_{\alpha 1}$-classes implied by Eq~\eqref{eqn:centriso}, the map $\sigma : A(\alpha)/\Delta_{\alpha 1} \rightarrow B^{\tau}$ given by $\sigma \left( \begin{bmatrix} \left\langle 0, x \right\rangle \\ \left\langle b, x \right\rangle \end{bmatrix}/\Delta_{\alpha 1} \right) := b$ is well-defined and bijective. Then for $f \in \tau$ with $n=\ar f$ we have 
\begin{align*}
f^{A(\alpha)/\Delta_{\alpha 1}} \left( \begin{bmatrix} \left\langle 0, x_{1} \right\rangle \\ \left\langle b_{1}, x_{1} \right\rangle \end{bmatrix}/\Delta_{\alpha 1}, \ldots, \begin{bmatrix} \left\langle 0, x_{n} \right\rangle \\ \left\langle b_{n}, x_{n} \right\rangle \end{bmatrix}/\Delta_{\alpha 1} \right) &= \begin{bmatrix} \left\langle T_{f}(\vec{x}) , f^{Q}(\vec{x}) \right\rangle \\ \left\langle f^{B}(\vec{b}) + T_{f}(\vec{x}) , f^{Q}(\vec{x}) \right\rangle \end{bmatrix}/\Delta_{\alpha 1} \\
&= \begin{bmatrix} \left\langle 0 , f^{Q}(\vec{x}) \right\rangle \\ \left\langle f^{B}(\vec{b}) , f^{Q}(\vec{x}) \right\rangle \end{bmatrix}/\Delta_{\alpha 1} \\
&\stackrel{\sigma}{\rightarrow} f^{B}(\vec{b})
\end{align*}
from which it follows that $\sigma$ is an isomorphism. We observe that the basic construction recovers all central extensions provided we have the correct conditions on the algebras $B^{\tau}$ and $Q$, and the 2-cocycle $T$ so that $B^{\tau} \otimes^{T} Q$ is in a variety with a difference term.

Let us examine the modified form 2-coboundaries, actions and compatibility with equational theories take in relation to the basic construction of central extensions. Write $\kappa: A(\alpha)/\Delta_{\alpha \alpha} \rightarrow A(\alpha)/\Delta_{\alpha 1}$ for the canonical epimorphism associated to the quotient by the congruence $\Delta_{\alpha 1}/\Delta_{\alpha \alpha}$. Then the isomorphism in Eq~\eqref{eqn:centriso} for central congruences in varieties with a difference term takes the form 
\begin{align}\label{eqn:isocentmod}
\psi(x) = \left\langle \kappa \circ \phi(x) , \pi(x) \right\rangle. 
\end{align}
Fix affine datum $(Q,A^{\alpha,\tau},\ast)$ in which the action is compatibile with $\mathcal V$, $Q \in \mathcal V$ and $T$ is a $\mathcal V$-compatible 2-cocycle appropriate for the datum where $\mathcal V$ has a difference term. Define $T^{\kappa}$ and $\ast^{\kappa}$ by
\begin{align*}
T_{f}^{\kappa}(\vec{x}) &:= \kappa \circ T_{f}(\vec{x})  &(\vec{x} \in Q^{\ar f}) \\
a^{\kappa}(f,i)(x_{1},\ldots,\kappa(a),\ldots,x_{\ar f}) &:= \kappa \circ a(f,i)(x_{1},\ldots,\kappa(a),\ldots,x_{\ar f})  &(a \in A(\alpha)/\Delta_{\alpha \alpha}, \vec{x} \in Q^{\ar f}) 
\end{align*}
Then for central extensions we have the isomorphism
\begin{align}\label{eqn:isorep}
A_{T}(Q,A^{\alpha,\tau},\ast) \approx A(\alpha)/\Delta_{\alpha 1} \otimes^{T^{\kappa}} Q
\end{align}
where we have written the algebra as $A = A_{T}(Q,A^{\alpha,\tau},\ast)$.

Let us consider what happens to the action pairings and the operations $f^{\Delta}$. Since $A(\alpha)/\Delta_{\alpha 1}$ is abelian with an idempotent element $\hat{\delta}$, there is a ring $R$ with identity such that $A(\alpha)/\Delta_{\alpha 1}$ is term-equivalent to the unital module $_{R} B$ where for the universe we are writing $B=A(\alpha)/\Delta_{\alpha 1}$. So for any $f \in \tau$ with $n=\ar f$ there exists ring elements $\{ r_{1},\ldots,r_{n} \} \subseteq R$ such that 
\begin{align*}
f^{A(\alpha)/\Delta_{\alpha 1}} \left(a_{1}, \ldots, a_{i} , \ldots, a_{n} \right) = r_{1} \cdot a_{1} + \cdots + r_{i} \cdot a_{i} + \cdots + r_{n} \cdot a_{n}.
\end{align*}
Then for the action pairings, by realization we have 
\begin{align*}
a^{\kappa}(f,i)(x_1,\ldots,x_{i-1},\kappa(a),x_{i+1},\ldots,x_n) &= \kappa \circ a(f,i)(x_1,\ldots,x_{i-1},a,x_{i+1},\ldots,x_n) \\
&= \kappa \circ f^{A(\alpha)/\Delta_{\alpha \alpha}} \left( \delta \circ l(x_{1}), \ldots, a , \ldots, \delta \circ l(x_{n}) \right)  \\
&= f^{A(\alpha)/\Delta_{\alpha 1}} \left( \kappa \circ \delta \circ l(x_{1}), \ldots, \kappa(a) , \ldots, \kappa \circ \delta \circ l(x_{n}) \right)  \\
&= r_{1} \cdot \hat{\delta} + \cdots + r_{i} \cdot \kappa(a) + \cdots + r_{n} \cdot \hat{\delta} \\
&= r_{i} \cdot \kappa(a). 
\end{align*}
The action pairings and operations $f^{\Delta}$ are encoded into the ring action of the module. In this way, compatibility with $\mathcal V$ reduces to the usual satisfaction of identities under the standard interpretation by module terms; that is, $A(\alpha)/\Delta_{\alpha 1} \in \mathcal V$. If we are starting from central datum $(Q,B^{\tau})$, this is the requirement that $B^{\tau} \in \mathcal V$.

The condition for the compatibility of a 2-cocycle with a variety (or equational theory) follows the case for affine datum, but the resulting equations will be simplified since it will be a composition of the operations of the 2-cocycle with the module terms. Fix central datum $(Q,B^{\tau})$ and appropriate 2-cocycle $T$. If we consider the term representation in Eq~\eqref{eqn:1} under the isomorphism Eq~\eqref{eqn:isocentmod} we arrive at the following: for each term $t(\vec{x})$ in the signature $\tau$, there is a term $t^{\partial, T} : Q^{\mathrm{var} \, t} \rightarrow B$ in the multisorted signature $T \cup \{r: r\in R\}$ such that for all evaluations $\epsilon: \mathrm{var} \, t \rightarrow B \times Q$, we have
\begin{align}\label{eqn:eqrepcentral}
F_{t} \left( \epsilon(\vec{x}) \right) = \left\langle \, f^{B}(p_{1} \circ \epsilon(\vec{x}) ) + t^{\partial,T}(p_{2} \circ \epsilon(\vec{x}) ) \ , \ t^{Q}(p_{2} \circ \epsilon(\vec{x})) \, \right\rangle
\end{align}
for the coordinate projections $p_{i}$ of the direct product set $B \times Q$. We say a 2-cocycle $T$ appropriate for central datum $(Q,B^{\tau})$ is \emph{compatible} with a variety $\mathcal V$ in the signature $\tau$ if for all $t(\vec{x})=g(\vec{y}) \in \mathrm{Id} \, \mathcal V$ and evaluations $\epsilon : \mathrm{var} \, t \cup \mathrm{var} \, g \rightarrow Q$ we have 
\begin{align}
	t^{\partial,T}(\epsilon(\vec{x})) = g^{\partial,T}(\epsilon(\vec{y})).
\end{align}

In a similar manner, 2-coboundaries under the isomorphism Eq~\eqref{eqn:isocentmod} leads to the following definition for central extensions: a sequence of operations $G=\{G_{f}: f \in \tau\}$ where $G_{f}: Q^{\ar f} \rightarrow B$ is a \emph{2-coboundary} of central datum $(Q,B^{\tau})$ if there is a function $h: Q \rightarrow B$ such that 
\begin{align}\label{eqn:centderiv}
G_{f}(\vec{x}) &= f^{B}(h(\vec{x})) - h(f^{Q}(\vec{x}))  &(f \in \tau, \vec{x} \in Q^{\ar f} ).
\end{align}
Then the analogous definitions for the groups of $\mathcal V$-compatible 2-cocycles $Z^{2}_{\mathcal V}(Q,B^{\tau})$ and 2-coboundaries $B^{2}(Q,B^{\tau})$ yields the abelian $2^{\mathrm{nd}}$-cohomology group 
\begin{align}
H^{2}_{\mathcal V}(Q,B^{\tau}) := Z^{2}_{\mathcal V}(Q,B^{\tau})/B^{2}(Q,B^{\tau})
\end{align}
which is in bijective correspondence with the equivalence classes of central extensions in the variety $\mathcal V$ realizing the central datum $(Q,B^{\tau})$ in which the zero corresponds to the class of the direct product. Given central datum $(Q,B^{\tau})$, if we assume $Q,B^{\tau} \in \mathcal V$ a variety with a difference term, then the central extensions in $\mathcal V$ which realize the datum are equivalent to algebras of the form $B^{\tau} \otimes^{T} Q$ for a $\mathcal V$-compatible 2-cocycle $T$. The definition of stabilizing isomorphisms characterizing equivalence classes of extensions remains unchanged but the defining condition for derivations (1-cocycles) witnessing trivial 2-coboundaries from Eq~\eqref{eqn:centderiv} reduces to stating homomorphisms from $Q$ to $B^{\tau}$.

\begin{example}
As an example, let us determine the derivations from stabilizing automorphism. Consider the extension $p_{2}: B^{\tau} \otimes^{T} Q \rightarrow Q$ in which $B^{\tau} \otimes^{T} Q \in \mathcal V$ has a difference term $m$; in particular, $m$ is idempotent in algebras in the variety. We calculate
\begin{align*}
\left\langle 0, x \right\rangle =  m \left( \left\langle 0, x \right\rangle, \left\langle 0, x \right\rangle, \left\langle 0, x \right\rangle \right) = \left\langle m^{B}(0,0,0) + m^{\partial, T}(x,x,x) , m^{Q}(x,x,x) \right\rangle = \left\langle m^{\partial, T}(x,x,x) , x \right\rangle
\end{align*}
from which we see that $m^{\partial, T}(x,x,x) = 0$.

Let $\gamma$ be a stabilizing automorphism of the extension $p_{2}: B^{\tau} \otimes^{T} Q \rightarrow Q$ and take the $\alpha$-trace $r(b,q) = \left\langle 0,q \right\rangle$. The first condition on a stabilizing automorphism $p_{2} \circ \gamma = p_{2}$ implies we can write $\gamma(b,q) = \left\langle g(b,q),q \right\rangle$ for some function $g: B \times Q \rightarrow Q$. Then the second condition yields
\begin{align*}
\gamma(b,q) = m \big( \gamma \circ r (b,q) , r(b,q) , \left\langle  b,q \right\rangle \big) = m \big( \left\langle  g(0,q), q \right\rangle , \left\langle 0,q \right\rangle , \left\langle b,q \right\rangle \big) &= \left\langle g(0,q) + b + m^{\partial, T}(q,q,q) , q \right\rangle \\
&= \left\langle g(0,q) + b , q \right\rangle .
\end{align*}
Calculating
\begin{align*}
&\left\langle   f^{B} \left( g(0,q_1),\ldots,g(0,q_n) \right) + f^{B}(b_1,\ldots,b_{n}) + T_{f}(q_1,\ldots,q_{n}), f^{Q}(q_1,\ldots,q_{n})  \right\rangle \\
&= F_{f} \left( \left\langle g(0,q_{1}) + b_{1}, q_{1} \right\rangle,\ldots, \left\langle g(0,q_{n}) + b_{n}, q_{n} \right\rangle \right) \\
&= F_{f} \left( \gamma(b_{1},q_{1}),\ldots,\gamma(b_{n},q_{n}) \right) \\
&= \gamma \left( F_{f}\left( \left\langle b_{1},q_{1} \right\rangle,\ldots,\left\langle b_{n},q_{n}\right\rangle \right) \right)\\
&= \gamma \left( \left\langle f^{B}(b_{1},\ldots,b_{n}) + T_{f}(q_{1},\ldots,q_{n}) , f^{Q}(q_{1},\ldots,q_{n}) \right\rangle    \right)  \\
&=  \left\langle  g(0,f^{Q}(q_{1},\ldots,q_{n})) + f^{B}(b_{1},\ldots,b_{n}) + T_{f}(q_{1},\ldots,q_{n}) , f^{Q}(q_{1},\ldots,q_{n}) \right\rangle .
\end{align*}
We recover the defining condition
\[
g \left( 0,f^{Q}(q_{1},\ldots,q_{n}) \right) = f^{B} \left( g(0,q_1),\ldots,g(0,q_n) \right)
\]
showing that $d_{\gamma}(x) := g(0,x)$ is a derivation (1-cocycle). 
\end{example}

An element $u \in A \in \mathcal V$ is \emph{idempotent} if $f(u,u,\cdots u) = u$ for all operation symbols $f \in \tau$. For any congruence $\alpha \in \Con A$, we let $I_{\alpha} = \{ a : a \mathrel{\alpha} u \}$ denote the $\alpha$-congruence class which contains $u$ and note it is a subalgebra of $A$ when $u$ is idempotent. If $m$ is a difference term for $\mathcal V$ and $\alpha \in \Con A$ is central, then we have an isomorphism
\begin{align}\label{eqn:idemideal}
A(\alpha)/\Delta_{\alpha 1} &\approx I_{\alpha} \ \ \ \text{witnessed by the map} \ \ \begin{bmatrix} a \\ b \end{bmatrix}/\Delta_{\alpha 1} \longmapsto m(a,b,u);
\end{align}
in addition, for $\beta \in \Con A$ the image of the subalgebra $A(\alpha \wedge \beta) \leq A(\alpha)$ under the canonical epimorphism $\kappa_{\alpha 1} : A(\alpha) \rightarrow A(\alpha)/\Delta_{\alpha 1}$ is given by $\kappa_{\alpha 1} ( A(\alpha \wedge \beta) ) \approx I_{\alpha \wedge \beta}$ (\cite[Prop 4.3]{wiresI}).

\vspace{0.3cm}

\section{Low-dimensional Hochschild-Serre sequence}\label{sec:4}
\vspace{0.3cm}

In this section, we prove an analogue of the Hochschild-Serre exact sequence (see Hochschild and Serre \cite[Thm ]{hochserre}) in dimension 1 which extends the special case in group theory of a central extension with abelian coefficients to a central extension with central datum in varieties with a difference term. We assume in our varieties with a difference term $m(x,y,z)$ that the term $m$ is a fundamental symbol of the signature.

\begin{definition}
Let $B^{\tau}$ and $E^{\tau}$ be kernel algebras. Define $\Hom(B^{\tau},E^{\tau})$ to be the set of homomorphisms $\phi: B^{\tau} \rightarrow E^{\tau}$ such that $\phi(0)=0$; that is, homomorphisms which preserve the zeros of the respective modules. $\Hom(B^{\tau},E^{\tau})$ becomes an abelian group under the induced abelian group operation in the codomain.
\end{definition}

Suppose $(Q,B^{\tau})$ and $(P,E^{\tau})$ are central datum contained in a variety $\mathcal V$ with a difference term $m(x,y,z)$. As central datum, by the (CD1)-condition there are operations $p : B^{3} \rightarrow B$ and $q : E^{3} \rightarrow E$ such that $p^{B}$ is the affine operation for a unital module $_{R} B$ and $q^{E}$ is the affine operation for a unital module $_{S} E$. Since the operations of $B^{\tau}$ and $E^{\tau}$ are interpreted by module terms, they are both abelian algebras and so the difference term $m$ interprets as a compatible Mal'cev operation in each algebra. Since a compatible Mal'cev operation on a set is unique, the difference term interprets as the affine operation in each kernel algebra; that is, $m^{B}(x,y,z) = p^{B}(x,y,z)$ and $m^{E}(x,y,z) = q^{E}(x,y,z)$. The difference term is idempotent, so it satisfies the (CD2)-condition for central datum; therefore, for central datum in a variety we can always assume the referenced operation is given by the same difference term of the variety. The abelian group operations in the modules are given by $x+y = m^{B}(x,0,y)$ in $_{R} B$ for the zero $0 \in B$ and $x+y = m^{E}(x,0',y)$ in $_{S} E$ for the zero $0' \in E$; as a result, any $\phi \in \Hom(B^{\tau},E^{\tau})$ is a homomorphism of the abelian group reduct of the modules since $\phi$ preserves the zeros by definition.

Let $(Q,B^{\tau})$ and $(Q,E^{\tau})$ be central datum and $\pi: A \rightarrow Q$ a central extension realizing $(Q,B^{\tau})$ with $A \in \mathcal U$ a variety with a difference term. The \emph{inflation} maps
\begin{align*}
\check{\sigma} &: \Hom(Q,E^{\tau}) \longrightarrow \Hom(A,E^{\tau})  \\
\check{\sigma} &: H^{2}_{\mathcal U}(Q,E^{\tau}) \longrightarrow H^{2}_{\mathcal U}(A,E^{\tau})
\end{align*}
are defined by $\check{\sigma}(\phi) := \phi \circ \pi$ and $\check{\sigma}([T]) := [T \circ \pi]$ where $(T \circ \pi)_{f}(a_{1},\ldots,a_{\ar f}) = T_{f}(\pi(a_{1}),\ldots,\pi(a_{\ar f}))$ for $f \in \tau, (a_1,\ldots,a_{\ar f}) \in A^{\ar f}$. It is straightforward to see that the inflation maps are well-defined homomorphisms. In order to see that $T \circ \pi$ is $\mathcal U$-compatible whenever $T$ is $\mathcal U$-compatible, note that we can represent $A = B^{\tau} \otimes^{S} Q$. Then for any term $t$ and evaluation $\epsilon : \mathrm{var} \, t \rightarrow A$ we see that $t^{\partial, T \circ \pi}(\epsilon(\bar{x})) = t^{\partial, T}(p_{2} \circ \epsilon(\bar{x}))$. It then follows from Eq~\eqref{eqn:eqrepcentral} that $T \circ \pi$ is $\mathcal U$-compatible.

The \emph{restriction} map 
\[
\check{r} : \Hom(A,E^{\tau})  \longrightarrow \Hom(B^{\tau},E^{\tau}) 
\]
is defined by $\check{r}(\phi)(b) := \phi(b,q) - \phi(0,q)$ for any choice of $q \in Q$. The definition explicitly utilizes the identification of $A=B^{\tau} \otimes^{T} Q$. The map is well-defined: for any $b \in B, q,p \in Q$ we have
\begin{align*}
\phi(b,q) - \phi(0,q) &= \phi(m^{A}((b,p),(0,p),(0,q))) - \phi(0,q) \displaybreak[0]\\
&= m^{E}\left( \phi(b,p), \phi(0,p), \phi(0,q) \right) \displaybreak[0]\\
&= \phi(b,p) - \phi(0,p) + \phi(0,q) - \phi(0,q) = \phi(b,p) - \phi(0,p). 
\end{align*}
We verify that $\check{r}(\phi)$ is in the codomain : take $f \in \tau$ with $n = \ar f$, for any $(b_{1},q_{1}),\ldots, (b_{n},q_{n}) \in B \times Q$,
\begin{align*}
f^{E} \left( \check{r}(\phi)(b_{1}),\ldots, \check{r}(\phi)(b_{n}) \right) &=  f^{E} \left( \phi(b_{1},q_{1}) - \phi(0,q_{1}),\cdots, \phi(b_{n},q_{n}) - \phi(0,q_{n} \right) \displaybreak[0]\\
&= f^{E} \left( \phi(b_{1},q_{1}),\ldots,\phi(b_{n},q_{n}) \right) -  f^{E} \left( \phi(0,q_{1}),\ldots,\phi(0,q_{n}) \right)  \displaybreak[0]\\
&= \phi \left( f^{A} \left( (b_{1},q_{1}),\ldots,(b_{n},q_{n}) \right)  \right)  - \phi \left( f^{A} \left( (0,q_{1}),\ldots,(0,q_{n}) \right) \right) \displaybreak[0]\\
&= \phi \left( f^{B}(b_{1},\ldots,b_{n}) + T_{f}(q_{1},\ldots,q_{n}), f^{Q}(q_{1},\ldots,q_{n})  \right) \displaybreak[0]\\
&\quad - \phi \left( T_{f}(q_{1},\ldots,q_{n}), f^{Q}(q_{1},\ldots,q_{n}) \right) \displaybreak[0]\\
&= \phi \left( f^{B}(b_{1},\ldots,b_{n}), f^{Q}(q_{1},\ldots,q_{n})  \right) - \phi \left(0, f^{Q}(q_{1},\ldots,q_{n})  \right) \displaybreak[0]\\
&\quad + \phi \left( T_{f}(q_{1},\ldots,q_{n}), f^{Q}(q_{1},\ldots,q_{n})  \right) -  \phi \left( T_{f}(q_{1},\ldots,q_{n}), f^{Q}(q_{1},\ldots,q_{n}) \right) \displaybreak[0]\\
&=  \phi \left( f^{B}(b_{1},\ldots,b_{n}), f^{Q}(q_{1},\ldots,q_{n})  \right) - \phi \left(0, f^{Q}(q_{1},\ldots,q_{n})  \right) \displaybreak[0]\\
&= \check{r}(\phi)(f^{B}(b_{1},\ldots,b_{n}));
\end{align*}
thus, $\check{r}(\phi) \in \Hom(B^{\tau},E^{\tau})$.

The \emph{transgression} is the map 
\[
\delta: \Hom (B^{\tau}, E^{\tau}) \times H^{2}_{\mathcal U}(Q,B^{\tau}) \longrightarrow H^{2}_{\mathcal U}(Q,E^{\tau})
\] 
defined by $\delta(\phi, [T]) := [\phi \circ T]$. If $[T]=[T']$, then by Eq~\eqref{eqn:centderiv} there is $h: Q \rightarrow B$ such that $T_{f} - T_{f}' = f^{B} \circ h^{n} - h \circ f^{Q}$. Then $\phi \circ T_{f} - \phi \circ T_{f}' = f^{E} ((\phi \circ h)^{n}) - (\phi \circ h ) \circ f^{Q}$ which shows $[\phi \circ T_{f}]=[\phi \circ T_{f}']$; thus, the transgression is well-defined.

Let $[T] \in H^{2}_{\mathcal U}(Q,B^{\tau})$ and take $t = s \in \mathrm{Id} \, \mathcal U$. Then $T$ is a 2-cocycle compatible with $\mathcal U$ which implies $t^{\partial, T} = s^{\partial, T}$. Note $t^{\partial,T}$ is a sum with summands of the form $f^{B} \left(T_{h_{1}}(\bar{g}_{1}^{Q}),\ldots,T_{h_{n}}(\bar{g}_{n}^{Q}) \right)$ where from the composition tree we can write $t = f(h_{1}(\bar{g}_{1}),\ldots,h_{n}(\bar{g}_{n}))$ for terms $f,g_{ij}$ and $h_i \in \tau$ where $\bar{g}_{i} = (g_{i1},\ldots,g_{i \ar h_{i}})$. A similar statement holds for $s^{\partial, T}$. Then $t^{\partial, \phi \circ T}$ is a sum with summands of the form $f^{E} \left(\phi \circ T_{h_{1}}(\bar{g}_{1}^{Q}),\ldots,\phi \circ T_{h_{n}}(\bar{g}_{n}^{Q}) \right)$; thus, $\phi(t^{\partial,T}) = t^{\partial, \phi \circ T}$ and similarly $\phi(s^{\partial,T}) = s^{\partial, \phi \circ T}$. Then by Eq~\eqref{eqn:eqrepcentral} we see that $\phi \circ T$ is $\mathcal U$-compatible. It is straightforward to see that the transgression is a homomorphism in both coordinates.

In order to give a more elegant sufficiency condition for exactness in Theorem~\ref{thm:inflationrestriction}, we will introduce the liftings of homomorphisms which will be addressed further in Section~\ref{sec:5}.

\begin{definition}
Let $\pi: A \rightarrow Q$ and $\rho: E \rightarrow P$ be extensions. A homomorphism $\tau: Q \rightarrow P$ can be \emph{lifted to} $A$ \emph{through }$\rho$ if there is a homomorphism $\bar{\tau} : A \rightarrow E$ such that $\tau \circ \pi = \rho \circ \bar{\tau}$. The extension $\pi: A \rightarrow Q$ has the $\rho$-\emph{lifting property} if every homomorphism $\tau: Q \rightarrow P$ can be lifted to $A$ through $\rho$.

Let $\mathcal U$ be a variety and $A \in \mathcal U$. An extension $\pi: A \rightarrow Q$ has the $\mathcal U$-\emph{central lifting property} for $Q$ if $\pi$ has the $\rho$-lifting property for all central extensions $\rho: E \rightarrow P$ with $E \in \mathcal U$.
\end{definition}

In Section~\ref{sec:5}, we will be most interested in the case when $\pi: A \rightarrow Q$ is a central extension with the $\mathcal U$-\emph{lifting property} for a variety $\mathcal U$ with a difference term.

\begin{lemma}\label{lem:100}
Let $Q \in \mathcal U$ a variety with a difference term and $F_{\mathcal U}(X)/\theta \approx Q$ a free presentation in $\mathcal U$.
\begin{enumerate}	

	\item The central extension $\pi': F_{\mathcal U}(X)/[\theta,1] \rightarrow Q$ has the $\mathcal U$-central lifting property for $Q$. 

	\item The following are equivalent:
	
		\begin{enumerate}
		
			\item $F_{\mathcal U}(X)/[\theta,1]$ has an idempotent element;
			
			\item every central extension of $Q$ in $\mathcal U$ has an idempotent element;
			
			\item there is an idempotent $0' \in Q$	such that for all $B^{\tau} \in \mathcal U$ and $[T] \in H^{2}_{\mathcal U}(Q,B^{\tau})$ there is a $\mathcal U$-compatible 2-cocycle $T'$ with $[T]=[T']$ where $T'_{f}(0',\ldots,0')=0$ for all $f \in \tau$.
		
		\end{enumerate}

\end{enumerate}
\end{lemma}
\begin{proof}
(1) Write $F = F_{\mathcal U}(X)$ and the canonical surjection $\pi: F \rightarrow Q$ with $\theta = \ker \pi$. By the homomorphism property, we see that
\begin{align*}
\Big[ \theta/[\theta , 1] , 1_{F/[\theta,1]} \Big] = \Big[ \theta/[\theta , 1] , 1/[\theta , 1]  \Big] = \Big( [\theta , 1] \vee [\theta , 1] \Big) / [\theta , 1] = 0_{F/[\theta,1]}
\end{align*}
which shows $\theta/[\theta,1]$ is central in $F/[\theta,1]$. By the $2^{\mathrm{nd}}$-isomorphism theorem, there is a central extension $\pi' : F/[\theta,1] \rightarrow Q$ with $\theta/[\theta, 1] = \ker \pi'$. We show $\pi'$ has the $\rho$-lifting property for all central extensions $\rho: A \rightarrow P$ with $A \in \mathcal U$. 

Suppose there is homomorphism $\gamma : Q \rightarrow P$. Since $\rho$ is surjective, there is a homomorphism $\sigma : F \rightarrow A$ such that $\sigma(x) \in \rho^{-1}(\gamma \circ \pi(x))$ for each $x \in X$. Then $\gamma \circ \pi = \rho \circ \sigma$ and $\theta = \ker \pi$ implies $\sigma(\theta) \subseteq \ker \rho$. Then applying the homomorphism property of the commutator to the surjection $\sigma : F \rightarrow \sigma(F)$, we have
\begin{align*}
\left( [\theta,1] \vee \ker \sigma  \right)/\ker \sigma = \left[ (\theta \vee \ker \sigma)/\ker \sigma, 1_{F}/\ker \sigma \right] = \left[ \sigma(\theta),1_{\sigma(F)} \right] &\leq  [\ker \rho|_{\sigma(F)},1_{\sigma(F)}] \\
&\leq [\ker \rho,1_{A}]|_{\sigma(F)} = 0_{\sigma(F)}
\end{align*}
which implies $[\theta,1] \leq \ker \sigma$. This induces a homomorphism $\sigma': F/[\theta,1] \rightarrow A$ such that $\gamma \circ \pi' = \rho \circ \sigma'$. 

(2) First, assume $u \in F/[\theta,1]$ is an idempotent element. If $\rho: A \rightarrow Q$ is a central extension in $\mathcal U$, then by part (1), there is a lifting $\phi: F/[\theta,1] \rightarrow A$ of the identity $\id_{Q}$. Then we see that $\phi(u)$ is an idempotent element in $A$ since $\phi$ is a homomorphism; clearly, (b) implies (a) as a particular central extension.

Now assume every central extensions of $Q$ in $\mathcal U$ has an idempotent element. Take idempotent $u \in F/[\theta,1]$ and set $0' = \pi'(u) \in Q$ which is also idempotent in $Q$. Now take $[T] \in H^{2}_{\mathcal U}(Q,B^{\tau})$. Then $\phi(u)$ is an idempotent element in $B^{\tau} \otimes^{T} Q$ where $\phi : F/[\theta,1] \rightarrow B^{\tau} \otimes^{T} Q$ is a lifting of $\id_{Q}$. Let $l': Q \rightarrow B \times Q$ be defined by $l'(x):=\left\langle 0, x \right\rangle$ for $x \neq 0'$ and $l'(0') = \phi(u)$. Then $l'$ is a lifting of $p_{2}:B^{\tau} \otimes^{T} Q \rightarrow Q$, and if $T'$ is the 2-cocycle determined by $l'$ we see that $[T]=[T']$. Then for all $f \in \tau$ we have 
\begin{align*}
T_{f}(0',\ldots,0') &= f^{B^{\tau} \otimes^{T} Q}(l(0'),\ldots,l(0')) - l'(f^{Q}(0')) \\
&= f^{B^{\tau} \otimes^{T} Q}(\phi(u),\ldots,\phi(u)) - l'(0') \\
&= \phi(u) - \phi(u) = 0 
\end{align*}

Now assume condition (c) holds. If $\rho : A \rightarrow Q$ is central extension in $\mathcal U$, then there is a $\mathcal U$-compatible 2-cocycle $T$, kernel algebra $B^{\tau}$ and an isomorphism $\psi : A \rightarrow B^{\tau} \otimes^{T} Q$. By assumption, there is a $\mathcal U$-compatible 2-cocycle $T'$ such that $[T]=[T']$ and $T'_{f}(0',\ldots,0')=0$ for all $f \in \tau$. It follows that $\left\langle 0, 0' \right\rangle$ is idempotent in $B^{\tau} \otimes^{T'} Q$. Since $[T]=[T']$, there is a stabilizing isomorphism $\phi: B^{\tau} \otimes^{T'} Q \rightarrow B^{\tau} \otimes^{T} Q$; therefore, $\psi^{-1} \circ \phi (0,0')$ is idempotent in $A$.
\end{proof}

Let us recall that a \emph{complex} of abelian groups is a sequence of homomorphisms
\[
 \cdots \stackrel{f_{-2}}{\longrightarrow} A_{-2} \stackrel{f_{-1}}{\longrightarrow} A_{-1} \stackrel{f_{0}}{\longrightarrow} A_{0} \stackrel{f_{1}}{\longrightarrow} A_{1} \stackrel{f_{2}}{\longrightarrow} A_{2} \stackrel{f_{3}}{\longrightarrow} \cdots
\]
such that $f_{i+1} \circ f_{i} = 0$. It is \emph{exact} at the i-th coordinate if $\im f_{i} = \ker f_{i+1}$ and \emph{finite} if only finitely many groups $A_{i}$ are non-trivial.

\begin{theorem}\label{thm:inflationrestriction}
Let $A \in \mathcal U$ a variety with a difference term and $\pi: A \rightarrow Q$ a central extension realizing central datum $(Q,B^{\tau})$. If $E^{\tau} \in \mathcal U$ is a kernel algebra, then 
\begin{align}\label{eq:hochserre}
0 \rightarrow \Hom(Q,E^{\tau}) \stackrel{\check{\sigma}}{\longrightarrow} \Hom(A,E^{\tau}) \stackrel{\check{r}}{\longrightarrow} \Hom ( B^{\tau}, E^{\tau} ) \stackrel{\delta}{\longrightarrow} H^{2}_{\mathcal U}(Q,E^{\tau}) \stackrel{\check{\sigma}}{\longrightarrow} H^{2}_{\mathcal U}(A,E^{\tau})
\end{align}
is a finite complex of abelian groups which is exact at the first three groups. If there is a free presentation $F_{\mathcal U}(X)/\theta \approx Q$ such that $F_{\mathcal U}(X)/[\theta,1]$ has an idempotent element, then the sequence is exact. 
\end{theorem}
\begin{proof}
We first show exactness at the first three groups. 

$\ker \check{\sigma} = 0$: Take a homomorphism $\phi: Q \rightarrow E^{\tau}$ such that $\phi \in \ker \check{\sigma}$. Then by surjectivity of $\pi$, for any $q \in Q$ there is $a \in A$ which gives $0 = \check{\sigma}(\phi)(a) = \phi \circ \pi (a) = \phi(q)$; thus, $\phi = 0$.

$\Hom(A,E^{\tau})$: Take a homomorphism $\phi : Q \rightarrow E^{\tau}$. Then 
\[
\check{r} \circ \check{\sigma}(\phi)(b)  = \check{r}(\phi \circ \pi)(b) = \phi \circ \pi (b,q) - \phi \circ \pi (0,q) = \phi(q) - \phi(q) = 0
\] 
which shows $\check{r} \circ \check{\sigma} = 0$. Now take $\phi \in \ker \check{r}$. Then $0 = \check{r}(\phi)(b) = \phi(b,q) - \phi(0,q)$ for all $b \in B, q \in Q$. Define $\gamma: Q \rightarrow E$  by $\gamma (q) := \phi(0,q)$. For $f \in \tau$ with $n=\ar f$, 
\begin{align*}
\gamma(f^{Q}(q_{1},\ldots,q_{n})) = \phi(0,f^{Q}(q_1,\ldots,q_{n})) &= \phi(T_{f}(q_{1},\ldots,q_{n}),f^{Q}(q_1,\ldots,q_{n}) ) \\
&= \phi \left( f^{A} \left( (0,q_{1}),\ldots,(0,q_{n}) \right) \right) \\
&= f^{E} \left( \phi(0,q_{1}),\ldots, (0,q_{n}) \right) \\
&= f^{E} (\gamma(q_{1}),\ldots,\gamma(q_{n}))
\end{align*}
and so $\gamma$ is a homomorphism. Now observe $\check{\sigma}(\gamma)(b,q) = \gamma \circ \pi (b,q) = \gamma(q) = \phi(0,q) = \phi(b,q)$ for any $b \in B, q \in Q$; thus, $\phi \in \im \check{\sigma}$.

$\Hom ( B^{\tau}, E^{\tau} )$: Take homomorphism $\phi : A \rightarrow E^{\tau}$. Then for any $f \in \tau$ with $n = \ar f$, $\vec{q} \in Q^{n}$, and any $p \in Q$ we have
\[
\left( \delta \circ \check{r}(\phi) \right)_{f}(\vec{q}) = \check{r}(\phi) \circ T_{f}(\vec{q}) = \phi(T_{f}(\vec{q}),p) - \phi(0,p).
\]
Define $\psi : Q \rightarrow E$ by $\psi(q) := \phi(0,q)$. Then 
\begin{align*}
f^{E}(\psi(q_{1}),\ldots,\psi(q_{n})) = f^{E}(\phi(0,q_{1}),\phi(0,q_{n})) &= \phi \left( f^{A} \left( (0,q_{1}),\ldots,(0,q_{n}) \right) \right) \\
&= \phi \left( T_{f}(\vec{q}),f^{Q}(\vec{q}) \right)  \\
&= \phi \left( T_{f}(\vec{q}),f^{Q}(\vec{q}) \right) -  \phi \left( 0,f^{Q}(\vec{q}) \right) +  \phi \left( 0,f^{Q}(\vec{q}) \right) \\
&= \left( \delta \circ \check{r}(\phi) \right)_{f}(\vec{q}) + \psi(f^{Q}(\vec{q})).
\end{align*}
This shows $\delta \circ \check{r}(\phi) = 0$ in $H^{2}_{\mathcal U}(Q,E^{\tau})$. Now suppose $\phi \in \ker \delta$. Then there exists $h: Q \rightarrow E$ such that for all $f \in \tau$ with $n=\ar f$, $\vec{q} \in Q^{n}$ we have
\[
\phi \circ T_{f}(\vec{q}) = f^{E}(h(q_{1}),\ldots,h(q_{n})) - h(f^{Q}(\vec{q})).
\]
Define $\psi: A \rightarrow E$ by $\psi(b,q) := \phi(b) + h(q)$. Then for $f \in \tau$ with $n = \ar f$ we have
\begin{align*}
f^{E} \left( \psi(b_{1},q_{1}),\ldots,(b_{n},q_{n})  \right) &= f^{E} \left( \phi(b_{1}),\ldots,\phi(b_{n}) \right) + f^{E} \left( h(q_{1}),\ldots,h(q_{n}) \right) \\
&= \phi \left(f^{B}(b_{1},\ldots,b_{n}) \right) +  \phi \circ T_{f}(\vec{q}) +  h(f^{Q}(\vec{q})) \\
&= \phi \left( f^{B}(b_{1},\ldots,b_{n}) + T_{f}(\vec{q}) \right) + h(f^{Q}(\vec{q})) \\
&= \psi \left( f^{B}(b_{1},\ldots,b_{n}) + T_{f}(\vec{q}), f^{Q}(\vec{q})  \right)  \\
&= \psi \left( f^{A} \left( (b_{1},q_{1}),\ldots,(b_{n},q_{n}) \right)  \right)
\end{align*}
This shows $\psi$ is a homomorphism. Then $\check{r}(\psi)(b) = \psi(b,q) - \psi(0,q) = \phi(b) + h(q) - \phi(0) - h(q) = \phi(b)$ for all $b \in B,q \in Q$; thus, $\phi \in \im \check{r}$.

To complete the demonstration that the sequence is a complex, we must show $\check{\sigma} \circ \delta = 0$. If we write $\alpha = \ker \pi$, then we have $A(\alpha)/\Delta_{\alpha 1} \approx B^{\tau}$. Fix a $\alpha$-trace $r: A \rightarrow A$. Define $s: A \rightarrow B$ by $s(x):= \begin{bmatrix} x \\ r(x) \end{bmatrix}/\Delta_{\alpha 1}$. We can write for any $f \in \tau$ with $n = \ar f$, 
\begin{align*} 
T_{f}(\pi(a_{1}),\ldots,\pi(a_{n})) &= \begin{bmatrix}  r \circ f (a_{1},\ldots,a_{n}) \\ f(r(a_{1}),\ldots,r(a_{n})) \end{bmatrix}/\Delta_{\alpha 1} \\
&= m^{A(\alpha)} \left( \begin{bmatrix} f(a_{1},\ldots,a_{n}) \\ f(r(a_{1}),\ldots,r(a_{n})) \end{bmatrix}, \begin{bmatrix} f(a_{1},\ldots,a_{n}) \\ r \circ f(a_{1},\ldots,a_{n}) \end{bmatrix}, \begin{bmatrix} r \circ f(a_{1},\ldots,a_{n}) \\ r \circ f(a_{1},\ldots,a_{n}) \end{bmatrix} \right)/\Delta_{\alpha 1} \\
&= \begin{bmatrix} f(a_{1},\ldots,a_{n}) \\ f(r(a_{1}),\ldots,r(a_{n})) \end{bmatrix}/\Delta_{\alpha 1} - \begin{bmatrix} f(a_{1},\ldots,a_{n}) \\ r \circ f(a_{1},\ldots,a_{n}) \end{bmatrix}/\Delta_{\alpha 1} \\
&= f^{B} \left( \begin{bmatrix} a_{1} \\ r(a_{1}) \end{bmatrix}/\Delta_{\alpha 1},\ldots, \begin{bmatrix} a_{1} \\ r(a_{1}) \end{bmatrix}/\Delta_{\alpha 1}  \right) -  \begin{bmatrix} f(a_{1},\ldots,a_{n}) \\ r \circ f(a_{1},\ldots,a_{n}) \end{bmatrix}/\Delta_{\alpha 1}   \\
&= f^{B} \left( s(a_{1}),\ldots,s(a_{n}) \right) - s(f^{A}(a_{1},\ldots,a_{n}))  
\end{align*}
Then for a homomorphism $\phi: B^{r} \rightarrow E^{r}$ we see for $f \in \tau$ with $n = \ar f$, $\vec{a} \in A^{n}$,
\begin{align*}
\left( \check{\sigma} \circ \delta ( \phi ) \right)_{f}(\vec{a}) = \phi \circ T_{f}(\pi(a_{1}),\ldots,\pi(a_{n})) = f^{E} \left( \phi \circ s(a_{1}), \ldots, \phi \circ s(a_{n})  \right) - \phi \circ s (f^{A}(a_{1},\ldots,a_{n})).
\end{align*}
Then $\phi \circ s : A \rightarrow E$ witnesses that $\check{\sigma} \circ \delta (\phi) = 0$ in $H^{2}_{\mathcal U}(A,E^{\tau})$.

Now assume there is a free presentation $F_{\mathcal U}(X)/\theta \approx Q$ such that $F_{\mathcal U}(X)/[\theta,1]$ has an idempotent element. We show $\im \delta = \ker \check{\sigma}$. Since the sequence is a complex, we already have $\im \delta \subseteq \ker \check{\sigma}$. Now take $[T'] \in \ker \check{\sigma}$. By Lemma~\ref{lem:100}(2c), we can assume there is idempotent $0' \in Q$ such that $T'_{f}(0',\ldots,0')=0$ for all $f \in \tau$. Now,  $[T'] \in \ker \check{\sigma}$ implies there exists $h: A \rightarrow E$ such that for $f \in \tau$ with $n=\ar f$ and $\vec{a} \in A^{n}$,
\[
T_{f}'(\pi(a_{1}),\ldots,\pi(a_{n})) = f^{E} \left( h(a_{1}),\ldots,h(a_{n}) \right) - h(f^{A}(\vec{a})).
\]
For $\vec{q} \in Q^{n}$, this yields
\[
T_{f}'(\vec{q}) = f^{E} \left( h(0,q_{1}),\ldots,h(0,q_{n}) \right) - h(T_{f}(\vec{q}),f^{Q}(\vec{q}))
\]
and for any $f \in \tau$ with $n = \ar f$, $\vec{b} \in B^{n}$, we have
\begin{align*}
0 = T_{f}' (0',\ldots,0') &= T_{f}' \left( \pi(b_{1},0'),\ldots,\pi(b_{n},0') \right) \\
&= f^{E} \left(h(b_{1},0'),\ldots,h(b_{n},0') \right) - h(f^{B}(\vec{b}) + T_{f}(0',\ldots,0'),f^{Q}(0',\ldots,0')) \\
&= f^{E} \left(h(b_{1},0'),\ldots,h(b_{n},0') \right) - h \left( f^{B}(\vec{b}),0' \right).
\end{align*}
Define $\phi: B \rightarrow E$ by $\phi(b):= h(b,0') - h(0,0')$. Then the previous equation implies $\phi$ is a homomorphism. For the difference term $m$ in $\mathcal U$, we can evaluate $T_{m}'$ in two ways:
\begin{align*}
0 = T_{m}'(0',0',f^{Q}(\vec{q})) &= T_{m}'\left( \pi(T_{f}(\vec{q}),0'),\pi(0,0'),\pi(0,f^{Q}(\vec{q})) \right) \\
&= m^{E} \left( h(T_{f}(\vec{q}),0'),h(0,0'),h(0,f^{Q}(\vec{q}))  \right) - h \left( m^{A} \left((T_{f}(\vec{q}),0'),(0,0'),(0,f^{Q}(\vec{q})) \right) \right) \\
&= h(T_{f}(\vec{q}),0') - h(0,0') + h(0,f^{Q}(\vec{q})) - h \left( T_{f}(\vec{q}) + T_{m}(0',0',f^{Q}(\vec{q})), f^{Q}(\vec{q}) \right)  \\
&= h(T_{f}(\vec{q}),0') - h(0,0') + h(0,f^{Q}(\vec{q})) - h \left( T_{f}(\vec{q}), f^{Q}(\vec{q}) \right)  
\end{align*} 
Define $s: Q \rightarrow E$ by $s(x):= h(0,x)$. Then we see that
\begin{align*}
(\phi \circ T)_{f}(\vec{q}) &= h(T_{f}(\vec{q}),0') - h(0,0') \\
&= h(T_{f}(\vec{q}),f^{Q}(\vec{q})) - h(0,f^{Q}(\vec{q})) \\
&= f^{E} \left(h(0,q_{1}),\ldots,h(0,q_{n}) \right) - T_{f}'(\vec{q}) - h(0,f^{Q}(\vec{q})) \\
&= f^{E} \left( s(q_{1}),\ldots,s(q_{n}) \right) - T_{f}'(\vec{q})  - s(f^{Q}(\vec{q})) .
\end{align*}
Since $\delta(\phi) = [\phi \circ T]$, we have shown $\delta(\phi) = [T']$.
\end{proof}

In the following, we present a simple example which shows that in general some condition on idempotence is necessary for the sequence in Eq~\eqref{eq:hochserre} to be exact at all positions.

\begin{example}
We write $\mathds{Z}_{m}$ for the finite cyclic group of order $m$ which has operations signature $\{+,-,0\}$ for addition, inverse and identity, respectively. We use a slightly informal notation to consider algebras which are finite cyclic groups expanded by a single unary operation $g(x)$. Define $E := \left\langle \mathds{Z}_{m}, g(x) \right\rangle$ and $Q := \left\langle \mathds{Z}_{k}, g(x) \right\rangle$ where $g^{E}(x) = g^{Q}(x) = x$. Define $B := \left\langle \mathds{Z}_{n}, g(x) \right\rangle$ where $g^{B}(x) = x+1$ and $\mathrm{gcd}(n,m) = 1$ with $m > n$. Then we can write the direct product $A = B \times Q \approx B' \otimes^{T} Q$ for the 2-cocycle $T$ such that $T_{g}(x) = -1$ and is trivial for the group operations, and $B' := \left\langle \mathds{Z}_{n}, g(x) \right\rangle$ with $g^{B}(x) = x$.

Define a 2-cocycle $S$ appropriate to $(Q,E)$ by $S_{g}(x) = -1$ for all $x \in Q$ and is trivial for the group operations. If $[S] = 0$, then there exists $r: Q \rightarrow E$ such that $S_{f}(\vec{x}) = g^{E}(r(\vec{x})) - r(f^{Q}(\vec{x}))$ for any fundamental operation $f$; however, $-1 = S_{g}(x) = g^{E}(r(x)) - r(g^{Q}(x)) = r(x) - r(x) = 0$, a contradiction. It must be that $[S] \neq 0$. Note $\Hom \left( B', E \right) = \Hom \left( \mathds{Z}_{n}, \mathds{Z}_{m} \right) = 0$ since $g$ interprets as the identity in both algebras and $\mathrm{gcd}(n,m) = 1$; thus, $\im \delta_{T} = 0$. We now show $[S] \in \ker \check{\sigma}$.

First note in the algebra $A$ we have $g^{A}( \left\langle b, x \right\rangle ) = \left\langle g^{B'}(b) + T_{g}(x) , g^{Q}(x) \right\rangle = \left\langle b + 1 , x\right\rangle$. Define $h : A \rightarrow E$ by $h(b,x) = b$. Then for the fundamental operations we have
\begin{align*}
S_{g} \circ \pi(b,x) = S_{g}(x) = -1 = b - (b+1) = h(b,x) - h(b+1,x) = g^{E}(h(b,x)) - h(g^{A}(b,x)).
\end{align*}
\begin{align*}
S_{+} \circ \pi \left( \left\langle b, x \right\rangle , \left\langle c, y \right\rangle \right) = S_{+}(x) = 0 = b + c - (b+c) &= h(b,x) + h(c,y) - h(b+c, x+y) \\
&= h(b,x) + h(c,y) - h(\left\langle b, x \right\rangle + \left\langle c, y \right\rangle)
\end{align*}
\begin{align*}
S_{-} \circ \pi \left( \left\langle b, x \right\rangle \right) = S_{-}(x) = 0 = -b - (-b) = - h(b,x) - h (\left\langle -b, -x \right\rangle ) = - h(b,x) - h (-\left\langle b, x \right\rangle ) \\ 
\end{align*}
and
\begin{align*}
S_{0} \circ \pi \left( \left\langle b, x \right\rangle \right) = S_{0}(x) = 0 = 0^{E}(b) - h ( \left\langle 0,0 \right\rangle) = 0^{E}(h(b,x)) - h(0^{A}(b,x)).
\end{align*}
This shows that $h$ witnesses $\check{\sigma}([S]) = 0$; altogether, we have shown $\im \delta_{T} \neq \ker \check{\sigma}$.
\end{example}

The action of the inflation map on second-cohomology classes can be given an alternative characterization.

\begin{lemma}\label{lem:30}
Let $A_{1}, A_{2} \in \mathcal U$ a variety with a difference term and suppose $\pi: Q_{1} \rightarrow Q_{2}$ is surjective. Let 
\[
\check{\sigma} : H^{2}_{\mathcal U}(Q_{2},B^{\tau}) \longrightarrow H^{2}_{\mathcal U}(Q_{1},B^{\tau})
\]
be the induced inflation homomorphism. Suppose $\rho_{1}: A_{1} \rightarrow Q_{1}$ realizes central datum $(B^{\tau},Q_{1})$ and is represented by the cohomology class $[T^{1}]$ and suppose $\rho_{2}: A_{2} \rightarrow Q_{2}$ realizes central datum $(B^{\tau},Q_{2})$ and is represented by the cohomology class $[T^{2}]$. Then $\check{\sigma}([T^{2}]) = [T^{1}]$ if and only if there is a homomorphism $\phi: A_{1} \rightarrow A_{2}$ such that
\begin{enumerate}

	\item $\begin{bmatrix} r(a) \\ a \end{bmatrix}/\Delta_{\ker \rho_{1} 1} = \begin{bmatrix} r'(\phi(a)) \\ \phi(a) \end{bmatrix}/\Delta_{\ker \rho_{2} 1}$ for any $\ker \rho_{1}$-trace $r$ and $\ker \rho_{2}$-trace $r'$;
	
	\item $\pi \circ \rho_{1} = \rho_{2} \circ \phi$.

\end{enumerate}
\end{lemma}
\begin{proof}
Suppose there is such a homomorphism $\phi:A_{1} \rightarrow A_{2}$. We have isomorphisms $\psi_{1} : A_{1} \rightarrow B^{\tau} \otimes^{T_{1}} Q_{1}$ and $\psi_{2} : A_{2} \rightarrow B^{\tau} \otimes^{T_{2}} Q_{2}$ with $A_{2}(\ker \rho_{2})/\Delta_{\ker \rho_{2} 1} \approx B^{\tau} \approx A_{1}(\ker \rho_{1})/\Delta_{\ker \rho_{1} 1}$. Then $\psi_{2} \circ \phi \circ \psi_{1}^{-1}\left( \begin{bmatrix} r(a) \\ a \end{bmatrix}/\Delta_{\ker \rho_{1} 1} , \rho_{1}(a) \right) = \left( \begin{bmatrix} r'(\phi(a)) \\ \phi(a) \end{bmatrix}/\Delta_{\ker \rho_{2} 1}, \rho_{2}(\phi(a)) \right)$. By replacing $\phi$ with $\psi_{2} \circ \phi \circ \psi_{1}^{-1}$ and relabeling, we can write the action of the homomorphism as $\phi(b,q) = (\alpha(b,q),\beta(b,q)) \in B \times Q_{2}$. Conditions (1) and (2) then imply $\phi(b,q) = (b,\pi(q))$. 

Since $\phi$ is a homomorphism, for $f \in \tau$ with $n=\ar f$,
\begin{align*}
\left( f^{B}(\vec{b}) + T^{1}_{f}(\vec{q}), \pi(f^{Q_{1}}(\vec{q}))  \right) &= \phi \left( f^{B}(\vec{b}) + T^{1}_{f}(\vec{q}), f^{Q_{1}}(\vec{q})   \right)  \\
&= \phi \left( f^{A_{1}} \left( (b_{1},q_{1}),\ldots,(b_{n},q_{n}) \right) \right) \\
&= f^{A_{2}} \left( \phi(b_{1},q_{1}),\ldots,\phi(b_{n},q_{n}) \right) \\
&= f^{A_{2}} \left( (b_{1},\pi(q_{1})),\ldots,(b_{n},\pi(q_{n}))  \right) \\
&= \left( f^{B}(\vec{b}) + T_{f}^{2}(\pi(q_{1}),\ldots,\pi(q_{1})), f^{Q_{2}} (\pi(q_{1}),\ldots,\pi(q_{n}) ) \right)
\end{align*}  
from which we can conclude $T^{1}_{f}(\vec{q}) = T_{f}^{2}(\pi(q_{1}),\ldots,\pi(q_{1}))$; thus, $\check{\sigma}([T^{2}]) = [T^{1}]$. 

Conversely, suppose $\check{\sigma}([T^{2}]) = [T^{1}]$. Then there exists $h: Q_{1} \rightarrow B$ such that for $f \in \tau$ with $n =\ar f$, 
\begin{align}\label{eq:36}
T^{2}_{f}(\pi(q_{1}),\ldots,\pi(q_{n})) - T^{1}_{f}(\vec{q}) = f^{B}(h(q_{1}),\ldots,h(q_{n}) ) - h(f^{Q_{1}}(\vec{q}))
\end{align}
for $\vec{q} \in Q_{1}^{n}$. Set $T_{f}'(\vec{q}):= T^{1}_{f}(\vec{q}) + f^{B}\left( h(q_{1}),\ldots,h(q_{n}) \right) - h(f^{Q}(\vec{q}))$ for $f \in \tau$ with $n = \ar f$ and $\vec{q} \in Q_{1}^{n}$. Then $T'$ and $T^{1}$ are equivalent. By utilizing the isomorphism $A_{1} \approx B^{\tau} \otimes^{T'} Q_{1}$ we may assume $A_{1}$ is defined by $T'$. Then Eq.~\ref{eq:36} implies the map $\phi : A_{1} \rightarrow A_{2}$ defined by $\phi(b,q):=(b,\pi(q))$ is a homomorphism such that $\pi \circ \rho_{1} = \rho_{2} \circ \phi$ holds. Condition (2) holds because $\Delta_{\ker \rho_{1} 1}$-classes do not depend on the choice of trace and under the isomorphism $A_{1}(\ker \rho_{1})/\Delta_{\ker \rho_{1} 1} \approx B^{\tau}$ we have $\begin{bmatrix} (0,q) \\ (b,q) \end{bmatrix}/\Delta_{\ker \rho_{1} 1} \longmapsto b$, similarly for $\ker \rho_{2}$.
\end{proof}

We write $\kappa_{\alpha \beta}: A(\alpha) \rightarrow A(\alpha)/\Delta_{\alpha \beta}$ for the canonical epimorphism for the quotient by $\Delta_{\alpha \beta}$.

\begin{lemma}\label{lem:31}
Let $\mathcal V$ be a variety with a difference term and $A \in \mathcal V$. Let $\alpha \in \Con A$ be central. There is an exact sequence 
\begin{align}\label{eqn:splitexact}
0 \longrightarrow \kappa_{\alpha 1} \left( \alpha \wedge [1,1] \right) \longrightarrow A(\alpha)/\Delta_{\alpha 1} \longrightarrow A/[1,1]\left( \beta \right)/\Delta_{\beta 1} \longrightarrow 0
\end{align}
where $\beta := (\alpha \vee [1,1])/[1,1]$. Let $F$ be a free algebra in $\mathcal V$ and $\theta \in \Con F$. If we set $A:= F/[\theta,1]$ and $\alpha:= \theta/[\theta,1]$, then the sequence in Eq~\eqref{eqn:splitexact} is split-exact.
\end{lemma}
\begin{proof}
By \cite[Cor 2.7]{wiresI}, we can identify $A$ by the isomorphism $A \approx B^\tau \otimes^{T} Q$ with $B^{\tau} \approx A(\alpha)/\Delta_{\alpha 1}$ and $Q \approx A/\alpha$. Consider $(a,c) \in \alpha \vee [1,1]$. Then there is a sequence $a=a_{1},a_{2},\ldots,a_{n}=c$ such that $(a_{i},a_{i+1}) \in \alpha$ for $i$ odd and $(a_{i},a_{i+1}) \in [1,1]$ for $i$ even. Then we can write $a_{i} = (b_{i},t_{i})$ with $t_{i}=t_{i+1}$ for $i$ odd. For $1 \leq i \leq n-2$ odd, by applying the difference term we see that   
\[
\begin{bmatrix} a_{i+2}/[1,1] \\ a_{i}/[1,1] \end{bmatrix} =   \begin{bmatrix} m \left( a_{i+1},a_{i+1},a_{i+2} \right)/[1,1]  \\ m \left( a_{i},a_{i+1},a_{i+1} \right)/[1,1] \end{bmatrix} \mathrel{\Delta_{\beta 1}} \begin{bmatrix} m \left( a_{i+2}, a_{i+2}, (0,t_{i}) \right)/[1,1] \\ m \left( a_{i}, a_{i+1}, (0,t_{i}) \right)/[1,1] \end{bmatrix} = \begin{bmatrix} (0,t_{i})/[1,1] \\ (b_{i} - b_{i+1},t_{i})/[1,1] \end{bmatrix}
\]
since $a_{i+1}/[1,1] = a_{i+2}/[1,1]$. Since $\alpha \geq \alpha \wedge [1,1]$, we also have 
\[
\begin{bmatrix} a_{i+2}/[1,1] \\ a_{i}/[1,1] \end{bmatrix} \mathrel{\Delta_{\beta 1}} \begin{bmatrix} (0,t_{i})/[1,1] \\ (b_{i} - b_{i+1},t_{i})/[1,1] \end{bmatrix} \mathrel{\Delta_{\beta 1}} \begin{bmatrix} (0,t)/[1,1] \\ (b_{i} - b_{i+1},t)/[1,1] \end{bmatrix}
\]
for all $t \in Q$ and $i$ odd. Inductively, for $i+2 \leq j \leq n-2$ odd and all $t \in Q$ we conclude
\begin{align*}
\begin{bmatrix} a_{j+2}/[1,1] \\ a_{i}/[1,1] \end{bmatrix} =   \begin{bmatrix} m \left( a_{j},a_{j},a_{j+2} \right)/[1,1]  \\ m \left( a_{i},a_{j},a_{j} \right)/[1,1] \end{bmatrix} &\mathrel{\Delta_{\beta 1}} \begin{bmatrix} m \left( (0,t), (0,t), (0,t) \right)/[1,1] \\ m \left( (b_{i} - b_{j},t), (0,t), (b_{j} - b_{j+1},t) \right)/[1,1] \end{bmatrix} \\
&= \ \ \ \begin{bmatrix} (0,t)/[1,1] \\ (b_{i} - b_{j+1},t)/[1,1] \end{bmatrix}.
\end{align*}
Altogether, 
\[
\begin{bmatrix} c/[1,1] \\ a/[1,1]  \end{bmatrix}/\Delta_{\beta 1} = \begin{bmatrix} (b_{n},t_{n})/[1,1] \\ (b_{1},t_{1})/[1,1]  \end{bmatrix}/\Delta_{\beta 1} = \begin{bmatrix} (0,t)/[1,1] \\ (b_{1} - b_{n},t)/[1,1]  \end{bmatrix}/\Delta_{\beta 1}.
\]
We also see that $\begin{bmatrix} (c,t)/[1,1] \\ (a,p)/[1,1]  \end{bmatrix}/\Delta_{\beta 1}  = \begin{bmatrix} (c,t)/[1,1] \\ (b,q)/[1,1]  \end{bmatrix}/\Delta_{\beta 1}$ implies 
\[
\begin{bmatrix} (0,t)/[1,1] \\ (a-c,t)/[1,1]  \end{bmatrix}/\Delta_{\beta 1} =  \begin{bmatrix} (c,t)/[1,1] \\ (a,p)/[1,1]  \end{bmatrix}/\Delta_{\beta 1}  = \begin{bmatrix} (c,t)/[1,1] \\ (b,q)/[1,1]  \end{bmatrix}/\Delta_{\beta 1} = \begin{bmatrix} (0,t)/[1,1] \\ (b-c,t)/[1,1]  \end{bmatrix}/\Delta_{\beta 1}
\]
and so $\left\langle (a-c,t)/[1,1],(b-c,t)/[1,1] \right\rangle \in [\beta,1]=0$. This implies $\left\langle (a-c,t),(b-c,t) \right\rangle \in \alpha \wedge [1,1]$ and so $a-b + \kappa_{\alpha 1} \left( \alpha \wedge [1,1] \right) = \kappa_{\alpha 1} \left( \alpha \wedge [1,1] \right)$. With the above, it is now easy to see that $\xi : \begin{bmatrix} (0,t) \\ (b,t) \end{bmatrix}/\Delta_{\alpha 1} \longmapsto \begin{bmatrix} (0,t)/[1,1] \\ (b,t)/[1,1] \end{bmatrix}/\Delta_{\beta 1}$ is a surjective homomorphism with kernel $\kappa_{\alpha 1} \left( \alpha \wedge [1,1] \right)$.

Now assume $F=F_{\mathcal V}(X)$ is an algebra in $\mathcal V$ freely generated by $X$ and take $\theta \in \Con F$ with $Q = F/\theta$. To ease notation, set $F' := F/[\theta,1], \theta' := \theta/[\theta,1]$, and $\beta:= \big( \theta \vee [1,1] \big)/[1,1] \in \Con F/[1,1]$ and $\beta':= \big( \theta' \vee [1',1'] \big)/[1',1'] \in \Con F'/[1',1']$ with the total congruences $1'=1_{F'}$ and $1''=1_{F'/[1',1']}$. We have the canonical surjections 
\begin{align*}
&\pi : F \rightarrow Q, &\pi' : F/[\theta,1] \rightarrow Q  \quad \quad \text{ and } \quad \quad &\pi'' : F/[1,1] \rightarrow F/(\theta \vee [1,1]).
\end{align*}
for the congruences $\theta$, $\theta'$ and $\beta$, respectively. By the $3^{\mathrm{rd}}$-isomorphism theorem, there are homomorphisms $\gamma: F/[\theta,1] \rightarrow F/[1,1]$ given by $\gamma(x/[\theta,1]) = x/[1,1]$ and $\rho: Q \rightarrow F/(\theta \vee [1,1])$ given by $\rho(x/\theta) = x/(\theta \vee [1,1])$ such that $\rho \circ \pi' = \pi'' \circ \gamma$. If we fix $l: Q \rightarrow F$ a lifting of $\pi$, then $l'(x):=l(x)/[\theta,1]$ is a lifting of $\pi'$. If $q: F/(\theta \vee [1,1]) \rightarrow Q$ is a lifting of $\rho$, then $l'' := \gamma \circ l' \circ q$ will be a lifting of $\pi''$.

We saw in Lemma~\ref{lem:100} that $\theta' \in \Con F'$ is central; similarly, $\beta$ is central in $F/[1,1]$. The lifting $l'$ determines a representation of the central extension 
\begin{align}
F' &\approx F'(\theta')/\Delta_{\theta' 1} \otimes^{T} Q &x/[\theta,1] \mapsto \left\langle b_{x}, x/\theta \right\rangle
\end{align}
with the identification $b_{x} = \begin{bmatrix} l' \circ \pi(x) \\ x/[\theta,1] \end{bmatrix}/\Delta_{\theta' 1'} = \begin{bmatrix} l(x/\theta)/[\theta,1] \\ x/[\theta,1] \end{bmatrix}/\Delta_{\theta' 1'} = \begin{bmatrix} \left\langle 0 , x/\theta \right\rangle \\ \left\langle b_{x} , x/\theta \right\rangle \end{bmatrix}/\Delta_{\theta' 1'}$; consequently, we note that $l'(x/\theta) \mapsto \left\langle 0 , x/\theta \right\rangle$ for all $x \in F$ which implies $b_{l(q)} = 0$ for all $q \in Q$. We also have the representation of the central extension 
\begin{align}\label{eqn:rep13}
F/[1,1] &= F/[1,1](\beta)/\Delta_{\beta 1} \otimes^{S} F/(\theta \vee [1,1]) &x/[1,1] \mapsto \left\langle c_{x}, x/(\theta \vee [1,1])\right\rangle
\end{align} 
determined by the lifting $l''$ with the identification $c_{x} = \begin{bmatrix} l'' \circ \pi''(x/[1,1]) \\ x/[1,1] \end{bmatrix}/\Delta_{\beta 1}$.

Since $F$ is freely generated, there is a homomorphism $\sigma : F \rightarrow F'(\theta')/\Delta_{\theta' 1}$ induced by the mapping on the generators $X \ni x \longmapsto b_{x}$. Since $F'(\theta')/\Delta_{\theta' 1}$ is an abelian algebra, there is an induced homomorphism $\check{\sigma} : F/[1,1] \rightarrow F'(\theta')/\Delta_{\theta' 1}$ such that $\check{\sigma}(\gamma(x/[\theta,1])) = \check{\sigma}( x/[1,1] ) = \sigma(x)$. From this we see that 
\begin{align}\label{eqn:eq33}
\check{\sigma}( \gamma \circ l'( x ) ) &= \check{\sigma}( \gamma \left( l( x )/[\theta,1] \right) ) = \sigma(l(x)) = b_{l(x)}=0 &(x \in Q).
\end{align}
The representation in Eq~\eqref{eqn:rep13} yields the restricted homomorphism $r(\check{\sigma}) : F/[1,1](\beta)/\Delta_{\beta 1} \rightarrow  F'(\theta')/\Delta_{\theta' 1}$ defined by $r(\check{\sigma})(c) = \check{\sigma}(c,t) - \check{\sigma}(0,t)$ for $c \in F/[1,1](\beta)/\Delta_{\beta 1}$ and any choice of $t \in F/(\theta \vee [1,1])$. It was shown in Section~\ref{sec:3}, that the definition of $r(\check{\sigma})$ is a homomorphism independent of the choice of $t \in F/(\theta \vee [1,1])$ and so is well-defined. If we have $x/[1,1] \mapsto \left\langle c_{x} , t \right\rangle$ in the representation Eq~\eqref{eqn:rep13}, then we can calculate
\begin{align*}
r(\check{\sigma})(c) = \check{\sigma}(c,t) - \check{\sigma}(0,t) &= \check{\sigma} \left( x/[1,1] \right) - \check{\sigma} \left( l'' \circ \pi''(x/[1,1]) \right) \\
&= \check{\sigma} \left( \gamma(x/[\theta,1]) \right) - \check{\sigma} \left( \gamma \circ l' \circ q \circ \pi'' \circ \gamma (x/[1,1]) \right) \\
&= \check{\sigma} \left( \gamma(x/[\theta,1]) \right) = \sigma(x)
\end{align*}
by Eq~\eqref{eqn:eq33}.

There is a homomorphism $\psi': F'/[1',1']\left( \beta' \right)/\Delta_{\beta' 1''} \rightarrow F/[1,1]\left( \beta \right)/\Delta_{\beta 1}$ induced by the third isomorphism theorem $\psi: F'/[1',1'] \rightarrow F/[1,1]$ with $\psi(\beta') = \beta$. Let us note $l'' \circ \pi''(x/[1,1]) = l'' \circ \pi'' \circ \gamma (x/[\theta,1]) = \gamma \circ l' \circ q \circ p \circ \pi' (x/[\theta,1]) = (l \circ q \circ p(x/\theta))/[1,1]$. We then calculate
\begin{align*}
\xi \circ r(\check{\sigma}) \circ \psi' \left( \begin{bmatrix} \left\langle 0,x/\theta \right\rangle/[1',1'] \\ \left\langle b_{x},x/\theta \right\rangle/[1',1']  \end{bmatrix}/\Delta_{\beta' 1''} \right) &= \xi \circ r(\check{\sigma}) \circ \psi' \left( \begin{bmatrix} \big( l(x/\theta)/[\theta,1] \big)/[1',1'] \\ \big( x/[\theta,1] \big)/[1',1']  \end{bmatrix}/\Delta_{\beta' 1''} \right) \displaybreak[0]\\
&= \xi \circ r(\check{\sigma}) \left( \begin{bmatrix} l(x/\theta)/[1,1] \\  x/[1,1]  \end{bmatrix}/\Delta_{\beta 1} \right) \displaybreak[0]\\
&= \xi \circ r(\check{\sigma}) \left( \begin{bmatrix} (l \circ q \circ p(x/\theta))/[1,1] \\ m \big( x, l(x/\theta), l \circ q \circ p(x/\theta) \big)/[1,1] \end{bmatrix}/\Delta_{\beta 1} \right) \displaybreak[0]\\
&= \xi \circ \check{\sigma} \Big( m \big( x, l(x/\theta), l \circ q \circ p(x/\theta) \big)/[1,1] \Big) \displaybreak[0]\\
&= \xi \left( \begin{bmatrix} (l \circ q \circ p(x/\theta)) \big)/[\theta,1] \\  m \big( x, l(x/\theta), l \circ q \circ p(x/\theta) \big)/[\theta,1] \end{bmatrix}/\Delta_{\theta' 1'} \right) \displaybreak[0]\\
&= \xi \left( \begin{bmatrix} \left\langle 0 , q \circ p(x/\theta) \right\rangle  \\  \left\langle b_{x} - b_{l(x/\theta)} + b_{l \circ q \circ p(x/\theta)} , q \circ p(x/\theta) \right\rangle \end{bmatrix}/\Delta_{\theta' 1'} \right) \displaybreak[0]\\
&= \xi \left( \begin{bmatrix} \left\langle 0 , q \circ p(x/\theta) \right\rangle \\ \left\langle b_{x} , q \circ p(x/\theta) \right\rangle \end{bmatrix}/\Delta_{\theta' 1'} \right) \displaybreak[0]\\
&= \xi \left( \begin{bmatrix} \left\langle 0 , x/\theta \right\rangle \\ \left\langle b_{x} , x/\theta \right\rangle \end{bmatrix}/\Delta_{\theta' 1'} \right) \displaybreak[0]\\ 
&= \begin{bmatrix} \left\langle 0 , x/\theta \right\rangle/[1',1'] \\ \left\langle b_{x} , x/\theta \right\rangle/[1',1'] \end{bmatrix}/\Delta_{\beta' 1''};
\end{align*}
therefore, $\xi \circ ( \chi \circ \psi') = \id$ which shows the sequence in Eq~\eqref{eqn:splitexact} is split-exact. 
\end{proof}

For varieties $\mathcal U$ with a difference term, the abelian extensions form a subgroup $\mathrm{Ext}_{\mathcal U}(Q, B^{\tau}) \leq H^{2}_{\mathcal U}(Q, B^{\tau})$ of second-cohomology according to \cite[Cor 3.45]{wiresI}. For the canonical homomorphism $\pi: Q \rightarrow Q/[1,1]$, the inflation homomorphism $\check{\sigma}: H^{2}_{\mathcal U}(Q/[1,1], B^{\tau}) \rightarrow H^{2}_{\mathcal U}(Q, B^{\tau})$ on cohomology can be restricted to the abelian extensions. The image of the inflation map on the subgroup of abelian extensions is characterized in the following lemma.

\begin{lemma}\label{lem:32}
Let $\mathcal U$ be a variety with a difference term and $(Q,B^{\tau})$ central datum in $\mathcal U$. For the restriction of the inflation map
\[
\check{\sigma}: \mathrm{Ext}_{\mathcal U}(Q/[1,1], B^{\tau}) \rightarrow H^{2}_{\mathcal U}(Q, B^{\tau}),
\]
$[T] \in \im \check{\sigma}$ if and only if $[T]$ is the cohomology class for an extension $\rho: A \rightarrow Q$ realizing $(Q,B^{\tau})$ such that $\ker \rho \wedge [1_{A},1_{A}] = 0_{A}$.
\end{lemma}
\begin{proof}
Suppose $[T] \in \im \check{\sigma}$. Then $T$ corresponds to a central extension $\rho_{1}: A_{1} \rightarrow Q$ and there exists a compatible 2-cocycle $S$ corresponding to central extension $\rho_{2}: A_{2} \rightarrow Q/[1,1]$ such that $[S \circ \pi] = \check{\sigma}([S]) = [T]$. By Lemma~\ref{lem:30}, there is a homomorphism $\phi: A_{1} \rightarrow A_{2}$ such that 
\begin{itemize}

	\item $\rho_{2} \circ \phi = \pi \circ \rho_{1}$, and
	
	\item $\begin{bmatrix} r(a) \\ a \end{bmatrix}/\Delta_{\ker \rho_{1} 1} = \begin{bmatrix} r'(\phi(a)) \\ \phi(a) \end{bmatrix}/\Delta_{\ker \rho_{2} 1}$ for all $\ker \rho_{1}$-traces $r$ and $\ker \rho_{2}$-traces $r'$.

\end{itemize}
Since $A_{2}$ is abelian, we have $[1,1] \leq \phi$. Take $(a,b) \in \ker \rho_{1} \wedge [1,1]$ in $A_{1}$. Then $(a,b) \leq [1,1] \leq \phi$ implies $\phi(a)=\phi(b)$ and $(a,b) \in \ker \rho_{1}$ implies $r(a)=r(b)$. We then have
\[
\begin{bmatrix} r(a) \\ a \end{bmatrix}/\Delta_{\ker \rho_{1} 1} = \begin{bmatrix} r'(\phi(a)) \\ \phi(a) \end{bmatrix}/\Delta_{\ker \rho_{2} 1} = \begin{bmatrix} r'(\phi(b)) \\ \phi(b) \end{bmatrix}/\Delta_{\ker \rho_{2} 1} = \begin{bmatrix} r(b) \\ b \end{bmatrix}/\Delta_{\ker \rho_{1} 1} 
\]  
which implies $a=b$ since $\ker \rho_{1}$ is central. We have $\ker \rho_{1} \wedge [1,1]=0$.

Now take a compatible 2-cocycle $T$ corresponding to a central extension $\rho: A \rightarrow Q$ realizing $(B^{\tau},Q)$ with $\ker \rho \mathrel{\wedge} [1,1] = 0$. Write $\alpha = \ker \rho$. By the homomorphism property in varieties with a difference term, we see that $[1/\alpha,1/\alpha] = [1,1] \vee \alpha/\alpha$ which implies $Q/[1,1] \approx A/([1,1] \vee \alpha)$. Then there is a central extension $\gamma: A/[1,1] \rightarrow Q/[1,1]$ with $\ker \gamma = ([1,1] \vee \alpha)/[1,1]$. By Lemma~\ref{lem:31}, $B^{\tau} \approx A(\alpha)/\Delta_{\alpha 1} \approx A/[1,1](\ker \gamma)/\Delta_{\ker \gamma 1}$. Then we have $A/[1,1] \approx B^{\tau} \otimes^{S} Q/[1,1]$ for some compatible 2-cocycle $S$. If we take the canonical homomorphism $\pi': A \rightarrow A/[1,1]$, then $\pi \circ \rho = \gamma \circ \pi'$ and it is easy to see that condition (2) in Lemma~\ref{lem:30} is satisfied; thus, $[T] = \check{\sigma}([S])$.
\end{proof}

\vspace{0.3cm}


\section{Schur's theorem in varieties with a difference term}\label{sec:5}
\vspace{0.3cm}

We begin with the definition of a regular kernel algebra.

\begin{definition}
A kernel algebra $E^{\tau}$ is \emph{regular} in the variety $\mathcal V$ if the following properties hold:
\begin{enumerate}

	\item[(RD1)] (injectivity) If $B^{\tau},C^{\tau} \in \mathcal V$ are kernel algebras such that $B^{\tau} \leq C^{\tau}$ and there is a homomorphism $\phi: B^{\tau} \rightarrow E^{\tau}$, then there is a homomorphism $\psi: C^{\tau} \rightarrow E^{\tau}$ such that $\psi|_{B} = \phi$.
	
	\item[(RD2)] (separation) If $B^{\tau} \in \mathcal V$ a kernel algebra and $(a,b) \not\in \alpha \in \Con B^{\tau}$, then there is a homomorphism $\phi: B^{\tau} \rightarrow E^{\tau}$ such that $\alpha \leq \ker \phi$ and $(a,b) \not\in \ker \phi$. 

\end{enumerate}
Central datum $(Q,E^{\tau})$ is \emph{regular datum} for the variety $\mathcal V$ if the kernel algebra is regular in $\mathcal V$.
\end{definition}

\begin{example}
$\mathds{C}^{\times}$ is regular for the variety $\mathcal G$ of groups. Since $\mathds{C}^{\times}$ is divisible, it has the injective mapping property for abelian groups (or $\mathds{Z}$-modules) and so condition (RD1) is satisfied. Since $\mathds{Z}, \mathds{Q}/\mathds{Z} \leq \mathds{C}^{\times}$, for any subgroup $S \leq G$ of an abelian group and $a \not\in S$, there is a homomorphism $\phi :G \rightarrow \mathds{C}^{\times}$ such that $\phi(S)=0$ and $\phi(a) \neq 0$ (\cite[Thm 10.58]{rotman}). This is the (RD2) condition.
\end{example}

\begin{example}\label{ex:multiext}
Let $\mathcal V$ be variety of modules over a division ring $D$ expanded by multilinear operations (see \cite{wiresII}) whose arities are at least binary; for example, $\mathcal V$ can be a variety of Leibniz algebras, dendriform algebras or diassociative algebras over $D$. If the multilinear operations are interpreted as trivial operations in $D$, then $D$ is regular in $\mathcal V$.

The signature $\tau$ of $\mathcal V$ is the union of the module operations with the additional multilinear operation indexed by $F$. For an algebra $M \in \mathcal V$, ideals $I \triangleleft M$ are in bijective correspondence with congruences by $\alpha_{I} = \{(a,b) \in M^{2}: a-b \in I \}$. There is an isomorphism 
\begin{align}\label{eqn:1000}
\psi : M(\alpha_{I})/\Delta_{\alpha_{I} 1} \rightarrow I/[I,M]
\end{align}
given by $\psi \left( \begin{bmatrix} b \\ a \end{bmatrix}/\Delta_{\alpha_{I} 1} \right) = (a-b)/[I,M]$; therefore, for a central ideal $I \triangleleft M$, there is an isomorphism $M(\alpha_{I})/\Delta_{\alpha_{I} 1} \approx I$.  We also see that for an abelian ideal $I \triangleleft M$, the multilinear operations which are at least binary are all trivial in $I$; therefore, the kernel algebras in $\mathcal V$ are $D$-modules. Given a submodule $S \leq I$, we can choose a basis $B_{I}$ for $S$ which can be extended to a basis $B_{I} \cup B'$ for $I$. With the aid of the basis, it is straightforward to construct the linear maps which verify conditions (RD1) and (RD2).
\end{example}

We next establish the existence of regular datum in a variety. We utilize the fact that any unital module over a ring with identity embeds in an injective module.

\begin{proposition}\label{prop:regexists}
If $\mathcal V$ is a variety with a difference term which contains a nontrivial abelian algebra, then $\mathcal V$ contains a kernel algebra which is regular for $\mathcal V$.
\end{proposition}
\begin{proof}
Let $\mathcal A$ be the subvariety generated by the abelian algebras in $\mathcal V$. Note the kernel algebras of $\mathcal V$ are just the algebras in $\mathcal A$ which have an idempotent element. Every algebra in $\mathcal A$ is polynomially equivalent to a unital module over a fixed ring with unity given by $R = \{t(x,y) \in F_{\mathcal A}(x,y): t(y,y)=y \}$. Take the 1-generated free algebra $J = F_{\mathcal A}(u)$ and its associated module $\tilde{J} = M(J,\id_{J})$. We have that the coslice categories $(J \downarrow \mathcal A)$ and $(\tilde{J} \downarrow \mathcal M_{R})$ are term-equivalent (a more general development is expounded in \cite[Cha 9]{commod} but we require only the more streamlined version found in \cite[Thm 7.38]{bergman}) and outlined in Section~\ref{sec:2}.

For any kernel algebra $B^{\tau} \in V$, if we take the homomorphism $\mu: J \rightarrow B^{\tau}$ such that $\mu(u)$ is the idempotent element of $B^{\tau}$ which corresponds to the zero of the module $B$, then $B=M(B^{\tau},\mu)$. Since the associated modules of kernel algebras are modules over the same ring $R$, we see that $\Hom(B^{\tau},E^{\tau}) = \Hom_{R}(B,E)$. Given $(M,\mu) \in (\tilde{J} \downarrow \mathcal M_{R})$ with $M \in \mathcal M_{R}$ and homomorphism $\mu : \tilde{J} \rightarrow M$, we have the corresponding abelian algebra $A(M,\mu) \in \mathcal A$ with $(A(M,\mu), \mu) \in (J \downarrow \mathcal A)$. We note that if $\hat{0}: \tilde{J} \rightarrow M$ is the homomorphism which maps all of $\tilde{J}$ to $0 \in M$, then in the constructed algebra $A(M,\hat{0})$ we see that $0$ will be an idempotent element; that is, $A(M,\hat{0})$ is a kernel algebra in $\mathcal A \leq \mathcal V$.

There is a divisible group $D$ with an embedding $\left( \bigoplus_{\theta \in \Con \, _{R} R} R/\theta \right)^{+} \hookrightarrow D$ of the abelian group reduct (\cite[Lem 3.10]{hungerford}). Then there is an $R$-module embedding $\epsilon : \bigoplus_{\theta \in \Con \, _{R} R} R/\theta \hookrightarrow \Hom_{\mathds{Z}}(R,D)$ where $\Hom_{\mathds{Z}}(R,D)$ is an injective $R$-module (\cite[Prop 3.12]{hungerford}). The claim is that $A(\Hom_{\mathds{Z}}(R,D), \hat{0})$ is regular in $\mathcal A$. It is easy to see that the injectivity property of $\Hom_{\mathds{Z}}(R,D)$ carries over to kernel algebras so that condition (RD1) holds for $A(\Hom_{\mathds{Z}}(R,D), \hat{0})$.

To verify the (RD2) condition, take a kernel algebra $B^{\tau} \in \mathcal V$ and $(a,b) \not\in \alpha \in \Con B^{\tau}$. Since the algebra $B^{\tau}$ is polynomially equivalent to the module $B$, there is a submodule $S \leq B$ such that $\alpha = \theta_{S} = \{ (x,y) \in B^{2} : x-y \in S \} \in \Con B = \Con B^{\tau}$. Then $(a,b) \not\in \alpha$ implies $a - b \not\in S$ and so $(a-b)/S \neq 0$ in $B/S$. There is the canonical surjection $\pi_{S} : B \rightarrow B/S$ and 1-generated submodule $R((a-b)/S) \leq B/S$. Then there is $\theta \in \Con \ _{R} R$ such that $R((a-b)/S) \approx R/\theta$ and so there is an inclusion $i: R((a-b)/S) \rightarrow \bigoplus_{\theta \in \Con \, _{R} R} R/\theta$ as a factor in the direct sum. So we have an embedding $\epsilon \circ i : R((a-b)/S) \rightarrow \Hom_{\mathds{Z}}(R,D)$ which by the injectivity property implies the existence of a homomorphism $\phi : B/S \rightarrow \Hom_{\mathds{Z}}(R,D)$ such that $\phi \circ \pi_{S} (a-b) = \phi ((a-b)/S) = \epsilon \circ i((a-b)/S) \neq 0$ and $\phi \circ \pi_{S}(S) = 0$. Then we have the corresponding homomorphism $\phi \circ \pi_{S}: B^{\tau} \rightarrow A(\Hom_{\mathds{Z}}(R,D), \hat{0})$ such that $\alpha = \theta_{S} \leq \ker \phi \circ \pi_{S}$ and $\phi \circ \pi_{S}(a) \neq \phi \circ \pi_{S}(b)$.
\end{proof}

\begin{lemma}\label{lem:regularnontriv}
Let $\mathcal V$ be a variety with a difference term which contains a nontrivial abelian algebra. If $E^{\tau}$ is regular in $\mathcal V$, then for all kernel algebras in $B^{\tau} \in \mathcal V$, $B^{\tau} \neq 0$ if and only if $\Hom \left( B^{\tau} , E^{\tau} \right) \neq 0$.
\end{lemma}
\begin{proof}
Clearly, if $B^{\tau}$ is trivial, then $\Hom \left( B^{\tau} , E^{\tau} \right) = 0$. If $B^{\tau} \neq 0$, then for any $0 \neq a \in B$ there is a homomorphism $\phi: B^{\tau} \rightarrow E^{\tau}$ such that $\phi(a) \neq 0$ by the (RD2) condition of regularity; thus, $\phi$ is a nontrivial homomorphism.
\end{proof}

For a group $Q$, there appears to be two different traditions in referring to the Schur Multiplier of $Q$: (1) as $\frac{R \cap [F,F] }{[R,F]}$ where $F/R \approx Q$ is a free-presentation, (2) as the $2^{\mathrm{nd}}$-cohomology group $H^{2}(Q,\mathds{C}^{\times})$. Using the second definition, one may develop the restriction, inflation and transgression maps of cohomology and establish the relationship between the Schur Multiplier and liftings of projective representations of $Q$, while the first definition is directly related to covers, universal central extensions and commutator relations. The isomorphism between the two in the case $Q$ is finite often referred to as the Schur or Schur-Hopf formula. In the next definition, we will follow tradition (1) in referring to the Schur multiplier of an algebra in a variety with a difference term. We shall see in Theorem~\ref{thm:schurmult} and Theorem~\ref{thm:coversexist} that second-cohomology of regular datum $(Q,E^{\tau})$ is characterized by the homset of the Schur multiplier of $Q$ with the regular algebra $E^{\tau}$. The definition is complicated because of the form of the kernel algebra of a central extension associated to a free presentation; nevertheless, we prove in Theorem~\ref{thm:invSchur} that the Schur multiplier does not depend on the choice of free-presentation nor the existence of idempotents.

If $F/\theta \approx Q$ is a free-presentation, then by Lemma~\ref{lem:100} the map $\pi: F/[\theta , 1] \rightarrow Q$ is a central extension with kernel $\theta/[\theta,1]$ and associated kernel algebra $F/[\theta , 1] \left( \theta/[\theta,1] \right)/\Delta_{\theta/[\theta,1], 1}$. For the canonical epimorphism $\kappa_{ \theta/[\theta,1] 1} : F/[\theta , 1] \left( \theta/[\theta,1] \right) \longrightarrow F/[\theta , 1] \left( \theta/[\theta,1] \right)/\Delta_{\theta/[\theta,1],1}$ we then consider the pushdown of the subalgebra $F/[\theta,1] \big( \theta \wedge [1,1]/[\theta,1] \big)$.

\begin{definition}\label{def:schurmult}
Let $Q \in \mathcal V$ a variety with a difference term. If $F/\theta \approx Q$ is a free-presentation in the variety, then the associated \emph{Schur multiplier} of $Q$ is the abelian algebra $\kappa_{ \theta/[\theta,1] 1} \Big( F/[\theta,1] \big( \theta \wedge [1,1]/[\theta,1] \big) \Big)$.
\end{definition}

\begin{proposition}\label{prop:schurmult}
Let $A \in \mathcal U$ a variety with a difference term and $E^{\tau} \in \mathcal U$ a regular kernel algebra. If $\pi: A \rightarrow Q$ is a central extension such that $A$ has an idempotent element, then for the associated transgression $\delta$,
\begin{align*}
\im \delta \approx \Hom \left( \kappa_{\ker \pi 1}\left( A \left( \ker \pi \wedge [1,1] \right) \right) , E^{\tau} \right) \approx \Hom( I_{\ker \pi \wedge [1,1]} , E^{\tau} ).
\end{align*}
\end{proposition}
\begin{proof}
Let $\pi: A \rightarrow Q$ realize datum $(B^{\tau},Q)$ where $B^{\tau} = A(\alpha)/\Delta_{\alpha 1}$ with $\alpha = \ker \pi$. By Theorem~\ref{thm:inflationrestriction}, we have
\[
\im \delta \approx \Hom \left( B^{\tau}, E^\tau \right)/\ker \delta \approx \Hom \left( B^{\tau}, E^\tau \right)/\im \check{r}.
\] 
There is the restriction homomorphism
\[
\Hom (B^{\tau},E^{\tau}) \stackrel{\nabla}{\longrightarrow} \Hom \left( \kappa_{\alpha 1} \left( A(\alpha \wedge [1,1]) \right), E^{\tau} \right)
\]
defined by the subalgebra inclusion $\kappa_{\alpha 1} \left( A(\alpha \wedge [1,1]) \right) \leq A(\alpha)/\Delta_{\alpha 1} = B^{\tau}$. By regularity of $E^{\tau}$, $\nabla$ is surjective. The proposition follows once we show $\im \check{r} = \ker \nabla$.

First, take $\phi \in \im \check{r}$. We use the isomorphism $\nu: A \rightarrow B^{\tau} \otimes^{T} Q$. There is a homomorphism $\psi: A \rightarrow E^{\tau}$ such that $\phi(b) = \psi(b,p) - \psi(0,p)$ for any $p \in Q$. Since $E^{\tau}$ is abelian, $[1,1] \leq \ker \psi$. Then for $\left\langle (a,p),(b,q) \right\rangle \in \alpha \wedge [1,1]$ we must have $p=\pi(a,p) = \pi(b,q)=q$. Since $[1,1] \leq \ker \psi$, we have $\psi(a,p)=\psi(b,q)$. Then $\phi(a)=\psi(a,p) - \psi(0,p) = \psi(b,p) - \psi(0,p) = \phi(b)$.

Note, by the isomorphism we have $B \ni a \longmapsto \begin{bmatrix} (0,p) \\ (a,p) \end{bmatrix}/\Delta_{\alpha 1}$ for any $p \in Q$. Then
\begin{align*}
a + \kappa_{\alpha 1} \left( A(\alpha \wedge [1,1]) \right) = b + \kappa_{\alpha 1} \left( A(\alpha \wedge [1,1]) \right) &\Leftrightarrow a-b + \kappa_{\alpha 1} \left( A(\alpha \wedge [1,1]) \right) \\
&\Leftrightarrow  \begin{bmatrix} (0,p) \\ (a-b,p) \end{bmatrix}/\Delta_{\alpha 1} \in \kappa_{\alpha 1} \left( A(\alpha \wedge [1,1]) \right) \\
&\Leftrightarrow  \begin{bmatrix} (b,p) \\ (a,p) \end{bmatrix}/\Delta_{\alpha 1} \in \kappa_{\alpha 1} \left( A(\alpha \wedge [1,1]) \right) \\
&\Leftrightarrow \left\langle (a,p),(b,p) \right\rangle \in \alpha \wedge [1,1]
\end{align*}
Then $a \in \kappa_{\alpha 1} \left( A(\alpha \wedge [1,1]) \right) \Longrightarrow \phi(a) = \phi(0) = \psi(0,p) - \psi(0,p) = 0$. This implies $\phi \in \ker \nabla$.

Now consider $\phi \in \ker \nabla$; thus, $\phi: B^{\tau} \rightarrow E^{\tau}$ with $\kappa_{\alpha 1} \left( A(\alpha \wedge [1,1]) \right) \leq \ker \phi$. If we are given $\begin{bmatrix} (b,t)/[1,1] \\ (a,t)/[1,1] \end{bmatrix}/\Delta_{11} \in \kappa_{11}(A/[1,1](\alpha/[1,1]))$, we apply the difference term $m$ to the sequence
\[
\begin{bmatrix} (b,t)/[1,1] \\ (a,t)/[1,1] \end{bmatrix} \Delta_{11} \begin{bmatrix} (b,t)/[1,1] \\ (a,t)/[1,1] \end{bmatrix}, \begin{bmatrix} (a,t)/[1,1] \\ (a,t)/[1,1] \end{bmatrix} \Delta_{11} \begin{bmatrix} (b,t)/[1,1] \\ (b,t)/[1,1] \end{bmatrix}, \begin{bmatrix} (a,t)/[1,1] \\ (a,t)/[1,1] \end{bmatrix} \Delta_{11} \begin{bmatrix} (0,t)/[1,1] \\ (0,t)/[1,1] \end{bmatrix}
\]
which yields 
\[
\begin{bmatrix} \left( b,t \right)/[1,1] \\ \left( a,t \right)/[1,1]  \end{bmatrix} \Delta_{11} \begin{bmatrix} \left( 0,t \right)/[1,1] \\ \left( b-a,t \right)/[1,1]  \end{bmatrix}
\]
since $A/[1,1]$ is abelian. Define $\phi': \kappa_{11}(A/[1,1](\alpha/[1,1])) \rightarrow E^{\tau}$ by $\phi' \left( \begin{bmatrix} (0,t)/[1,1] \\ (a,t)/[1,1] \end{bmatrix}/\Delta_{11} \right) := \phi(a)$. Note that $\begin{bmatrix} (0,t)/[1,1] \\ (a,t)/[1,1] \end{bmatrix}/\Delta_{11} = \begin{bmatrix} (0,t)/[1,1] \\ (b,t)/[1,1] \end{bmatrix}/\Delta_{11}$ implies $(a,t)/[1,1] = (b,t)/[1,1]$; therefore, we have $\left\langle (a,t), (b,t) \right\rangle \in \alpha \wedge [1,1]$ and so $\phi(a) = \phi(b)$. This shows $\phi'$ is well-defined. We also see that
\begin{align*}
\phi' &\left( f \left( \begin{bmatrix} (0,t_{1})/[1,1] \\ (a_{1},t_{1})/[1,1] \end{bmatrix}/\Delta_{11}, \ldots, \begin{bmatrix} (0,t_{n})/[1,1] \\ (a_{n},t_{n})/[1,1] \end{bmatrix}/\Delta_{11} \right) \right) \displaybreak[0]\\
&= \phi' \left( \begin{bmatrix}  \left( T_{f}(\vec{t}), f^{Q}(\vec{t}) \right)/[1,1] \\ \left( f^{B}(\vec{a}) + T_{f}(\vec{t}), f^{Q}(\vec{t}) \right)/[1,1] \end{bmatrix}/\Delta_{11} \right) \displaybreak[0]\\
&= \phi' \left( \begin{bmatrix}  \left( 0, f^{Q}(\vec{t}) \right)/[1,1] \\ \left( f^{B}(\vec{a}), f^{Q}(\vec{t}) \right)/[1,1] \end{bmatrix}/\Delta_{11} \right) \displaybreak[0]\\
&= \phi \left( f^{B}(\vec{a}) \right) \displaybreak[0]\\
&= f^{E} \left( \phi(a_{1}),\ldots,\phi(a_{n}) \right) \displaybreak[0]\\
&= f^{E} \left( \phi'\left( \begin{bmatrix} (0,t_{1})/[1,1] \\ (a_{1},t_{1})/[1,1] \end{bmatrix}/\Delta_{11} \right),\ldots, \phi' \left( \begin{bmatrix} (0,t_{1})/[1,1] \\ (a_{1},t_{1})/[1,1] \end{bmatrix}/\Delta_{11} \right) \right)
\end{align*}
and so $\phi'$ is a homomorphism.

By \cite[Prop 2.5]{wiresI}, there is an isomorphism $\sigma: A/[1,1] \approx \left( A/[1,1] \times A/[1,1] \right)/ \Delta_{11} \otimes^{T'} \{ \mathrm{pt} \}$. Now we take $\bar{0}$ to be the idempotent element of $A$ and $r: A \rightarrow A$ the 1-trace $r(A) = \{ \bar{0}\}$. Then using idempotence we have for any $f \in \tau$,
\begin{align*}
\left\langle 0, \mathrm{pt} \right\rangle = \left\langle \begin{bmatrix} \bar{0}/[1,1] \\ \bar{0}/[1,1] \end{bmatrix}/\Delta_{11} , \mathrm{pt} \right\rangle &= \sigma \left( \bar{0}/[1,1] \right) \displaybreak[0]\\
&= \left\langle \begin{bmatrix} f(r(\bar{0}),\ldots,r(\bar{0}))/[1,1] \\ f(\bar{0},\ldots,\bar{0})/[1,1] \end{bmatrix}/\Delta_{11} + T'_{f}(\mathrm{pt},\ldots,\mathrm{pt}), f(\mathrm{pt},\ldots,\mathrm{pt}) \right\rangle \displaybreak[0]\\
&= \left\langle \begin{bmatrix} f(\bar{0},\ldots,\bar{0})/[1,1] \\ f(\bar{0},\ldots,\bar{0})/[1,1] \end{bmatrix}/\Delta_{11} + T'_{f}(\mathrm{pt},\ldots,\mathrm{pt}),\mathrm{pt} \right\rangle \displaybreak[0]\\
&= \left\langle T'_{f}(\mathrm{pt},\ldots,\mathrm{pt}), \mathrm{pt} \right\rangle
\end{align*}
which yields $T'_{f}(\mathrm{pt},\ldots,\mathrm{pt}) = 0$. This implies $A/[1,1] \approx \left( A/[1,1] \times A/[1,1] \right)/\Delta_{11}$.

By regularity of $E^{\tau}$, there is a homomorphism $\phi'': \left( A/[1,1] \times A/[1,1] \right)/\Delta_{11} \rightarrow E^{\tau}$ which extends $\phi'$ on the subalgebra $\kappa_{11} \left( A/[1,1](\alpha/[1,1]) \right) \leq \left( A/[1,1] \times A/[1,1] \right)/\Delta_{11}$. Now taking the canonical homomorphism $\rho: A \rightarrow A/[1,1]$, we can define the homomorphism $\psi: A \rightarrow E^{\tau}$ by $\psi := \phi'' \circ \sigma \circ \rho$. Let $\nu(\bar{0}) = (0',p')$. We can then verify for $b \in B$, 
\begin{align*}
\check{r}(\psi)(b) &= \psi(b,t) - \psi(0,t) \displaybreak[0]\\
&= \phi'' \left( \begin{bmatrix} (0',p')/[1,1] \\ (b,t)/[1,1] \end{bmatrix}/\Delta_{11} \right) - \phi'' \left( \begin{bmatrix} (0',p')/[1,1] \\ (0,t)/[1,1] \end{bmatrix}/\Delta_{11} \right)\displaybreak[0]\\
&= m \left( \phi'' \left( \begin{bmatrix} (0',p')/[1,1] \\ (b,t)/[1,1] \end{bmatrix}/\Delta_{11} \right), \phi'' \left( \begin{bmatrix} (0',p')/[1,1] \\ (0,t)/[1,1] \end{bmatrix}/\Delta_{11} \right), \phi'' \left( \begin{bmatrix} (0,t)/[1,1] \\ (0,t)/[1,1] \end{bmatrix}/\Delta_{11} \right) \right)  \displaybreak[0]\\
&= \phi'' \left(  \begin{bmatrix} m \left( (0',p'),(0',p'),(0,t) \right)/[1,1] \\ m \left( (b,t),(0,t),(0,t) \right)/[1,1] \end{bmatrix}/\Delta_{11} \right)  \displaybreak[0]\\
&= \phi'' \left( \begin{bmatrix} (0,t)/[1,1] \\ (b,t)/[1,1] \end{bmatrix}/\Delta_{11}  \right) \displaybreak[0]\\
&=  \phi' \left( \begin{bmatrix} (0,t)/[1,1] \\ (b,t)/[1,1] \end{bmatrix}/\Delta_{11}  \right) \displaybreak[0]\\
&= \phi(b)
\end{align*}
because $\left\langle (0,t),(b,t) \right\rangle \in \alpha$ an abelian congruence and $\begin{bmatrix} (0,t)/[1,1] \\ (b,t)/[1,1] \end{bmatrix}/\Delta_{11} \in \kappa_{11} \left( A/[1,1](\alpha/[1,1]) \right)$. We can conclude $\phi \in \im \check{r}$. 
\end{proof}

With regular datum, it is possible to give necessary and sufficient conditions for when the transgression homomorphism is injective or surjective. In Proposition~\ref{prop:surjtrans}, parts (a) and (c) will have the most relevance for the present manuscript; however, this is a decent point to introduce representative datum in Definition~\ref{def:representative} and prove its connection to liftings and surjectivity of the transgression map. It will play a role in future investigations. Any homomorphism $\phi: Q \rightarrow P$ induces an inflationary map $\check{\sigma}:H^{2}_{\mathcal U}\left( Q,E^{\tau} \right) \rightarrow H^{2}_{\mathcal U}\left( P,E^{\tau} \right)$. If the central extension $\pi: C \rightarrow P$ is represented by $[S]$, then the next definition describes the situation when there are enough inflationary maps through which $[S]$ can represent every extension in $\mathcal U$ realizing datum $(Q,E^{\tau})$.

\begin{definition}\label{def:representative}
Let $P,Q \in \mathcal U$ a variety with a difference term. Let $\rho: C \rightarrow P$ realize central datum $(P,E^{\tau})$ and be represented by the compatible 2-cocycle $[S] \in H^{2}_{\mathcal U}\left( P,E^{\tau} \right)$. The extension $\rho$ is $(Q,E^{\tau})$-\emph{representative} if for all $[R] \in H^{2}_{\mathcal U}\left( Q,E^{\tau} \right)$ there is a homomorphism $\phi: Q \rightarrow P$ such that $[R]=[S \circ \phi]$.
\end{definition}

\begin{proposition}\label{prop:surjtrans}
Let $\mathcal U$ a variety with a difference term and $(Q,B^{\tau})$ central datum in $\mathcal U$. Let $\pi: A \rightarrow Q$ an extension represented by $[T] \in H^{2}_{\mathcal U}(Q,B^{\tau})$. 
\begin{enumerate}

	\item If the associated transgression $\delta_{T} : \Hom (B^{\tau}, E^{\tau}) \longrightarrow H^{2}_{\mathcal U}(Q,E^{\tau})$ is surjective, then $\pi$ has the $\rho$-lifting property for any extension $\rho: C \rightarrow P$ in $\mathcal U$ realizing $(P,E^{\tau})$. 
	
	\item Assume $F_{\mathcal U}(X)/\theta \approx Q$ is a free-presentation in which $F_{\mathcal U}(X)/[\theta,1]$ has an idempotent element. If $\pi$ has the $\rho$-lifting property for $\rho: C \rightarrow Q$ which is $(Q,E^{\tau})$-representative, then the associated transgression $\delta_{T}$ is surjective.
	
	\item Assume $F_{\mathcal U}(X)/\theta \approx Q$ is a free-presentation in which $F_{\mathcal U}(X)/[\theta,1]$ has an idempotent element. If $\pi$ has the $\mathcal U$-lifting property for $Q$, then every transgression $\delta_{T}$ is surjective.
	
\end{enumerate}	
\end{proposition}
\begin{proof}
(1) Assume the transgression $\delta_{T}$ is surjective. Let $\rho: C \rightarrow P$ be an extension in $\mathcal U$ realizing $(E^{\tau},P)$. We make use of the identifications $A \approx B^{\tau} \otimes^{T} Q$ and $C \approx E^{\tau} \otimes^{S} P$ for some 2-cocycle $S$. Let $\phi: Q \rightarrow P$ be a homomorphism. Then $[S \circ \phi] \in H_{\mathcal U}^{2} \left(Q,E^{\tau} \right)$. Since $\delta_{T}$ is surjective, there exists $\psi \in \Hom (B^{\tau},E^{\tau})$ such that $\delta_{T}(\psi) = [\psi \circ T] = [S \circ \phi]$. This implies there is a map $h : Q \rightarrow E$ such that for all $f \in \tau$, $\vec{q} \in Q^{\ar f}$ we have 
\[
\psi(T_{f}(\vec{q})) - S_{f}(\phi(\vec{x})) = f^{E}(h(\vec{q})) - h(f^{Q}(\vec{q})).
\]
We can then use the identifications to define $\bar{\phi}:A \rightarrow C$ by $\bar{\phi}(b,q):= \left\langle \psi(b) + h(q), \phi(q) \right\rangle$. Then by calculating for $f \in \tau$ with $\ar f = n$, 
\begin{align*}
\bar{\phi} \left( F^{A}_{f} \left( \left\langle b_{1},q_{1} \right\rangle,\ldots, \left\langle b_{n},q_{n} \right\rangle \right) \right) &= \bar{\phi} \left( \left\langle f^{B}(\vec{b}) + T_{f}(\vec{q}), f^{Q}(\vec{q}) \right\rangle \right) \\
&= \left\langle  \psi(f^{B}(\vec{b})) + \psi(T_{f}(\vec{q})) + h(f^{Q}(\vec{q})), \phi(f^{Q}(\vec{q})) \right\rangle \\
&= \left\langle f^{E}(\psi(\vec{b})) + f^{E}(h(\vec{q})) + S_{f}(\phi(\vec{q})), f^{P}(\phi(\vec{q})) \right\rangle  \\
&= \left\langle f^{E}\left( \psi(b_{1}) + h(q_{1}),\ldots, \psi(b_{1}) + h(q_{1}) \right) +  S_{f}(\phi(\vec{q})), f^{P}(\phi(\vec{q})) \right\rangle  \\
&= F^{C}_{f} \left( \left\langle \psi(b_{1}) + h(q_{1}), \phi(q_{1}) \right\rangle, \ldots, \left\langle \psi(b_{n}) + h(q_{n}), \phi(q_{n}) \right\rangle \right) \\
&= F^{C}_{f}\left( \bar{\phi} \left( b_{1},q_{1} \right),\ldots,\bar{\phi} \left( b_{n},q_{n} \right) \right);
\end{align*}
therefore, $\bar{\phi}$ is a homomorphism. By definition, we see that $\rho \circ \bar{\phi} = \phi \circ \pi$ and so $\bar{\phi}$ is a lifting of $\phi$.

(2) Assume $\pi$ has the $\rho$-lifting property for $\rho:C \rightarrow P$ which is $(E^{\tau},Q)$-versatile. By Lemma~\ref{lem:100}(2c), $Q$ has an idempotent element $0'$. We show $\delta_{T}$ is surjective. Again we make use of the identifications $A \approx B^{\tau} \otimes^{T} Q$ and $C \approx E^{\tau} \otimes^{S} P$ for some 2-cocycle $S$. Take $[R] \in H^{2}_{\mathcal U}\left(Q,E^{\tau} \right)$. Then by representation, there is a homomorphism $\phi:Q \rightarrow P$ such that $[R] = [S \circ \phi]$. Then by assumption, $\phi$ can be lifted to a homomorphism $\bar{\phi}: A \rightarrow C$ such that $\rho \circ \bar{\phi} = \phi \circ \pi$. Then there is a map $\lambda: A \rightarrow E$ such that
\[
\bar{\phi}(b,q) = \left\langle \lambda(b,q),\phi(q) \right\rangle.
\] 
For any $f \in \tau$ with $n=\ar f$, $\vec{q} \in Q^{\ar f}, \vec{b} \in B^{\ar f}$, we can calculate the operations
\begin{align*}
\bar{\phi} \left( F^{A}_{f} \left( \left\langle  b_{1},q_{1} \right\rangle,\ldots,\left\langle b_{n},q_{n} \right\rangle \right) \right) &= \bar{\phi} \left( \left\langle f^{B}(\vec{b}) + T_{f}(\vec{q}), f^{Q}(\vec{q}) \right\rangle  \right)  \\
&=  \left\langle \lambda \left( f^{B}(\vec{b}) + T_{f}(\vec{q}), f^{Q}(\vec{q}) \right), \phi(f^{Q}(\vec{q})) \right\rangle
\end{align*}
and 
\begin{align*}
F^{C}_{f} \left( \bar{\phi}(b_{1},q_{1}),\ldots, \bar{\phi}(b_{n},q_{n}) \right) &= F^{C}_{f} \left( \left\langle \lambda(b_{1},q_{1}),\phi(q_{1})  \right\rangle,\ldots,\left\langle \lambda(b_{n},q_{n}),\phi(q_{n}) \right\rangle  \right) \\
&= \left\langle f^{E} \left( \lambda(b_{1},q_{1}),\ldots, \lambda(b_{n},q_{n}) \right) + S_{f}(\phi(\vec{q})), f^{P}(\phi(\vec{q})) \right\rangle.
\end{align*}
Since $\bar{\phi}$ is a homomorphism, this yields the equation
\begin{align}\label{eq:40}
\lambda \left( f^{B}(\vec{b}) + T_{f}(\vec{q}), f^{Q}(\vec{q}) \right) =  f^{E} \left( \lambda(b_{1},q_{1}),\ldots, \lambda(b_{n},q_{n}) \right) + S_{f}(\phi(\vec{q})).
\end{align}
Define $h: Q \rightarrow E$ by $h(q):= \lambda(0,q)$ and $\psi: B \rightarrow E$ by $\psi := \lambda(b,0')$. By taking the constant $\vec{0} \in B^{\ar f}$, Eq.(\ref{eq:40}) yields 
\begin{align}\label{eq:41}
\lambda \left(T_{f}(\vec{q}), f^{Q}(\vec{q}) \right) = f^{E}(h(\vec{q})) + S_{f}(\phi(\vec{q})).
\end{align}
Applying Eq.(\ref{eq:40}) to the difference term we see that
\begin{align*}
\lambda(b,q) &= \lambda \left( m^{B}(b,0,0) + T_{m}(0',0',q), m^{Q}(0',0',q) \right) \\
&= m^{E} \left( \lambda(b,0'),\lambda(0,0'),\lambda(0,q) \right) + S_{m}(\phi(0'),\phi(0'),\phi(q)) \\
&= m^{E}(\psi(b),\psi(0),h(q)) \\
&= \psi(b) - \psi(0) + h(q)
\end{align*}
because the identity $m(x,x,y)=y$ implies $T_{m}(x,x,y) = 0$ in any compatible 2-cocycle $T$. Then in Eq.(\ref{eq:40}) we have 
\[
f^{E}(h(\vec{q})) + S_{f}(\phi(\vec{q})) = \lambda \left(T_{f}(\vec{q}), f^{Q}(\vec{q}) \right) = \psi(T_{f}(\vec{q}) - \psi(0) + h(f^{Q}(\vec{q}))
\]
which shows $[\psi \circ T] = [S \circ \phi]$ provided we can show $\psi$ is a homomorphism - this is the only point were we use the assumption that $Q$ has an idempotent element $0'$. In Eq.(\ref{eq:40}) taking the constant $0'$-tuples we have 
\[
\psi(f^{B}(\vec{b})) = \lambda \left( f^{B}(\vec{b}), 0' \right) =  f^{E} \left( \lambda(b_{1},0'),\ldots, \lambda(b_{n},0') \right) + S_{f}(\phi(\vec{0'})) = f^{E}\left( \psi(\vec{b})) \right);
\]
thus, $\psi$ is a homomorphism and so $\psi(0)=0$. Then $\delta_{T}(\psi) = [S \circ \phi]=[R]$.

(3) This is precisely the same argument as in (2), but given $[R] \in H^{2}_{\mathcal U}\left(Q,E^{\tau} \right)$ we use the lifting property on the identity $\id : Q \rightarrow Q$ to yield a homomorphism $\bar{\phi} : A \rightarrow E^{\tau} \otimes^{R} Q$ such that $p_{2} \circ \bar{\phi} = \rho$. The rest follows.
\end{proof}

\begin{lemma}\label{lem:precover}
Let $\mathcal U$ a variety in the signature $\tau$ with a difference term and $F_{\mathcal U}(X)/\theta \approx Q$ a free presentation such that $F_{\mathcal U}(X)/[\theta,1]$ has an idempotent element. Then for central datum $(Q,E^{\tau})$ in $\mathcal U$, the transgression map
\[
\delta : \Hom( F'(\theta')/\Delta_{\theta' 1}, E^{\tau}) \longrightarrow H^{2}_{\mathcal U}(Q,E^{\tau})
\]
is surjective where $F' = F_{\mathcal U}(X)/[\theta,1]$ and $\theta' = \theta/[\theta,1]$.
\end{lemma}
\begin{proof}
Surjectivity of the transgression map follows by Lemma~\ref{lem:100}(1) and Proposition~\ref{prop:surjtrans}(3c).  
\end{proof}

In Theorem~\ref{thm:schurmult}, Eq~\eqref{eqn:schurhopf} is our version of the Schur-Hopf formula characterizing the Schur Multiplier of regular datum in varieties with a difference term. Note that we require an idempotent element in the central extension derived from a free-presentation.

\begin{theorem}\label{thm:schurmult}
Let $\mathcal U$ a variety with a difference term and $(Q,E^{\tau})$ regular datum in $\mathcal U$. Let us assume $F_{\mathcal U}/\theta \approx Q$ is a free-presentation such that $F_{\mathcal U}/[\theta,1]$ has an idempotent element. Then
\begin{align}\label{eqn:schurhopf}
H^{2}_{\mathcal U}(Q,E^{\tau}) \approx \Hom \left( \kappa_{ \theta/[\theta,1] 1} \Big( F_{\mathcal U}/[\theta,1] \big( \theta \wedge [1,1]/[\theta,1] \big) \Big) , E^{\tau} \right) \approx \Hom \left( I_{\frac{\alpha \wedge [1,1]}{[\alpha,1]}} , E^{\tau} \right).
\end{align}
\end{theorem}
\begin{proof}
We have seen that $F_{\mathcal U}/[\theta,1] \rightarrow Q$ is a central extension with kernel $\theta/[\theta,1]$. We also see that $[1/[\theta,1],1/[\theta,1]] = [1,1]/[\theta,1]$ and so $\theta/[\theta,1] \wedge [1/[\theta,1],1/[\theta,1]] = \theta \wedge [1,1]/[\theta,1]$. The result now follows from Proposition~\ref{prop:schurmult} and Lemma~\ref{lem:precover}.
\end{proof}

Here we should remark how the classical result for groups can be derived from Theorem~\ref{thm:schurmult}. We also consider how the analogous result can be derived for certain modules expanded by multilinear operations.

\begin{example}
We show how the classic characterization of the Schur Multiplier for groups can be derived from Theorem~\ref{thm:schurmult}. For a group $G$, we have the bijective correspondence between normal subgroups $K\triangleleft G$ and congruences $\alpha_{K} = \{(a,b) \in M^{2}: ab^{-1} \in K \}$. By \cite[Lem 2.10]{wiresI}, there is an isomorphism $G(\alpha_{K})/\Delta_{\alpha_{K} \alpha_{G}} \approx K/[K,G]$.

Let $F/R \approx Q$ be a free-presentation in groups. Then  
\begin{align*}
F/[R, F](\alpha_{R/[R,F]})/\Delta_{\alpha_{R/[R,F]} 1} \approx \frac{R/[R,F]}{[R/[R,F],F/[R,F]]} = R/[R,F]
\end{align*}
and 
\begin{align}
\kappa_{ \alpha_{R}/[\alpha_{R},1] 1} \Big( F/[R,F] \big( \alpha_{R} \wedge [1,1]/[\alpha_{R},1] \big) \Big) \approx \frac{R \wedge [F,F]}{[R,F]}.
\end{align}
If $Q$ is finite, then $K \wedge [G,G]$ is finite for any central extension $K \hookrightarrow G \rightarrow Q$. Applied to the central extension $R/[R,F] \hookrightarrow F/[R,F] \rightarrow Q$ we conclude $R \wedge [F,F]/[R,F]$ is finite. Since finite abelian groups are self-dual, we recover from a free-presentation $F/R \approx Q$ of a finite group the characterization of the Schur Multiplier
\begin{align}
\begin{split}
M(Q)= H^{2}(Q, \mathds{C}^{\times}) &\approx \Hom \Big( \kappa_{ \alpha_{R}/[\alpha_{R},1] 1} \Big( F/[R,F] \big( \alpha_{R} \wedge [1,1]/[\alpha_{R},1] \big) \Big) , \mathds{C}^{\times} \Big) \\
&\approx \Hom \left(\frac{R \wedge [F,F]}{[R,F]} , \mathds{C}^{\times} \right) \\
&\approx \frac{R \wedge [F,F]}{[R,F]}.
\end{split}
\end{align}
\end{example}

\begin{theorem}\label{thm:multischur}
Let $\mathcal V$ be a variety of vector spaces over the field $\mathds{F}$ expanded by multilinear operations whose arities are at least binary. If $F/R \approx Q$ is a free-presentation in $\mathcal V$, then  
\begin{align}\label{eqn:schur1000}
H^{2}_{\mathcal V}(Q,\mathds{F}) \approx \frac{R \wedge [F,F]}{[R,F]}.
\end{align}
\end{theorem}
\begin{proof}
 We write $\alpha_{R}$ for the congruence determined by the ideal $R \triangleleft F$. Using Example~\ref{ex:multiext} and the isomorphism Eq~\eqref{eqn:1000}, we see that
\begin{align*}
F/[R,F](\alpha_{R/[R,F]})/\Delta_{\alpha_{R/[R,F]} 1} \approx \frac{R/[R,F]}{[R/[R,F],F/[R,F]]} = R/[R,F]
\end{align*}
and
\begin{align*}
\kappa_{ \alpha_{R}/[\alpha_{R},1] 1} \Big( F/[R,F] \big( \alpha_{R} \wedge [1,1]/[\alpha_{R},1] \big) \Big) \approx \frac{R \wedge [F,F]}{[R,F]}.
\end{align*}
We also noted that the multilinear operations must be trivial which implies each kernel algebra $B^{\tau} \in \mathcal V$ is just a vector space over $\mathds{F}$ and so isomorphic to a direct sum $\bigoplus_{v \in S} \mathds{F}_{v}$ with each $\mathds{F}_{v} = \mathds{F}$ for a choice of basis $S$. Then distributivity of homsets over direct sums and $\Hom(\mathds{F},\mathds{F}) \approx \mathds{F}$ implies 
\begin{align*}
\Hom \left( B^{\tau}, \mathds{F} \right) \approx \Hom \left( \bigoplus_{v \in S} \mathds{F}_{v}, \mathds{F} \right) \approx \bigoplus_{v \in S} \Hom \left( \mathds{F}_{v}, \mathds{F} \right) \approx \bigoplus_{v \in S} \mathds{F}_{v} \approx B^{\tau}
\end{align*} 
for any kernel algebra $B^{\tau} \in \mathcal V$. Now altogether, Eq~\eqref{eqn:schurhopf} in Theorem~\ref{thm:schurmult} takes the form 
\begin{align*}
\begin{split}
H^{2}_{\mathcal V}(Q,\mathds{F}) \approx \Hom \Big( \kappa_{ \alpha_{R}/[\alpha_{R},1] 1} \Big( F/[R,F] \big( \alpha_{R} \wedge [1,1]/[\alpha_{R},1] \big) \Big), \mathds{F} \Big) &\approx \Hom \Big(\frac{R \wedge [F,F]}{[R,F]} , \mathds{F} \Big) \\
&\approx \frac{R \wedge [F,F]}{[R,F]}
\end{split}
\end{align*}
\end{proof}

\begin{example}
In any particular variety, we can breakdown Eq~\eqref{eqn:schurhopf} or Eq~\eqref{eqn:schur1000} further by introducing any knowledge of the generators of the commutator; for example, in the case of modules expanded by multilinear operations we have the ideal closures determined by the multilinear operations. Let us observe the form it may take in three particular examples.

In the case of Leibniz algebras $\left\langle L , + , 0, \mathds{F}, \cdot \right\rangle$ over a field $\mathds{F}$ the commutator is generated by $[R,F] = R\cdot F + F \cdot R$ and so Eq~\eqref{eqn:schur1000} becomes
\begin{align}
H^{2}_{\mathcal V}(Q,\mathds{F}) \approx \frac{R \wedge (F \cdot F )}{R\cdot F + F \cdot R}.
\end{align}
Diassociative algebras $\left\langle L , + ,0,\mathds{F}, \vdash, \dashv \right\rangle$ are vector spaces over $\mathds{F}$ with two bilinear maps satisfying additional identities \cite{mainellis2}. If we write $R \diamond F = R \vdash F + R \dashv F$, then the commutator is generated by $[R,F] = R \diamond F + F \diamond R$ and so Eq~\eqref{eqn:schur1000} takes the form
\begin{align}
H^{2}_{\mathcal V}(Q,\mathds{F}) \approx \frac{R \wedge (F \diamond F )}{R \diamond F + F \diamond R}.
\end{align}
An n-Lie algebra $\left\langle L , + ,0,\mathds{F}, \{-,\ldots,-\}_{\mathds{F}} \right\rangle$ is a vector space over $\mathds{F}$ with an n-ary skew-symmetric multilinear map satisfying an additional identity (see Filippov \cite{filipov}). Then Eq~\eqref{eqn:schur1000} becomes 
\begin{align}
H^{2}_{\mathcal V}(Q,\mathds{F}) \approx \frac{R \wedge \{F,\ldots,F \}_{\mathds{F}}}{ \sum_{i=1}^{n} \{F,\ldots,R,\ldots,F \}_{\mathds{F}}}.
\end{align}
\end{example}

We may consider what happens in Theorem~\ref{thm:multischur} when the expanded multilinear operations also include linear operations; in this case, centrality of the ideal does not  trivialize the linear operations and so they are naturally folded into the ring of operators.

\begin{example}\label{ex:manyunary}
Let $S$ be a ring with unit and $\mathcal V$ a variety of unital $S$-modules expanded by multilinear operations indexed by $F$. Let $F_{1}$ index the linear operations in $F$ and assume $F_{1} \neq \emptyset$. Let us observe that since each module has the idempotent $0$, abelian and kernel algebras are the same. As previously seen, the multilinear operations of at least binary arity in a kernel algebra are all trivial; therefore, each abelian algebra in $\mathcal V$ is an $S$-module with additional linear operations index by $F_{1}$. If we take $S \left\langle  F_{1} \right\rangle$ to be the non-commutative free $S$-algebra in the indeterminates $F_{1}$, then each abelian algebra $B^{\tau} \in \mathcal V$ becomes an $S \left\langle  F_{1} \right\rangle$-module by evaluation $L \ast a = L^{B}(b)$ for $L \in F_{1}, b \in B^{\tau}$; therefore, the abelian algebras in $\mathcal V$ form a variety of $S \left\langle  F_{1} \right\rangle$-modules. We can be more precise. The abelian algebras of $\mathcal V$ form a subvariety $\mathcal A$. Since the nonlinear multilinear operations are all trivial in $\mathcal A$ and we have an idempotent element $0$, any identity $f = g \in \mathrm{Id} \, \mathcal V$ when interpreted in an abelian algebra $B^{\tau} \in \mathcal A$ is equivalent to a $S \left\langle  F_{1} \right\rangle$-module equation
\begin{align*}
(f_{1} - g_{1}) \ast x_{1} + \cdots + (f_{n} - g_{n}) \ast x_{n} = 0.
\end{align*}
By making substitutions $x_{i}=0$, we see how the identities of $\mathcal A$ are determined by the 1-variable module identities which define an ideal $I_{\mathcal V} =  \{ f-g \in S \left\langle  F_{1} \right\rangle : (f - g) \ast x = 0 \text{ for all } B^{\tau} \in \mathcal A \} \triangleleft S \left\langle  F_{1} \right\rangle$. In this manner, we observe that upto term-equivalence, the abelian algebras of $\mathcal V$ are precisely $S \left\langle  F_{1} \right\rangle/I_{\mathcal V}$-modules.

We then take a divisible abelian group $D$ which embeds the 1-generated $S \left\langle  F_{1} \right\rangle/I_{\mathcal V}$-modules as abelian groups. Then by Proposition~\ref{prop:regexists}, $\Hom_{\mathds{Z}}(S \left\langle  F_{1} \right\rangle/I_{\mathcal V}, D )$ is a regular algebra where the higher multilinear operations are interpreted as trivial operations. If $F/R \approx Q$ is a free presentation in $\mathcal V$, then Eq~\ref{eqn:schurhopf} takes the form
\begin{align}\label{eqn:manylinear}
H^{2}_{\mathcal V} \left( Q, \Hom_{\mathds{Z}} \big( S \left\langle  F_{1} \right\rangle/I_{\mathcal V}, D \big) \right) \approx \Hom \left( \frac{R \wedge [F,F]}{[R,F]} , \Hom_{\mathds{Z}} \big( S \left\langle  F_{1} \right\rangle/I_{\mathcal V}, D \big) \right).
\end{align}
\end{example}

In the case in which $F_{1}$ consists of a single operation and the ring is commutative, then a general simplification is possible.

\begin{theorem}\label{thm:singleunary}
Let $k$ be a commutative ring with unity and $\mathcal V$ a variety of unital $k$-modules expanded by multilinear operations in which the linear operations consist of but a single operation $L$. If $F/R \approx Q$ is a free-presentation in $\mathcal V$, then  
\begin{align}\label{eqn:singleunary}
H^{2}_{\mathcal V} \left( Q, \Hom_{\mathds{Z}}(k[L]/I_{\mathcal V}, D ) \right) \approx \Hom_{\mathds{Z}} \left( \frac{R \wedge [F,F]}{[R,F]} , D \right)
\end{align}
where $D$ is any divisible abelian group which embeds the 1-generated $k[L]/I_{\mathcal V}$-modules. 
\end{theorem}
\begin{proof}
From the discussion in Example~\ref{ex:manyunary}, the abelian algebras in $\mathcal V$ are term-equivalent to unital $k[L]/I_{\mathcal V}$-modules over the commutative ring $k[L]/I_{\mathcal V}$ with a unit. Then we can treat $\frac{R \wedge [F,F]}{[R,F]}$ as a right-module over $k[L]/I_{\mathcal V}$ and $k[L]/I_{\mathcal V}$ as a left-module over itself. As abelian groups, both $k[L]/I_{\mathcal V}$ and $D$ are right $\mathds{Z}$-modules. We can then simplify the nested homsets using the tensor identities
\begin{align*}
H^{2}_{\mathcal V} \left( Q, \Hom_{\mathds{Z}} \big(k[L]/I_{\mathcal V}, D \big) \right) &\approx \Hom \left( \frac{R \wedge [F,F]}{[R,F]} , \Hom_{\mathds{Z}} \big( k[L]/I_{\mathcal V}, D \big) \right) \\
&\approx \Hom_{\mathds{Z}} \left( \frac{R \wedge [F,F]}{[R,F]} \otimes k[L]/I_{\mathcal V} , D \right) \\
&\approx \Hom_{\mathds{Z}} \left( \frac{R \wedge [F,F]}{[R,F]} , D \right).
\end{align*}
\end{proof}

\begin{example}
A classical \emph{Rota-Baxter algebra of weight} $\lambda$ is an associative algebra $A = \left\langle A, +, -, 0 , k, \cdot, L \right\rangle$ over the commutative ring $k$ with a linear operator $L$ and constant $\lambda \in k$ such that 
\begin{align*}
L(x) \cdot L(y) &= L \left( L(x)\cdot y + x \cdot L(y) + \lambda L(x \cdot y) \right) &(x,y \in A).
\end{align*}
Let $\mathcal R \mathcal B_{\lambda}$ be the variety of Rota-Baxter algebras over $k$ of weight $\lambda$. Since in abelian algebras the bilinear product is trivial, the above equation puts no restrictions on the linear operator $L$ in abelian algebras; therefore, $I_{\mathcal R \mathcal B_{\lambda}} = 0$. If $F/R \approx Q$ is a presentation of $Q$ in the variety $\mathcal R \mathcal B_{\lambda}$, then it follows from Eq~\ref{eqn:singleunary} in Theorem~\ref{thm:singleunary} that
\begin{align*}
H^{2}_{\mathcal R \mathcal B_{\lambda}} \left( Q, \Hom_{\mathds{Z}}(k[L], D ) \right) \approx \Hom_{\mathds{Z}} \left( \frac{R \wedge [F,F]}{[R,F]} , D \right)
\end{align*}
where $D$ is a divisible group which embeds the 1-generated $k[L]$-modules.
\end{example}

\begin{lemma}\label{lem:34}
Let $\mathcal U$ be a variety with a difference term with $A \in \mathcal U$ and $\rho : A \rightarrow Q$ a central extension realizing regular datum $(Q,B^{\tau})$. If the induced restriction map 
\[
\check{r} : \Hom(A,E^{\tau}) \rightarrow \Hom(B,E^{\tau})
\] 
is surjective, then $\ker \rho \wedge [1,1] = 0$. If $A$ has an idempotent element, then the converse holds.
\end{lemma}
\begin{proof}
In the second paragraph of the proof of Proposition~\ref{prop:schurmult}, it is shown that $\im \check{r} \leq \ker \nabla$ for the restriction homomorphism
\[
\Hom (B^{\tau},E^{\tau}) \stackrel{\nabla}{\longrightarrow} \Hom \left( \kappa_{\ker \rho 1} \left( A(\ker \rho \wedge [1,1]) \right), E^{\tau} \right).
\] 
Because $E^{\tau}$ separates points, $\ker \rho \wedge [1,1] \neq 0 $ if and only if $\kappa_{\ker \rho 1} \left( A(\ker \rho \wedge [1,1]) \right) \neq 0$ if and only if $\ker \nabla \neq \Hom(B^{\tau},E^{\tau})$; consequently, if $\check{r}$ is surjective, then $\ker \rho \wedge [1,1] = 0$. In the case $A$ has an idempotent element, Proposition~\ref{prop:schurmult} establishes the equality $\im \check{r} = \ker \nabla$ which yields the reverse implication. 
\end{proof}

\begin{proposition}\label{prop:traninject}
Let $\mathcal U$ a variety with a difference term and $(Q,B^{\tau})$ central datum in $\mathcal U$. Let $\rho: A \rightarrow Q$ an extension represented by $[T] \in H^{2}_{\mathcal U}(Q,B^{\tau})$ and $\delta_{T} : \Hom (B^{\tau}, E^{\tau}) \longrightarrow H^{2}_{\mathcal U}(Q,E^{\tau})$ the associated transgression for the regular kernel algebra $E^{\tau} \in \mathcal U$.
\begin{enumerate}

	\item If $\ker \pi \leq [1,1]$, then the transgression $\delta_{T}$ is injective. 
	
	\item If $A$ has an idempotent element and $\delta_{T}$ is injective, then $\ker \pi \leq [1,1]$. 

\end{enumerate}	
\end{proposition}
\begin{proof}
We can identify $A \approx B^{\tau} \otimes^{T} Q$ and write $\alpha=\ker \pi$.  By Theorem~\ref{thm:inflationrestriction}, $\delta_{T}$ is injective if and only if $\im \check{r} = 0$; therefore, it is enough to show $\im \check{r} = 0$ if and only if $\alpha \leq [1,1]$.

(1) Since $E^{\tau}$ is affine, $[1,1] \leq \ker \phi$ for all $\phi \in \Hom \left( A,E^{\tau} \right)$; thus, $\alpha \leq [1,1] \leq \ker \phi$. Since $\left\langle (b,t),(0,t) \right\rangle \in \alpha$ for any $b \in B, t \in Q$, we see that $\check{r}(\phi)(b) = \phi(b,t) - \phi(0,t) = 0$; therefore, $\im \check{r} = 0$.

(2) Now assume $A$ has an idempotent element and $\im \check{r} = 0$. In the proof of Proposition~\ref{prop:schurmult}, it is shown that $\im \check{r} = \ker \nabla$ for the restriction $\nabla: \Hom(B^{\tau},E^{\tau}) \rightarrow \Hom ( \kappa_{\alpha 1}(\alpha \wedge [1,1]),E^{\tau})$. Suppose $\alpha \nleq [1,1]$. If $(a,b) \in \alpha - [1,1]$, then $\begin{bmatrix} b \\ a \end{bmatrix}/\Delta_{\alpha 1} \not\in \kappa_{\alpha 1}\left( \alpha \wedge [1,1] \right)$; thus, $\kappa_{\alpha 1}\left( \alpha \wedge [1,1] \right) < B^{\tau}$. Since the terms of $B^{\tau}$ are module terms, the submodule $\kappa_{\alpha 1}\left( \alpha \wedge [1,1] \right)$ determines a congruence $\theta \in \Con B^{\tau}$ with $\theta < 1$. By regularity of $E^{\tau}$, there is a homomorphism $\phi: B^{\tau} \rightarrow E^{\tau}$ such that $\left( \begin{bmatrix} b \\ a \end{bmatrix}/\Delta_{\alpha 1} , 0 \right) \not\in \ker \phi$ but $\theta \leq \ker \phi$. This implies $\phi \neq 0$ but $\phi \in \ker \nabla$, and so produces the contradiction $\im \check{r} \neq 0$. It must be that $\alpha \leq [1,1]$. 
\end{proof}

\begin{remark}
$\delta_{T}$ injective. Then $0=\ker \delta_{T} = \im \check{r}$. So for all $\phi: A \rightarrow E^{\tau}$, $0 = \check{r}(\phi) = \phi(b,q) - \phi(0,q)$. Then $\left\langle (b,q),(0,q) \right\rangle \in \ker \phi$ which implies $\ker \pi \leq \ker \phi$ for all $\phi: A \rightarrow E^{\tau}$.
\end{remark}

\begin{definition}
Let $\mathcal U$ be a variety with $A,Q \in \mathcal U$. The algebra $A$ is a $\mathcal U$-\emph{cover} of the algebra $Q$ if there is an extension $\pi: A \rightarrow Q$ which has the $\mathcal U$-lifting property for $Q$ and $\ker \pi \leq [1_{A},1_{A}] \wedge \zeta_{A}$. 
\end{definition}

The following theorem concerns the existence of covers. In Theorem~\ref{thm:schurmult}, we observed that the Schur Multiplier of regular datum is isomorphic to the homset of a subalgebra of a kernel algebra and the regular algebra. In Theorem~\ref{thm:coversexist}, we construct a central extension so that the Schur Multiplier is isomorphic to homset of the induced kernel algebra and the regular algebra.

\begin{theorem}\label{thm:coversexist}
Let $\mathcal U$ a variety with a difference term which contains regular datum $(Q, E^{\tau})$. Let $F_{\mathcal U}(X)/\theta \approx Q$ be a free-presentation such that $F_{U}(X)/[\theta,1]$ has an idempotent element. There is a cover $\pi : A \rightarrow Q$ which realizes central datum $\left( Q,\kappa_{\theta/[\theta,1] 1} \left(F/[\theta,1] \big( \theta \wedge [1,1]/[\theta,1] \big) \right) \right) \approx (Q, I_{\frac{\alpha \wedge [1,1]}{[\alpha , 1]}})$.
\end{theorem}
\begin{proof}
We are given the free algebra $F=F_{\mathcal U}(X)$ and $\theta \in \Con F$ such that $F/\theta = Q$. For the central congruence $\theta/[\theta,1] \in \Con F/[\theta,1]$, there is a representation $F/[\theta,1] \approx F/[\theta,1](\theta/[\theta,1])/\Delta_{\theta/[\theta,1] 1} \otimes^{T} Q$. For simplicity, we write $F' = F/[\theta,1]$, $\theta'=\theta/[\theta,1]$ and $F''=F'/[1',1']$ with the total congruences $1' = 1_{F'}$, $1'' = 1_{F''}$. By Lemma~\ref{lem:31}, there is the direct sum
\begin{align*}
F'(\theta')/\Delta_{\theta' 1'} \approx \kappa_{\theta' 1'} \left(F' \big( \theta' \wedge [1',1'] \big) \right) \bigoplus F'/[1',1']( \beta' )/\Delta_{\beta' 1''}
\end{align*}
where $\beta' = (\theta' \vee [1',1'])/[1',1']$. If $p_{1}$ is the first-projection in the above direct sum, then we see that $S:=p_{1} \circ T$ is a $\mathcal U$-compatible 2-cocycle. The claim is that $A := \kappa_{\theta' 1'} \left(F' \big( \theta' \wedge [1',1'] \big) \right) \otimes^{S} Q$ is a cover of $Q$.  According to Proposition~\ref{prop:surjtrans} and Proposition~\ref{prop:traninject}, the claim follows if we show the transgression $\delta_{S}$ is bijective.

From Lemma~\ref{lem:precover}, we see that $\delta_{T}$ is surjective; therefore, $\im \delta_{T} = H^{2}_{\mathcal U}(Q,E^{\tau})$. In the proof of Proposition~\ref{prop:schurmult}, we showed that $\ker \delta_{T} = \ker \nabla$ for the restriction homomorphism
\begin{align*}
\Hom (F'(\theta')/\Delta_{\theta' 1'} ,E^{\tau}) \stackrel{\nabla}{\longrightarrow} \Hom \left( \kappa_{\theta' 1'} \left(F' \big( \theta' \wedge [1',1'] \big) \right), E^{\tau} \right).
\end{align*}
It follows from the isomorphism
\begin{align}\label{eqn:directsumhom}
\begin{split}
\Hom \big( \kappa_{\theta' 1'} \left(F' \big( \theta' \wedge [1',1'] \big) \big) , E^{\tau} \right) \bigoplus \Hom \big( F'/[1',1']( \beta' )/\Delta_{\beta' 1''} , E^{\tau} \big) &\approx \Hom( F'(\theta')/\Delta_{\theta' 1'}, E^{\tau} ) \\
(\phi_{1},\phi_{2}) &\longmapsto \phi_{1} + \phi_{2}
\end{split}
\end{align}
that $\ker \delta_{T} = \ker \nabla = \Hom ( F'/[1',1']( \beta' )/\Delta_{\beta' 1''} , E^{\tau})$. From the definition of $S$, it follows from Eq~\eqref{eqn:directsumhom} that $\ker \delta_{S} = 0$ and $H^{2}_{\mathcal U}(Q,E^{\tau})= \im \delta_{T} = \im \delta_{S}$.
\end{proof}

Let us recall that the transgression is a homomorphism in both coordinates. In Theorem~\ref{thm:inflationrestriction}, the choice of 2-cocycle is fixed and the transgression is considered as a function of the corresponding homset. The next theorem treats the transgression as a function of $2^{\mathrm{nd}}$-cohomology.

\begin{theorem}\label{thm:47}
Let $\mathcal U$ a variety with a difference term containing central datum $(Q,B^{\tau})$. Let $E^{\tau} \in \mathcal U$ be a regular kernel algebra. Let $F_{\mathcal U}(X)/\theta \approx Q$ be a free-presentation such that $F_{\mathcal U}(X)/[\theta,1]$ has an idempotent element. Then the sequence 
\begin{align}\label{eq:product}
0 \longrightarrow \mathrm{Ext}_{\mathcal U} \big( Q/[1,1], B^{\tau} \big) \stackrel{\check{\sigma}}{\longrightarrow} H^{2}_{\mathcal U}(Q, B^{\tau}) \stackrel{\delta}{\longrightarrow} \Hom \Big( \Hom (B^{\tau} ,E^{\tau}),H^{2}_{\mathcal U}(Q, E^{\tau}) \Big) 
\end{align}
is exact. 
\end{theorem}
\begin{proof}
Note $\check{\sigma} :  \mathrm{Ext}_{\mathcal U} \big( Q/[1,1], B^{\tau} \big) \rightarrow  H^{2}_{\mathcal U}(Q,B^{\tau})$ is the restriction of the inflation map to the subgroup $\mathrm{Ext}_{\mathcal U} \big( Q/[1,1] , B^{\tau} \big) \leq H^{2}_{\mathcal U} \big( Q/[1,1],B^{\tau} \big)$. We first show $\check{\sigma}$ is injective.

Suppose $[S] \in \ker \check{\sigma}$. Let $\pi: Q \rightarrow Q/[1,1]$ be the canonical homomorphism. Let $S$ be the 2-cocycle corresponding to the central extension $\rho_{2}: A_{2} \rightarrow Q/[1,1]$ and $\rho_{1} : A_{1} \rightarrow Q$ the central extension corresponding to $[S \circ \pi] = \check{\sigma}([S]) = 0$. By realization, there is a lifting $l: Q \rightarrow A_{1}$ which defines $S \circ \pi$: for all $f \in \tau$, $\vec{q} \in Q^{\ar f}$,
\[
S_{f}(\pi(\vec{q})) = \begin{bmatrix} l(f^{Q}(\vec{q})) \\ f^{A_{1}}(l(\vec{q})) \end{bmatrix}/\Delta_{\ker \rho_{1} 1}.
\]
Since $[S \circ \pi] = 0$, we see that $S_{f}(\pi(\vec{q})) = \delta$ which implies $l$ is a homomorphism and $A_{1} \approx B^{\tau} \times Q$. By Lemma~\ref{lem:30}, there is a homomorphism $\phi: A_{1} \rightarrow A_{2}$ such that $\rho_{2} \circ \phi = \pi \circ \rho_{1}$. For the homomorphism $\phi \circ l : Q \rightarrow A_{2}$, since $A_{2}$ is abelian, $[1,1] \leq \ker \phi \circ l$ and so there is an induced homomorphism $\eta: Q/[1,1] \rightarrow A_{2}$ such that $\pi \circ \pi = \phi \circ l$. Then because $\pi$ is an epimorphism and $\rho_{2} \circ \eta \circ \pi = \rho_{2} \circ \phi \circ l = \pi \circ \rho_{1} \circ l = \pi$ we conclude $\rho_{2} \circ \eta$. This implies $\eta$ is a lifting for $\rho_{2}$ and so defines a compatible 2-cocycle $T'$ for $A_{2}$ by 
\[
T'_{f}(\vec{q}) = \begin{bmatrix} \eta(f^{Q/[1,1]}(\vec{q})) \\ f^{A_{2}}(\eta(\vec{q})) \end{bmatrix}/\Delta_{\ker \rho_{2} 1} \quad \quad \quad \quad \quad (f \in \tau, \vec{q} \in Q/[1,1]^{\ar f}).
\] 
Then $T' \sim S$ and because $\eta$ is a homomorphism, we conclude that $[S]=[T']=0$; thus, $\check{\sigma}$ is injective.

We now show exactness at $H^{2}_{\mathcal U}(Q, B^{\tau})$. This is the only point where we utilize the hypothesis of an idempotent element in the extensions in $\mathcal U$ realizing datum $(Q,B^{\tau})$. According to Theorem~\ref{thm:inflationrestriction}, $[T] \in \ker \delta$ if and only if $\Hom(B^{\tau},E^{\tau}) = \ker \delta(-,[T])$ if and only if the restriction $\check{r} : \Hom(A,E^{\tau}) \rightarrow \Hom(B^{\tau},E^{\tau})$ is surjective where $\rho: A \rightarrow Q$ is a central extension represented by $[T]$. By Lemma~\ref{lem:34}, this occurs precisely when $\ker \rho \mathrel{\wedge} [1,1] = 0$ in $A$ which is equivalent to $[T] \in \im \check{\sigma}$ according to Lemma~\ref{lem:32}. 
\end{proof}

\vspace{0.2cm}

\section{Universal central extensions and perfect algebras}\label{sec:6}

The uniqueness of the solution in the mapping diagram for the lifting property is intimately connected to the topic of universal central extensions and extensions of perfect algebras.

\begin{definition}
Let $\mathcal U$ be variety and $\pi: A \rightarrow Q$ be a central extension in $\mathcal U$. The extension $\pi: A \rightarrow Q$ is a $\mathcal U$-\emph{universal central extension} if for all central extensions $\rho: E \rightarrow P$ in $E \in \mathcal U$ and homomorphisms $\tau : Q \rightarrow P$, there is a unique homomorphisms $\bar{\tau}: A \rightarrow E$ such that $\rho \circ \bar{\tau} = \tau \circ \pi$.
\end{definition}

A universal central extension $\pi: A \rightarrow Q$ in $\mathcal U$ has the $\mathcal U$-central lifting property with unique solutions.

\begin{lemma}\label{lem:perfectext}
Let $A \in \mathcal U$ a variety with a difference term and $\pi: A \rightarrow Q$ a central extension with $Q$ perfect. Then $\ker \pi \vee [1_{A},1_{A}]=1_{A}$ and $[1_{A},1_{A}]$ is neutral.
\end{lemma}
\begin{proof}
Let $\alpha = \ker \pi$. By the homomorphism property
\begin{align*}
1 = [1,1] = [1/\alpha,1/\alpha] = \left( [1,1] \vee \alpha \right)/\alpha
\end{align*}
which implies $1 = [1,1] \vee \alpha$ in $A$ since $\alpha \leq [1,1] \vee \alpha$. Then 
\begin{align*}
[1,1] = [[1,1] \vee \alpha,[1,1] \vee \alpha] = [[1,1],[1,1]] \vee [[1,1],\alpha] \vee [\alpha,[1,1]] \vee [\alpha,\alpha] = [[1,1],[1,1]]
\end{align*}
by centrality of $\alpha$.
\end{proof}

The following is our analogue for varieties with a difference term of the characterization of perfect groups and universal central extensions by $2^{\mathrm{nd}}$-cohomology (see Milnor \cite{milnor}).

\begin{theorem}\label{thm:perfectchar}
Let $\mathcal U$ be a variety with a difference term and $Q \in \mathcal U$. Let $F/\theta \approx Q$ be a free presentation in $\mathcal U$ such that $F/[\theta,1]$ has an idempotent element.
\begin{enumerate}

	\item $Q$ has a universal central extension in $\mathcal U$ if and only if $Q$ is perfect.
	
	\item Assume $Q \in \mathcal U$ is perfect and let $A \in \mathcal U$. A central extension $\pi: A \rightarrow Q$ is a universal central extension in $\mathcal U$ if and only if $A$ is perfect and $H^{2}_{\mathcal U}(A,B^{\tau})=0$ for all datum $(A,B^{\tau})$ in $\mathcal U$.

\end{enumerate}
\end{theorem}
\begin{proof}
(1) First, we show necessity. Suppose $Q$ is not perfect. Then $Q/[1,1]$ is abelian and nontrivial. We claim $p_{2} : Q/[1,1] \times Q \rightarrow Q$ is a central extension. To see this, note $\ker p_{1} = 0 \times 1_{Q}$ and $\ker p_{2} = 1/[1,1] \times 0_{Q}$. Then 
\begin{align*}
[\ker p_{2},1] = [\ker p_{2} , \ker p_{2} \vee \ker p_{1}] &= [\ker p_{2} , \ker p_{2}] \vee [\ker p_{2} , \ker p_{1}] \\
&\leq \big( [1/[1,1],1/[1,1]] \times [0_{Q},0_{Q}] \big) \vee \left( \ker p_{2} \wedge \ker p_{1} \right) \\
&= \left( [1,1] \vee [1,1] \right)/[1,1] \times 0_{Q}   \\
&= 0 . 
\end{align*}
The set $A= Q/[1,1] \times Q$ and we have central extension $A \approx A(\alpha)/\Delta_{\alpha 1} \times Q \stackrel{p_{2}}{\rightarrow} Q$ if we choose the lifting $l(x):= \left\langle 0, x \right\rangle$ for $p_{2}$ and we calculate $T_{f} = \hat{\delta}$ where this lifting is used to define it. Take the canonical epimorphism $\pi: Q \rightarrow Q/[1,1]$ and $\kappa : A \rightarrow A(\alpha)/\Delta_{\alpha 1}$. Now let $\rho: E \rightarrow Q$ be any central extension. There is a homomorphism $\phi: E \rightarrow A(\alpha)/\Delta_{\alpha 1}$ given by $\phi(x):= \hat{\delta}$ for all $x \in E$. Define $\psi, \psi': E \rightarrow  A(\alpha)/\Delta_{\alpha 1} \times Q$ by $\psi := \left( \kappa \circ (\pi \circ \rho , \rho), \rho \right)$ and $\psi' := ( \phi , \rho)$. Then $p_{2} \circ \psi = \rho$ and $p_{2} \circ \psi' = \rho$; thus, $\rho: E \rightarrow Q$ cannot be a universal central extension since the uniqueness condition is violated.

Now, let us show sufficiency. Assume $Q$ is perfect. By Theorem~\ref{thm:coversexist}, there is a cover $\pi: E \rightarrow Q$ of $Q$; that is, $\ker \pi \leq \zeta_{E} \wedge [1_{E},1_{E}]$ and has the $\rho$-lifting property for any central extension $\rho : A \rightarrow Q$ in $\mathcal U$. Since it has the lifting property, any central extension $\sigma: G \rightarrow E$ is a semidirect product; thus, $H^{2}_{\mathcal U}(E,B^{\tau})=0$ for all datum $(E,B^{\tau})$ in $\mathcal U$. Since $Q$ is perfect, we see that $[1_{E},1_{E}]$ in $E$ is neutral by Lemma~\ref{lem:perfectext}. We also have $\ker \pi \vee [1,1] = 1_{E}$ by Eq~\eqref{eqn:calc} whch implies $[1_{E},1_{E}]=1_{E}$ since $\ker \pi \leq [1_{E},1_{E}]$; thus, $E$ is perfect. Then by part (2) it must be that $\pi: E \rightarrow Q$ is a universal central extension.

(2) We are given that $Q$ is perfect. Assume $\pi: A \rightarrow Q$ is a universal central extension in $\mathcal U$. We can write $A = C^{\tau} \otimes^{T} Q$. For any central extensions $\rho: B^{\tau} \otimes^{T} A \rightarrow A$, then the universal central extension property gives a unique homomorphism $\gamma : A \rightarrow B^{\tau} \otimes^{T'} Q$ such that $ \pi \circ \rho \circ \gamma = \pi$. Then the uniqueness condition on the maps yields $\rho \circ \gamma = \id_{A}$. This implies every extension of the datum $(A,B^{\tau})$ is a direct product; thus, $H^{2}_{\mathcal U}(A,B^{\tau}) = 0$.

Now we show $A$ must be perfect; for a contradiction, assume $[1,1] < 1$ in $A$. Then by Lemma~\ref{lem:perfectext}, $[1,1]$ is perfect. Write $\alpha = \ker \pi$. Then
\begin{align}\label{eqn:calc}
\left( [1,1] \vee \alpha \right)/ \alpha = [ \left( 1 \vee \alpha \right)/\alpha, \left( 1 \vee \alpha \right)/\alpha ] = [1/\alpha,1/\alpha] = [1_{Q},1_{Q}] = 1_{Q}
\end{align}
which implies $\alpha \vee [1,1] = 1$ in $A$ by the Correspondence Theorem. It must be that $\alpha \nleq [1,1]$. Then by Proposition~\ref{prop:traninject}(2), the transgression $\delta_{T}: \Hom(C^{\tau},B^{\tau}) \rightarrow H^{2}_{\mathcal U}(Q,B^{\tau})$ is not injective (here we have to assume $A$ has an idempotent element). So there is $\sigma \in \Hom(C^{\tau},B^{\tau})$ and $h: Q \rightarrow B$ such that 
\begin{align*}
\sigma \circ T_{f}(\vec{x}) &= f^{B}(h(\vec{x})) - h(f^{Q}(\vec{x}))  &(f \in \tau).
\end{align*} 
Define $\gamma : A \rightarrow B^{\tau} \times Q$ by $\gamma(c,x) := \left\langle \sigma(c) + h(x) , x \right\rangle$. We check $\gamma$ is a homomorphism: 
\begin{align*}
\gamma \circ F_{f}\left( \left\langle b_{1}, x \right\rangle,\ldots, \left\langle b_{n}, x \right\rangle \right) &= \gamma \left( f^{C}(\vec{b}) + T_{f}(\vec{x}) , f^{Q}(x) \right) \\
&= \left\langle \sigma( f^{C}(\vec{b}) ) + \sigma \circ T_{f}(\vec{x}) +h( f^{Q}(x) ) , f^{Q}(x) \right\rangle \\
&= \left\langle f^{B}( \sigma(\vec{b}) ) + f^{B}(h(\vec{x})) , f^{Q}(x) \right\rangle \\
&= \left\langle f^{B} \circ (\sigma + h)(\vec{b}) , f^{Q}(x) \right\rangle \\
&= F_{f} \left(\gamma (b_{1},x),\ldots, \gamma (b_{1},x) \right)
\end{align*}
Then $\pi = p_{2} \circ \gamma$. We also have $\pi = p_{2} \circ \kappa$ where $\kappa(c,x) := \left\langle 0, x \right\rangle$; however, $\sigma + h \not\equiv 0$, so that $\gamma \neq \kappa$. This violates the uniqueness condition in the universal central extension. It must be that $A$ is perfect.

Conversely, let assume $A$ is perfect and $H^{2}_{\mathcal U}(A,B^{\tau}) = 0$ for all central datum $(A,B^{\tau})$ in $\mathcal U$. The triviality of the second-cohomology implies the transgression $\delta_{T}: \Hom(C^{\tau},B^{\tau}) \rightarrow H^{2}_{\mathcal U}(Q,B^{\tau})$ is surjective by Theorem~\ref{thm:inflationrestriction}. This implies by Proposition~\ref{prop:surjtrans}(1) that $\pi: A \rightarrow Q$ has the $\rho$-lifting property for all central extensions $\rho: E \rightarrow Q$ in $\mathcal U$. Since $A$ is perfect, we have $\alpha \leq 1 = [1,1]$; therefore, Proposition~\ref{prop:traninject}(1) implies $\delta_{T}$ is injective.

Now suppose we have $\psi_{1},\psi_{2} : A \rightarrow B^{\tau} \otimes^{T'} Q$ such that $\pi = p_{2} \circ \psi_{1} =  p_{2} \circ \psi_{2}$ and $\psi_{1} \neq \psi_{2}$. Then we can write $\psi_{i}\left( c,x \right) = \left\langle \sigma_{i}(c,x), x \right\rangle$; thus, $\sigma_{1} \neq \sigma_{2}$.  Since $\psi_{i}$ are homomorphisms, we calculate
\begin{align*}
\psi_{i} \circ F_{f}\left( \left\langle b_{1}, x_{1} \right\rangle, \ldots, \left\langle b_{n}, x_{n} \right\rangle \right) &= \psi_{i} \left( \left\langle f^{C}(\vec{b}) + T_{f}(\vec{x}) , f^{Q}(\vec{x}) \right\rangle \right)\\
&= \left\langle  \sigma \left( f^{C}(\vec{b}) + T_{f}(\vec{x}) , f^{Q}(\vec{x}) \right) , f^{Q}(\vec{x}) \right\rangle
\end{align*}
and
\begin{align*}
F_{f} \left( \psi_{i}(b_{1},x_{1}), \ldots, \psi_{i}(b_{n},x_{n}) \right) &= F_{f} \left( \left\langle \sigma_{i}(b_{1},x_{1}) , x_{1} \right\rangle, \ldots, \left\langle \sigma_{i}(b_{n},x_{n}) , x_{n} \right\rangle  \right) \\
&= \left\langle f^{B}( \sigma_{i}(b_{1},x_{1}),\ldots,\sigma_{i}(b_{1},x_{1}) ) + T_{f}'(\vec{x}), f^{Q}(\vec{x}) \right\rangle ;
\end{align*}
therefore, we conclude that
\begin{align}\label{eqn:homunivcent}
\sigma_{i}( f^{C}(\vec{c}) + T_{f}(\vec{x}) , f^{Q}(\vec{x}) ) = f^{B}( \sigma_{i}(b_{1},x_{1}),\ldots,\sigma_{i}(b_{n},x_{n}) ) + T_{f}'(\vec{x}).
\end{align}
Now define $\gamma : A=C^{\tau} \otimes^{T} Q \rightarrow B^{\tau} \times Q$ by $\gamma(c,x) := \left\langle \sigma_{1}(c,x) - \sigma_{2}(c,x) , x\right\rangle$. Then we see that for $f \in \tau$ with $n=\ar f$, Eq~\eqref{eqn:homunivcent} yields
\begin{align*}
\gamma &\circ F_{f} \left( \left\langle b_{1}, x_{1} \right\rangle,\ldots,\left\langle b_{n}, x_{n} \right\rangle \right) \\
&= \gamma \left( f^{C}(\vec{b}) + T_{f}(\vec{x}), f^{Q}(\vec{x}) \right) \\
&= \left\langle \sigma_{1}( f^{C}(\vec{c}) + T_{f}(\vec{x}) , f^{Q}(\vec{x}) ) - \sigma_{2}( f^{C}(\vec{c}) + T_{f}(\vec{x}) , f^{Q}(\vec{x}) )  \ , f^{Q}(\vec{x}) \right\rangle  \\
&= \left\langle  f^{B}( \sigma_{1}(b_{1},x_{1}),\ldots,\sigma_{1}(b_{n},x_{n}) ) + T_{f}'(\vec{x}) - f^{B}( \sigma_{2}(b_{1},x_{1}),\ldots,\sigma_{2}(b_{n},x_{n}) ) - T_{f}'(\vec{x})  , f^{Q}(\vec{x}) \right\rangle  \\
&= \left\langle f^{B} \left( \sigma_{1}(b_{1},x_{1}) - \sigma_{2}(b_{1},x_{1}), \ldots, \sigma_{1}(b_{n},x_{n}) - \sigma_{2}(b_{n},x_{n}) \right) , f^{Q}(\vec{x}) \right\rangle  \\
&= F_{f} \left( \gamma(b_1,x_{1}), \ldots, \gamma(b_1,x_{1}) \right) ;
\end{align*}
thus, $\gamma$ is a homomorphism such that $\pi = p_{2} \circ \gamma$. If we write $\gamma(c,x) = \left\langle \kappa(c,x), x \right\rangle$, we then have 
\begin{align}\label{eqn:homunivcent1}
\kappa( f^{C}(\vec{b}) + T_{f}(\vec{x}) , f^{Q}(\vec{x}) ) = f^{B}(\kappa(b_{1},x_{1}),\ldots,\kappa(b_{n},x_{n})) .
\end{align}
We note that $\kappa(c,x) - \kappa(0,x)$ is independent of the choice of $x \in Q$; to see this, we apply Eq~\eqref{eqn:homunivcent1} to the difference term to see that
\begin{align*}
\kappa(c,x) &= \kappa \left( m^{C}(b,0,0) + T_{m}(y,y,x) , m^{Q}(y,y,x) \right) \\
&= m^{B} \left( \kappa(b,y), \kappa(0,y), \kappa(0,x) \right) \\
&= \kappa(b,y) - \kappa(0,y) + \kappa(0,x) .
\end{align*}
This yields $\kappa(c,x) - \kappa(0,x) = \kappa(b,y) - \kappa(0,y)$. We also see from Eq~\eqref{eqn:homunivcent1} that
\begin{align}\label{eqn:addcent}
\kappa(b + c,x) = \kappa( m^{C}(b,0,c) + T_{m}(x,x,x), m^{Q}(x,x,x) ) = \kappa(b,x) - \kappa(0,x) + \kappa(c,x) .
\end{align}
Define $\phi : C^{\tau} \rightarrow B^{\tau}$ by $\phi(c) := \kappa(b,x) - \kappa(0,x)$ for any choice of $x \in Q$. Then using Eq~\eqref{eqn:addcent} and Eq~\eqref{eqn:homunivcent1} we see that for $f \in \tau$ with $n =\ar f$,
\begin{align*}
\phi( f^{C}(\vec{b})) &= \kappa (f^{C}(\vec{b}),y) - \kappa(0,y) \\
&= \kappa (f^{C}(\vec{b}),f^{Q}(\vec{x})) - \kappa(0,f^{Q}(\vec{x})) \\
&= \kappa (f^{C}(\vec{b}),f^{Q}(\vec{x})) - \kappa ( T_{f}(\vec{x}) , f^{Q}(\vec{x})) + \kappa ( T_{f}(\vec{x}) , f^{Q}(\vec{x})) - \kappa(0,f^{Q}(\vec{x})) \\
&= \kappa (f^{C}(\vec{b}) + T_{f}(\vec{x}) ,f^{Q}(\vec{x})) - \kappa ( T_{f}(\vec{x}) , f^{Q}(\vec{x})) \\
&= f^{B} \left( \kappa(b_{1},x_{1}), \ldots, \kappa(b_{n},x_{n}) \right) - f^{B} \left( \kappa(0,x_{1}), \ldots, \kappa(0,x_{n}) \right) \\
&= f^{B}\left( \phi(b_{1}), \ldots, \phi(b_{n}) \right) ;
\end{align*}
thus, $\phi$ is a nontrivial homomorphism. Then Eq~\eqref{eqn:homunivcent1} yields 
\begin{align*}
\delta_{T} (\phi)_{f}(\vec{x}) = \phi(T_{f}(\vec{x})) &= \kappa( T_{f}(\vec{x}),f^{Q}(\vec{x}) ) - \kappa(0,f^{Q}(\vec{x})) \\
&= f^{B}\left( \kappa(0,x_{1}),\ldots,\kappa(0,x_{n}) \right) - \kappa(0,f^{Q}(\vec{x}))
\end{align*}
which shows $\delta_{T}(\phi) = 0$. This contradicts injectivity of the transgression. It must be that $\psi_{1}=\psi_{2}$; therefore, we have the unique lifting property for $\pi: A \rightarrow Q$ which shows it is a universal central extension. 
\end{proof}

\begin{proposition}
Let $\mathcal U$ be a variety with a difference term and $Q \in \mathcal U$ with a free-presentation $F_{\mathcal U}(X)/\theta \approx Q$ such that $F_{\mathcal U}(X)/[\theta,1]$ has an idempotent element. A universal central extension of $Q$ in $\mathcal U$ is a cover of $Q$.
\end{proposition}
\begin{proof}
Assume $\pi: A \rightarrow Q$ is a universal central extension in $\mathcal U$. Then by Theorem~\ref{thm:perfectchar}, $A$ is perfect and $H^{2}(Q,B^{\tau}) = 0$ for all kernel algebras $B^{\tau} \in \mathcal U$. We can write $A = E^{\tau} \otimes^{T} Q$. Then by Theorem~\ref{thm:inflationrestriction}, the associated transgression $\delta_{T}$ is always surjective which implies $\pi$ has the $\mathcal U$-lifting property by Proposition~\ref{prop:surjtrans}. Since $Q$ is perfect, $\ker \pi < 1_{A}=[1_{A},1_{A}]$; altogether, $\pi$ is a cover of $Q$.
\end{proof}

\begin{corollary}
Let $\mathcal U$ be a variety with a difference term and $Q \in \mathcal U$ a perfect algebra such that $F_{\mathcal U}(X)/[\theta,1]$ has an idempotent element where $F_{\mathcal U}(X)/\theta \approx Q$ is a free presentation in $\mathcal U$. Assume $H^{2}_{\mathcal V}(Q,E^{\tau})=0$ for some regular datum $(Q,E^{\tau})$ in $\mathcal U$.
\begin{enumerate}

	\item $H^{2}_{\mathcal U}(Q,B^{\tau})=0$ for all kernel algebras $B^{\tau} \in \mathcal U$.
	
	\item If $\gamma \in \Con Q$ is central, then the canonical surjection $\pi : Q \rightarrow Q/\gamma$ is a covering and so $\Hom ( Q(\gamma)/\Delta_{\gamma 1} , E^{\tau}) \approx H^{2}_{\mathcal U}(Q/\gamma,E^{\tau})$.

\end{enumerate}
\end{corollary}
\begin{proof}
(1) By Theorem~\ref{thm:perfectchar}(1), there is a universal central extension $\pi: A \rightarrow Q$ in $\mathcal U$; thus, $A$ is a cover of $Q$ which implies the transgression $\delta_{T}$ is bijective where $T$ is the 2-cocycle associated to the central extension $\pi: A \rightarrow Q$. Then $H^{2}_{\mathcal V}(Q,E^{\tau})=0$ implies $\Hom\left( A(\ker \pi)/\Delta_{\ker \pi 1} , E^{\tau} \right) = 0$; therefore, it must be that $A(\ker \pi)/\Delta_{\ker \pi 1}$ is trivial by the Lemma~\ref{lem:regularnontriv}. Then $A \approx A(\ker \pi)/\Delta_{\ker \pi 1} \otimes^{T} Q \approx 0 \otimes^{T} Q \approx Q$ which implies $Q$ is its own universal central extension. Then $H^{2}_{\mathcal U}(Q,B^{\tau}) = 0$ for all kernel algebras $B^{\tau} \in \mathcal U$ by Theorem~\ref{thm:perfectchar}(2).

(2) Since $Q$ is perfect, we have $[1,1]=1$ in $Q$. Then for the canonical surjection $\pi: Q \rightarrow Q/\gamma$ associated to the central congruence, $\gamma = \ker \pi = \gamma \leq \zeta_{Q} = \zeta_{Q} \wedge [1,1]$. Let $T$ be the 2-cocycle associated to the central extensions $\pi$. By Theorem~\ref{thm:inflationrestriction}, $\im \delta_{T} = \ker \check{\sigma} = H^{2}_{\mathcal U}(Q/\gamma,E^{\tau})$ since $H^{2}_{\mathcal V}(Q,E^{\tau})=0$ by hypothesis. According to Proposition~\ref{prop:surjtrans}(1), the extension $\pi$ has the $\rho$-lifting property for all extensions $\rho: C \rightarrow P$ in $\mathcal U$ realizing regular datum $(P,E^{\tau})$. Altogether, we have shown $\pi$ is a cover of $Q/\gamma$. 
\end{proof}

\vspace{0.2cm}

\section{Invariance of the Schur multiplier}\label{sec:7}

In this section we prove the following theorem concerning the Schur multiplier from Definition~\ref{def:schurmult}.

\begin{theorem}\label{thm:invSchur}
Let $Q \in \mathcal V$ a variety with a difference term. The Schur multiplier of $Q$ in $\mathcal V$ is invariant under presentations.
\end{theorem}

The proof of the theorem depends on Lemma~\ref{lem:mapscent} and Lemma~\ref{lem:commbase} which concerns the recursive generation of the TC-commutator. Fix an algebra $A$ and $\alpha,\beta \in \Con A$. Recall the subalgebra of matrices 
\begin{align}
M(\alpha,\beta) = \left\{ \begin{bmatrix} t(\bar{a},\bar{u}) & t(\bar{a},\bar{v}) \\ t(\bar{b},\bar{u}) & t(\bar{b},\bar{v}) \end{bmatrix} : \bar{a} \mathrel{\alpha} \bar{b}, \bar{u} \mathrel{\beta} \bar{v}, t(\bar{x},\bar{y}) \text{ a term } \right\}.
\end{align}
We will define a countable sequence of linearly ordered congruences whose union equals the commutator.

Set $\tau^{0}(\alpha,\beta) = 0_{A}$ the identity relation. Assume $\tau^{i}(\alpha,\beta)$ is defined and a congruence and let $R^{i+1}(\alpha,\beta)$ be the binary relation defined in the following manner: $(u,v) \in R^{i+1}(\alpha,\beta)$ if there is a matrix $\begin{bmatrix} r & s \\ u & v \end{bmatrix} \in M(\alpha,\beta)$ such that $(u,v) \in \tau^{i}(\alpha,\beta)$. It is easy to see that $R^{i+1}(\alpha,\beta)$ is reflexive, symmetric and closed under unary polynomials. Let $\tau^{i+1}(\alpha,\beta)$ be the transitive closure of $R^{i+1}(\alpha,\beta)$; thus, $\tau^{i+1}(\alpha,\beta)$ is a congruence. The same argument which is applied to the strong and strongly rectangular commutator in \cite[Lem 3.8]{shape} also shows that $\tau^{i}(\alpha,\beta) \leq \tau^{i+1}(\alpha,\beta)$ and 
\begin{align}\label{eqn:705}
[\alpha,\beta] = \bigcup_{i < \omega} \tau^{i}(\alpha,\beta).
\end{align} 
With stronger equational theories, Eq~\eqref{eqn:705} may take more manageable forms. In the case of a Mal'cev variety, we have the following two facts: specializing \cite[Thm 4.9]{commod} from congruence modular varieties
\begin{align}
(a,b) \in [\alpha,\beta] \quad &\Longleftrightarrow \quad \exists u \in A \text{ such that }  \begin{bmatrix} u & u \\ a & b \end{bmatrix} \in \Delta_{\alpha \beta}
\end{align}
and $\Delta_{\alpha \beta} = M(\alpha,\beta)$ \cite[Lemma 7.24]{bergman}. Together these facts show that $[\alpha,\beta] = R^{1}(\alpha,\beta)$. The following lemma can be seen as a weaker analogue which holds in varieties with a difference term.

\begin{lemma}\label{lem:commbase}
Let $A \in \mathcal V$ a variety with a difference term. If $\alpha,\beta \in \Con A$ with $\alpha$ abelian, then $\alpha \wedge [\beta,\beta] = \alpha \wedge R^{1}(\beta,\beta)$.
\end{lemma}
\begin{proof}
The inclusion $\alpha \wedge [\beta,\beta] \supseteq \alpha \wedge R^{1}(\beta,\beta)$ is immediate from definition; for the reverse inclusion, the first step is the reduction in the transitive closure of $R^{m}(\beta,\beta)$.

\begin{claim}\label{claim:transit}
$\alpha \wedge \tau^{m}(\beta,\beta) = \alpha \cap R^{m}(\beta,\beta)$ for each $m$.
\end{claim}
\begin{proof}
We show $\alpha \wedge \tau^{m}(\beta,\beta) \subseteq \alpha \cap R^{m}(\beta,\beta)$ since the reverse inclusion is immediate by definition. Let $(a,b) \in \alpha \wedge \tau^{m}(\beta,\beta)$. Then there exists matrices $\begin{bmatrix} c_{1} & d_{1} \\ a & x_{1}  \end{bmatrix},\begin{bmatrix} c_{2} & d_{2} \\ x_{1} & x_{2}  \end{bmatrix} ,\ldots , \begin{bmatrix} c_{n} & d_{n} \\ x_{n} & b \end{bmatrix}$ in $M(\beta,\beta)$ such that each $(c_{i},d_{i}) \in \tau^{m-1}(\beta,\beta)$. By introducing a constant matrix we can assume $n=2k$ is even. The difference term $m$ satisfies the Mal'cev identities on the pair $(a,b)$ since $\alpha$ is abelian. Using this we can construct a sequence of matrices 
\begin{align*}
\begin{bmatrix} m( c_{1},d_{n} , b) & m(d_{1} ,c_{n} ,b ) \\ a & m(x_{1},x_{n},b ) \end{bmatrix} &= \begin{bmatrix} m( c_{1},d_{n} , b) & m(d_{1} ,c_{n} ,b ) \\ m(a,b,b) & m(x_{1},x_{n},b ) \end{bmatrix}, \begin{bmatrix} m(c_{2} , d_{n},b ) & m(d_{2} ,c_{n},b ) \\ m( x_{1},x_{n} ,b ) & m( x_{2},x_{n-1} ,b ) \end{bmatrix}, \\
&\ldots,\begin{bmatrix} m( c_{k},d_{k+1},b ) & m(d_{k},c_{k+1},b ) \\ m(x_{k-1} , x_{k+1},b ) & m(x_{k},x_{k},b) \end{bmatrix} = \begin{bmatrix} m( c_{k},d_{k+1},b ) & m(d_{k},c_{k+1},b ) \\ m(x_{k-1} , x_{k+1},b ) & b  \end{bmatrix}
\end{align*}
where the top rows of the matrices are still $\tau^{m-1}(\beta,\beta)$-related since $\tau^{m-1}(\beta,\beta)$ is a congruence. Inductively, we can shorten the sequences until there is $\begin{bmatrix} c & d \\ a & b  \end{bmatrix} \in M(\beta,\beta)$ with $(c,d) \in \tau^{m-1}(\beta,\beta)$; thus, $(a,b) \in R^{m}(\beta,\beta)$.
\end{proof}

Suppose we are given $\begin{bmatrix} c & d \\ a & b  \end{bmatrix} = \begin{bmatrix} t(\bar{w},\bar{u}) & t(\bar{w},\bar{v}) \\ t(\bar{z},\bar{u}) & t(\bar{z},\bar{v}) \end{bmatrix} \in M(\beta,\beta)$ with $(a,b) \in \alpha \cap R^{m}(\beta,\beta)$. Define the polynomial $s(\bar{x},\bar{y}) := m(t(\bar{x},\bar{y}) ,t(\bar{x},\bar{u}) ,t(\bar{z},\bar{u}) )$ and observe that 
\begin{align*}
m(d,c,a) = s(\bar{w},\bar{v}) = m(t(\bar{w},\bar{v}),t(\bar{w},\bar{u}),t(\bar{z},\bar{u})) \mathrel{\tau^{m-1}(\beta,\beta)} m( t(\bar{w},\bar{u}) ,t(\bar{w},\bar{u}) ,t(\bar{z},\bar{u})) = a.
\end{align*} 
Then the inclusion of $(a,b) \in \alpha \cap R^{m}(\beta,\beta$ is witnessed by a matrix of the form 
\begin{align}\label{eqn:740}
\begin{bmatrix} a & m(d,c,a) \\ a & b \end{bmatrix} = \begin{bmatrix} s(\bar{w},\bar{u}) & s(\bar{w},\bar{v}) \\ s(\bar{z},\bar{u}) & s(\bar{z},\bar{v}) \end{bmatrix} \in M(\beta,\beta).
\end{align}
Finally, let $(a,b) \in \alpha \wedge [\beta,\beta]$. According to Eq~\eqref{eqn:705}, let $m$ be smallest such that $(a,b) \in \tau^{m}(\beta,\beta)$; for a contradiction, suppose $m > 1$. By Claim~\ref{claim:transit} and Eq~\eqref{eqn:740} there is $\begin{bmatrix} a & y \\ a & b \end{bmatrix} \in M(\beta,\beta)$ with $(a,x) \in \tau^{m-1}(\beta,\beta)$. Then there are matrices $\begin{bmatrix} c_{1} & d_{1} \\ a & x_{1}  \end{bmatrix},\begin{bmatrix} c_{2} & d_{2} \\ x_{1} & x_{2}  \end{bmatrix} ,\ldots , \begin{bmatrix} c_{n} & d_{n} \\ x_{n} & y \end{bmatrix}$ in $M(\beta,\beta)$ with $\tau^{m-2}(\beta,\beta)$. Since $\tau^{0}(\beta,\beta) \leq \tau^{m-2}(\beta,\beta)$ we can adjoin the transpose $\begin{bmatrix} a & a \\ y & b \end{bmatrix}$ to the above sequence of matrices to witness that $(a,b) \in \tau^{m-1}(\beta,\beta)$ which contradicts the assumption $m > 1$. It must be that $\alpha \wedge [\beta,\beta] \subseteq \alpha \wedge \tau^{1}(\beta,\beta) = \alpha \cap R^{1}(\beta,\beta)$ using Claim~\ref{claim:transit}.
\end{proof}

\begin{lemma}\label{lem:mapscent}
Let $Q \in \mathcal V$ a variety with a difference term. Suppose $\pi : F \rightarrow Q$ and $\rho: G \rightarrow Q$ are central extensions in $\mathcal V$ and there exists homomorphisms $\sigma, \lambda : F \rightarrow G$ such that $\rho \circ \lambda = \pi = \rho \circ \sigma$. Then there are homomorphisms $\hat{\sigma}, \hat{\lambda} : F(\ker \pi)/\Delta_{\ker \pi 1} \rightarrow G(\ker \lambda)/\Delta_{\ker \lambda 1}$ such that
\begin{align*}
\sigma(a) &= \left\langle \hat{\sigma} \left( \begin{bmatrix} r(a) \\ a \end{bmatrix}/\Delta_{\ker \pi 1} \right) , \pi(a) \right\rangle ,& \lambda(a) &= \left\langle \hat{\lambda} \left( \begin{bmatrix} r(a) \\ a \end{bmatrix}/\Delta_{\ker \pi 1} \right) , \pi(a) \right\rangle
\end{align*}
and $\hat{\sigma} = \hat{\lambda}$ when restricted to the subalgebra $\kappa_{\ker \pi 1} \big( F(\ker \pi \wedge [1,1]) \big)$. 
\end{lemma}
\begin{proof}

Note that $\rho \circ \lambda = \pi = \rho \circ \sigma$ implies $\sigma(\ker \pi) \subseteq \ker \rho$ and $\lambda(\ker \pi) \subseteq \ker \rho$. Fix a section $l:Q \rightarrow F$ for $\pi$ with associated $\ker \pi$-trace $r = l \circ \pi$. Define $l' := \sigma \circ l : Q \rightarrow G$. Then $\rho \circ l' = \rho \circ \sigma \circ l = \pi \circ l = \id_{Q}$ shows $l'$ is a section for $\rho$ with associated $\ker \rho$-trace $s = l' \circ \rho$. With this choice of section the central extension $\rho : G \rightarrow Q$ is represented by $G \approx G(\ker \rho)/\Delta_{\ker \rho 1} \otimes^{T} Q$ given by $a \longmapsto \left\langle \begin{bmatrix} s(a) \\ a \end{bmatrix}/\Delta_{\ker \rho 1}  , \rho(a) \right\rangle$. Note $s \circ \sigma(a) = l' \circ \rho \circ \sigma(a) = \sigma \circ l \circ \rho \circ \sigma(a) = \sigma \circ l \circ \pi (a) = \sigma \circ r(a)$. If we identify $G$ with its representation, then the image of homomorphism $\sigma$ is given by $\sigma(a) = \left\langle \begin{bmatrix} s (\sigma(a)) \\ \sigma(a) \end{bmatrix}/\Delta_{\ker \rho 1}  , \rho(\sigma(a)) \right\rangle = \left\langle \begin{bmatrix} \sigma(r(a)) \\ \sigma(a) \end{bmatrix}/\Delta_{\ker \rho 1}  , \pi(a) \right\rangle$. The homomorphism $\hat{\sigma}$ is then defined by
\begin{align}\label{eqn:700}
\hat{\sigma} \left( \begin{bmatrix} a \\ b \end{bmatrix}/\Delta_{\ker \pi 1} \right) = \hat{\sigma} \left( \begin{bmatrix} r(a) \\ m(b,a,r(a)) \end{bmatrix}/\Delta_{\ker \pi 1} \right) = \begin{bmatrix} \sigma \circ r(a) \\ \sigma \circ m(b,a,r(a)) \end{bmatrix}/\Delta_{\ker \rho 1}
\end{align}
The homomorphism $\hat{\lambda}$ is defined in the same manner. Note $s \circ \lambda (a) = l' \circ \rho \circ \lambda(a) = \sigma \circ l \circ \pi(a) = \sigma \circ r(a)$. From this we see that 
\begin{align}\label{eqn:702} 
\hat{\lambda} \left( \begin{bmatrix} a \\ b \end{bmatrix}/\Delta_{\ker \pi 1} \right) = \begin{bmatrix} \sigma \circ r(a) \\  m(\lambda(b),\lambda(a), \sigma \circ r(a)) \end{bmatrix}/\Delta_{\ker \rho 1}
\end{align}
and $\sigma(a) \mathrel{\beta} \lambda(a)$ for all $a \in F$.

If we have a matrix $M = \begin{bmatrix} t(\bar{a},\bar{x}) & t(\bar{a},\bar{y}) \\ t(\bar{b},\bar{x}) & t(\bar{b},\bar{y}) \end{bmatrix} \in M(1,1)$ in $F$, we define a new matrix $\sigma(M) = \begin{bmatrix} t(\sigma(\bar{a}),\sigma(\bar{x})) & t(\sigma(\bar{a}),\sigma(\bar{y})) \\ t(\sigma(\bar{b}),\sigma(\bar{x})) & t(\sigma(\bar{b}),\sigma(\bar{y})) \end{bmatrix}$ with $\sigma$ applied entry-wise to $M$. If $(u,v) \in [1,1]$, then by Eq~\eqref{eqn:705} there is some sequence of matrices $M_{1},\ldots,M_{n}$ in $F$ where the entries are related to each other in some complicated manner according to the relations $\tau^{i}(1,1)$; therefore, it follows that the matrices $\sigma(M_{1}),\ldots,\sigma(M_{n})$ witness the inclusion $(\sigma(u),\sigma(v)) \in [1,1]$ in $G$. This implies that $\hat{\sigma}$ and $\hat{\lambda}$ restrict to homomorphisms
\begin{align*}
\hat{\sigma} , \hat{\lambda} : \kappa_{\ker \pi 1} \big( F(\ker \pi \wedge [1,1]) \big) \longrightarrow \kappa_{\ker \rho 1} \big( F(\ker \rho \wedge [1,1]) \big).
\end{align*}

According to Eq~\eqref{eqn:700} and Eq~\eqref{eqn:702}, in order to show $\hat{\sigma} = \hat{\lambda}$ when restricted to the subalgebra $\kappa_{\ker \pi 1} \big( F(\ker \pi \wedge [1,1]) \big)$, we must prove that 
\begin{align}\label{eqn:701}
m (\sigma(v),\sigma(u),\sigma \circ r(u) ) = m (\lambda(v),\lambda(u),\sigma \circ r(u) )
\end{align} 
whenever $(u,v) \in \ker \pi \wedge [1,1]$.

By Lemma~\ref{lem:commbase}, there exists $\begin{bmatrix} t(\bar{a},\bar{x}) & t(\bar{a},\bar{y}) \\ t(\bar{b},\bar{x}) & t(\bar{b},\bar{y}) \end{bmatrix} \in M(1,1)$ in $F$ with $t(\bar{a},\bar{x}) = t(\bar{a},\bar{y})$ and $v=t(\bar{b},\bar{x})$, $u=t(\bar{b},\bar{y})$. Then $t(\lambda(\bar{a}),\lambda(\bar{x})) = t(\lambda(\bar{a}),\lambda(\bar{y}))$ and centrality of $\beta$ yields
\begin{align}\label{eqn:712}
t(\sigma(\bar{a}),\lambda(\bar{x})) = t(\sigma(\bar{a}),\lambda(\bar{y})).
\end{align}
For any $e \in G$, equality in the first and second coordinates  
\begin{align*}
m \Big( t(\sigma(\bar{a}),\lambda(\bar{x})) , t(\sigma(\bar{a}),\lambda(\bar{y})) , e \Big) = m \Big( t(\sigma(\bar{a}),\sigma(\bar{x})) , t(\sigma(\bar{a}),\sigma(\bar{y})) , e \Big)
\end{align*}
yields
\begin{align}\label{eqn:711}
m \Big( t(\sigma(\bar{b}),\lambda(\bar{x})) , t(\sigma(\bar{b}),\lambda(\bar{y})) , e \Big) = m \Big( t(\sigma(\bar{b}),\sigma(\bar{x})) , t(\sigma(\bar{b}),\sigma(\bar{y})) , e \Big)
\end{align}
by centrality of $\beta$. Again the identity of the difference term gives
\begin{align*}
m \Big( t(\sigma(\bar{b}),\lambda(\bar{x})) , t(\sigma(\bar{b}),\lambda(\bar{x})) , e \Big) = m \Big( t(\lambda(\bar{b}),\lambda(\bar{x})) , t(\lambda(\bar{b}),\lambda(\bar{x})) , e \Big)
\end{align*}
which yields
\begin{align}\label{eqn:760}
m \Big( t(\sigma(\bar{b}),\lambda(\bar{x})) , t(\sigma(\bar{b}),\lambda(\bar{y})) , e \Big) = m \Big( t(\lambda(\bar{b}),\lambda(\bar{x})) , t(\lambda(\bar{b}),\lambda(\bar{y})) , e \Big)
\end{align}
by centrality of $\beta$. Then Eq~\eqref{eqn:711} and Eq~\eqref{eqn:760} together yields
\begin{align*}
m \Big( \sigma(v), \sigma(u) , e \Big) &= m \Big( t(\sigma(\bar{b}),\sigma(\bar{x})) , t(\sigma(\bar{b}),\sigma(\bar{y})) , e \Big) \\
&= m \Big( t(\lambda(\bar{b}),\lambda(\bar{x})) , t(\lambda(\bar{b}),\lambda(\bar{y})) , e \Big) \\
&= m \Big( \lambda(v), \lambda(u) , e \Big) .
\end{align*}
By taking $e = \sigma \circ r(u)$ we have Eq~\eqref{eqn:701}.
\end{proof}

\begin{proof}[(proof of Theorem~\ref{thm:invSchur})] Let $F/\alpha \approx Q \approx G/\beta$ be free-presentations in $\mathcal V$. By Lemma~\ref{lem:100}(1) there are central extensions $\pi : F/[\alpha,1] \rightarrow Q$, $\lambda : G/[\beta,1] \rightarrow Q$ and homomorphisms $\sigma: F/[\alpha,1] \rightarrow G/[\beta,1], \lambda: G/[\beta,1] \rightarrow F/[\alpha,1]$ such that $\rho \circ \sigma = \pi$ and $\rho = \pi \circ \lambda$. From the definition of the induced maps we see that $\widehat{\lambda \circ \sigma} = \hat{\lambda} \circ \hat{\sigma}$ and $\hat{\id} = \id$ for the induced map of the identity $\id:F/[\alpha,1] \rightarrow F/[\alpha,1]$. Then Lemma~\ref{lem:mapscent} applied to the central extension $\pi: F/[\alpha,1] \rightarrow Q$ and the two homomorphisms $\id , \lambda \circ \sigma : F/[\alpha,1] \rightarrow F/[\alpha,1]$ yields $\id = \hat{\lambda} \circ \hat{\sigma}$; similarly, the central extension $\rho: G/[\beta,1] \rightarrow Q$ and the homomorphisms $\id , \sigma \circ \lambda : G/[\beta,1] \rightarrow G/[\beta,1]$ yields $\id = \hat{\sigma} \circ \hat{\lambda}$.
\end{proof}

\begin{acknowledgments}
The research in this manuscript was supported by NSF China Grant \#12071374.  I would also like to express my sincere gratitude to Professor Shiqing Zhang and the mathematics department at Sichuan University for a very rewarding and productive research term. 
\end{acknowledgments}


\end{document}